\documentclass{article}

\usepackage{arxiv}

\usepackage[utf8]{inputenc} 
\usepackage[T1]{fontenc}    
\usepackage{url}            
\usepackage{booktabs}       

\usepackage{amsmath,amssymb,amsfonts,latexsym,amsthm,bm,breqn,dsfont,gensymb}
\usepackage{latexsym,breqn}
\usepackage{nicefrac}       
\usepackage{microtype}      
\usepackage{natbib}

\usepackage{graphicx,subfigure}
\usepackage{tikz}
\definecolor{greenmoi}{RGB}{0,100,0}

\usepackage{algorithm}
\usepackage[noend]{algpseudocode}
\makeatletter
\def\BState{\State\hskip-\ALG@thistlm}

\usepackage{ctable,multicol,multirow}

\usepackage{color}
\definecolor{green}{rgb}{0, 0.5, 0}
\usepackage[pdftex,
			bookmarks, bookmarksnumbered=true, bookmarksopen=true,
			colorlinks, citecolor=blue, filecolor=black, linkcolor=red, urlcolor=green,
			pdfauthor={López-Lopera, Andrés Felipe}, pdftitle={PhDEMSE},
			pdftoolbar=true, pdfmenubar=true, pdffitwindow = true
			]{hyperref}

\usepackage{pifont}
\definecolor{green}{rgb}{0, 0.5, 0}

\newcommand{\ensnombre}[1]{\mathbb{#1}}
\newcommand{\R}{\ensnombre{R}}
\newcommand{\beq} {\begin{eqnarray*}}
	\newcommand{\eeq} {\end{eqnarray*}}
\newcommand{\N}{\ensnombre{N}}
\newcommand{\tbf} {\textbf}
\newcommand{\norme}[1]{\left\lVert #1 \right\rVert}

\DeclareMathOperator{\argmax}{argmax}

\DeclareMathOperator{\Var}{Var}
\DeclareMathOperator{\Cov}{Cov}

\newcommand{\noi} {\noindent}

\def \E{\mathbb{E}}
\def \P{\mathbb{P}}

\newcommand{\defeq}{\mathrel{\mathop:}=}

\newcommand{\abs}[1]{\left\lvert #1 \right\rvert}
\newcommand{\eqdef}{\mathrel{=}:}

\newtheorem{theorem}{Theorem}[section]
\newtheorem{corollary}[theorem]{Corollary}
\newtheorem{lemma}[theorem]{Lemma}
\newtheorem{proposition}[theorem]{Proposition}

\title{Maximum likelihood estimation for Gaussian processes under inequality constraints}
\date{}
\author{
  Fran\c{c}ois Bachoc \\
  Institut de Math\'ematiques de Toulouse\\
	Universit\'e Paul Sabatier  \\  
  31062 Toulouse, France \\
  \texttt{Francois.Bachoc@math.univ-toulouse.fr} \\  
  \And
  Agn\`es Lagnoux\\
  Institut de Math\'ematiques de Toulouse\\
	Universit\'e Toulouse Jean Jaur\`es\\  
  31062 Toulouse, France \\	
  \texttt{lagnoux@univ-tlse2.fr} \\
  \And  
  Andr\'es F. L\'opez-Lopera\\
  Mines Saint-\'Etienne\\
  42000 Saint-\'Etienne, France \\
  \texttt{andres-felipe.lopez@emse.fr} \\  
}

\begin{document}
\maketitle

\begin{abstract}
	We consider covariance parameter estimation for a Gaussian process under inequality constraints (boundedness, monotonicity or convexity) in fixed-domain asymptotics. We address the estimation of the variance parameter and the estimation of the microergodic parameter of the Mat\'ern and Wendland covariance functions. First, we show that the (unconstrained) maximum likelihood estimator has the same asymptotic distribution, unconditionally and conditionally to the fact that the Gaussian process satisfies the inequality constraints. Then, we study the recently suggested constrained maximum likelihood estimator. We show that it has the same asymptotic distribution as the (unconstrained) maximum likelihood estimator. In addition, we show in simulations that the constrained maximum likelihood estimator is generally more accurate on finite samples. Finally, we provide extensions to prediction and to noisy observations.
\end{abstract}

\section{Introduction}\label{sec:intro}

Kriging \citep{stein99interpolation,rasmussen06gaussian} consists in inferring the values of a Gaussian random field given observations at a finite set of points.
It has become a popular method for a
large range of applications, such as geostatistics \citep{matheron70theorie}, numerical code approximation \citep{sacks89design,santner03design,bachoc16improvement} and calibration \citep{paulo12calibration,bachoc14calibration}, global optimization \citep{jones98efficient}, and machine learning \citep{rasmussen06gaussian}.

When considering a Gaussian process, one has to deal with the estimation of its covariance function.
Usually, it is assumed that the covariance function belongs to a given parametric family \citep[see][for a review of classical families]{abrahamsen97review}.
In this case, the estimation boils down to estimating the corresponding covariance parameters. The main estimation techniques are based on maximum likelihood \citep{stein99interpolation}, cross-validation \citep{zhang10kriging,Bachoc13cross,Bac2014} and variation estimators \citep{istas97quadratic,And2010,ABKLN18}. 

In this paper, we address maximum likelihood estimation of covariance parameters under fixed-domain asymptotics \citep{stein99interpolation}. The fixed-domain asymptotics setting corresponds to observation points for the Gaussian process that become dense in a fixed bounded domain. 
Under fixed-domain asymptotics, two types of covariance parameters can be distinguished: microergodic and non-microergodic parameters \citep{ibragimov78gaussian,stein99interpolation}. A covariance parameter is said to be microergodic if, when it takes two different values, the two corresponding Gaussian measures are orthogonal \citep{ibragimov78gaussian,stein99interpolation}. It is said to be non-microergodic if, even for two different values, the corresponding Gaussian measures are equivalent. Although non-microergodic parameters cannot be estimated consistently, they have an asymptotically negligible impact on prediction \citep{AEPRFMCF,BELPUICF,UAOLPRFUISOS,zhang04inconsistent}. On the contrary, it is at least possible to consistently estimate microergodic covariance parameters, and misspecifying them can have a strong negative impact on predictions.

It is still challenging to obtain results on maximum likelihood estimation of microergodic parameters that would hold for very general classes of covariance functions. Nevertheless, significant contributions have been made for specific types of covariance functions. In particular, when considering the isotropic Mat\'ern family of covariance functions, for input space dimension $d=1,2,3$, a reparameterized quantity obtained from the variance and correlation length parameters is microergodic \citep{zhang04inconsistent}. It has been shown in \citep{ShaKau2013}, from previous results in \citep{DuZhaMan2009} and \citep{WanLoh2011}, that the maximum likelihood estimator of this microergodic parameter is consistent and asymptotically Gaussian distributed. Anterior results on the exponential covariance function have been also obtained in \citep{ying91asymptotic,ying93maximum}. 

In this paper, we shall consider the situation where the trajectories of the Gaussian process are known to satisfy  either boundedness, monotonicity or convexity constraints. Indeed, Gaussian processes with inequality constraints provide suitable regression models in application fields such as computer networking (monotonicity) \citep{Golchi2015MonotoneEmulation}, social system analysis (monotonicity) \citep{Riihimaki2010GPwithMonotonicity} and econometrics (monotonicity or positivity) \citep{Cousin2016KrigingFinancial}. Furthermore, it has been shown that taking the constraints into account may considerably improve the predictions and the predictive intervals for the Gaussian process \citep{DaVeiga2012GPineqconst,Golchi2015MonotoneEmulation,Riihimaki2010GPwithMonotonicity}.

Recently, a constrained maximum likelihood estimator (cMLE) for the covariance parameters has been suggested in \citep{LLBDR17}. Contrary, to the (unconstrained) maximum likelihood estimator (MLE) discussed above, the cMLE explicitly takes into account the additional information brought by the inequality constraints. In \citep{LLBDR17}, it is shown, essentially, that the consistency of the MLE implies the consistency of the cMLE under boundedness, monotonicity or convexity constraints. 

The aim of this paper is to study the asymptotic conditional distributions of the MLE and the cMLE, given that the Gaussian process satisfies the constraints. We consider the estimation of a single variance parameter and the estimation of the microergodic parameter in the isotropic Mat\'ern family of covariance functions. In both cases, we show that the asymptotic conditional distributions of the MLE and the cMLE are identical to the unconditional asymptotic distribution of the MLE. Hence, it turns out that the impact of the constraints on covariance parameter estimation is asymptotically negligible.
To the best of our knowledge, this paper is the first work on the asymptotic distribution of covariance parameter estimators for constrained Gaussian processes. The proofs involve tools from asymptotic spatial statistics, extrema of Gaussian processes and reproducing kernel Hilbert spaces. These proofs bring a significant level of novelty compared to these in \citep{LLBDR17}, where only consistency is addressed.
In simulations, we confirm that for large sample sizes, the MLE and the cMLE have very similar empirical distributions, that are close to the asymptotic Gaussian distribution. For small or moderate sample sizes, we observe that the cMLE is generally more accurate than the MLE, so that taking the constraints into account is beneficial. 
Finally, we explore three extensions: to prediction, to the Wendland covariance model and to the framework of noisy observations.

The rest of the manuscript is organized as follows. In Section \ref{sec:framework}, we introduce in details the constraints, the MLE, and the cMLE. In Section \ref{sec:var} we provide the asymptotic results for the estimation of the variance parameter, while the asymptotic results for the isotropic Mat\'ern family of covariance functions are given in Section \ref{sec:matern}. In Section \ref{sec:num}, we report the simulation outcomes. The extensions are presented in Section \ref{sec:extensions}. Concluding remarks are given in Section \ref{sec:ccl}. 
All the proofs are postponed to the appendix.

\section{Gaussian processes under inequality constraints}\label{sec:framework}

\subsection{Framework and purpose of the paper}\label{ssec:framework}

We consider a parametric set of functions $ \{ k_{\theta}; \theta \in \Theta \}$ defined from $\R^d$ to $\R$, where $\Theta $ is a compact set of $\mathbb{R}^p$. We also assume that, for each $\theta \in \Theta$, there exists a Gaussian process with continuous realizations having mean function zero and covariance function $\widetilde k_{\theta}$ on $[0,1]^d\times [0,1]^d$ defined by $\widetilde k_{\theta}(u,v)=k_{\theta}(u-v)$ for $u, v\in [0,1]^d$. We refer to, e.g., \citep{adler1990introduction} for mild smoothness conditions on $k_{\theta}$ ensuring this.
We consider an application 
\begin{align*}
Y : ( \Omega , \mathcal{A} )  \to (\mathcal{C}( [0,1]^d , \R), \mathcal{B} ),
\end{align*}
where $( \Omega , \mathcal{A} )$ is a measurable space, $\mathcal{C}( [0,1]^d , \R )$ is the set of continuous functions from $[0,1]^d$ to $\mathbb{R}$, and $\mathcal{B}$ is the Borel Sigma algebra on $\mathcal{C}( [0,1]^d  , \R)$ corresponding to the $L^{\infty}$ norm. For each $\theta \in \Theta$, let $\mathbb{P}_{\theta}$ be a probability measure on $\Omega$ for which
\begin{align*}
Y : ( \Omega , \mathcal{A} , \mathbb{P}_{\theta} )  \to (\mathcal{C}( [0,1]^d  , \R), \mathcal{B} ),
\end{align*}
has the distribution of a Gaussian process with mean function zero and covariance function $\widetilde k_{\theta}$. 

Now consider a triangular array $\big( x^{(n)}_i \big)_{n \in \mathbb{N},i=1,\ldots,n}$ of observation points in $[0,1]^d$, where we write for concision
$(x_1,\ldots,x_n) = \big(x^{(n)}_1,\ldots,x^{(n)}_n\big)$.
We assume that $\big( x^{(n)}_i \big)$ is dense, that is $\sup_{x \in [0,1]^d} \inf_{i=1, \ldots ,n} |x - x_i^{(n)}| \to 0$ as $n \to \infty$.
Let $y$ be the Gaussian vector defined by $y_i = Y(x_i)$ for $i=1,\ldots,n$. For $\theta \in \Theta$, let $R_{\theta} = [ k_{\theta}(x_i - x_j) ]_{1 \leqslant i,j \leqslant n }$ and
\begin{equation} \label{eq:log:lik}
\mathcal{L}_n(\theta)
=
- \frac{n}{2} \ln( 2 \pi )
- \frac{1}{2} \ln( | R_{\theta} | )
- \frac{1}{2} y^{\top} R_{\theta}^{-1} y,
\end{equation}
be the log likelihood function.  Here, $| R_{\theta} |$ stands for $ \det(R_{\theta}) $. Maximizing $\mathcal{L}_n(\theta)$ with respect to $\theta$ yields the widely studied and applied MLE \citep{santner03design,stein99interpolation,ying93maximum,zhang04inconsistent}.

In this paper, we assume that the information $ \{Y \in \mathcal E_{\kappa} \}$ is available where $\mathcal E_{\kappa }$ is a convex set of functions defined by inequality constraints. We will consider
\begin{align*}
\begin{array}{lll}
\mathcal E_{0}=\{ f\in \mathcal C([0,1]^d,\R) \quad \textrm{s.t.}   \; \ell \leqslant f(x) \leqslant u,\, \forall x\in [0,1]^d\},\\
\mathcal E_{1}=\{ f\in \mathcal C^1([0,1]^d,\R) \quad \textrm{s.t.} \;  \partial f(x) / \partial x_i \geqslant 0, \; \forall x\in [0,1]^d, i \in \{ 1,\ldots,d\}  \},\\
\mathcal E_{2}=\{ f\in \mathcal C^2([0,1]^d,\R) \quad \textrm{s.t. $f$ is convex}\},\\
\end{array}
\end{align*}
which correspond to boundedness, monotonicity and convexity constraints respectively 
. For $\mathcal E_{0}$, the bounds $- \infty \leqslant \ell < u \leqslant + \infty$ are fixed and known. 

First, we will study the conditional asymptotic distribution of the (unconstrained) MLE obtained by maximizing \eqref{eq:log:lik}, given $ \{Y \in \mathcal E_{\kappa} \}$. Nevertheless, a drawback of this MLE is that it does not exploit the information $ \{Y \in \mathcal E_{\kappa} \}$. Then we study the cMLE introduced in \citep{LLBDR17}. This estimator is obtained by maximizing the logarithm of the probability density function of $y$, conditionally to $\{Y \in \mathcal E_{\kappa} \}$, with respect to the probability measure $\mathbb{P}_{\theta}$ on $\Omega$. This logarithm of conditional density is given by
\begin{align} \label{eq:const:log:lik}
\mathcal{L}_{n,c}(\theta) 
= 
\mathcal{L}_n(\theta) - \ln \mathbb{P}_{\theta} ( Y \in \mathcal E_\kappa)
+
\ln  \mathbb{P}_{\theta} ( Y \in \mathcal E_\kappa | y ) 
= 
\mathcal{L}_n(\theta) 
+A_n(\theta)
+ B_n(\theta),
\end{align}
say, where $\P_{\theta}(\cdot{})$ and $\P_{\theta}(\cdot{}|\cdot{})$ are defined in Section \ref{ssec:nota}. 
In \citep{LLBDR17}, the cMLE  is studied and compared to the MLE. The authors show that the cMLE is consistent when the MLE is. In this paper, we aim at providing more quantitative results regarding the asymptotic distribution of the MLE and the cMLE, conditionally to $\{Y \in \mathcal E_{\kappa} \}$.

\subsection{Notation}\label{ssec:nota}

In the paper, $0 < c < + \infty$ stands for a generic constant that may differ from one line to another. It is convenient to have short expressions for terms that converge in probability to zero. Following \citep{van2000asymptotic}, the notation $o_{\P}(1)$ (respectively $O_{\P}(1)$) stands for a sequence of random variables (r.v.'s) that converges to zero in probability (resp. is bounded in probability) as $n \to \infty$. More generally, for a sequence of r.v.'s $R_n$,
\beq
X_n&=&o_{\P}(R_n) \quad \textrm{means} \quad X_n=Y_nR_n \quad \textrm{with} \quad Y_n\overset{\P}{\rightarrow}0,\\
X_n&=&O_{\P}(R_n) \quad \textrm{means} \quad X_n=Y_nR_n \quad \textrm{with} \quad Y_n=O_{\P}(1).
\eeq 
For deterministic sequences $X_n$ and $R_n$, the stochastic notation reduce to the usual $o$ and $O$.
For a sequence of random vectors or variables $(X_n)_{n \in \mathbb{N}}$ on $\mathbb{R}^l$, that are functions of $Y$, and for a probability distribution $\mu$ on $\mathbb{R}^l$, we write 
\[
X_n \xrightarrow[n \to \infty]{\mathcal{L}| Y \in \mathcal E_{\kappa}}
\mu,
\] 
when, for any bounded continuous function $g: \mathbb{R}^l \to \mathbb{R}$, we have
\[
\mathbb{E} 
\left[  \left. g(X_n)  \right| Y \in \mathcal E_{\kappa}  \right]
\underset{n\to\infty}{\longrightarrow}
\int_{\mathbb{R}^l} g(x) \mu(dx).
\]
We also write $X_n = o_{ \P | Y \in \mathcal E_{\kappa} }(1)$
when for all $\varepsilon >0$ we have $\P( |X_n| \geqslant \varepsilon | Y \in \mathcal E _{\kappa} ) \to 0$ as $n \to \infty$. Finally, we write $X_n = O_{ \P | Y \in \mathcal E_{\kappa} }(1)$
when we have $ \limsup_{n \to \infty} \P( |X_n| \geqslant K | Y \in \mathcal E _{\kappa} ) \to 0$ as $K \to \infty$.

For any two functions $f(Y)$ and $g(Y)$, let $\mathbb{E}_{\theta}[f(Y)]$  (respectively $\mathbb{E}_{\theta}[f(Y)|g(Y)]$) be the expectation  (resp. the conditional expectation) with respect to the measure $\mathbb{P}_{\theta}$ on $\Omega$. We define similarly $\mathbb{P}_{\theta}(A(Y))$ and $\mathbb{P}_{\theta}(A(Y)|g(Y))$ when $A(Y)$ is an event with respect to $Y$. Let $\theta_0 \in \Theta$ be fixed. We consider $\theta_0$
as the true unknown covariance parameter and we let $\mathbb{E}[\cdot]$, $\mathbb{E}[\cdot|\cdot]$, $\mathbb{P}(\cdot)$, and $\mathbb{P}(\cdot|\cdot)$ be shorthands for $\mathbb{E}_{\theta_0}[\cdot]$, $\mathbb{E}_{\theta_0}[\cdot|\cdot]$, $\mathbb{P}_{\theta_0}(\cdot)$, and $\mathbb{P}_{\theta_0}(\cdot|\cdot)$. When a quantity is said to converge, say, in probability or almost surely, it is also implicit that we consider the measure $\mathbb{P}_{\theta_0}$ on $\Omega$. 

\subsection{Conditions on the observation points}\label{ssec:doe}

In some cases, we will need to assume that as $n \to \infty$, the triangular array of observation points contains finer and finer tensorized grids. \\

\tbf{Condition-Grid}. There exist $d$ sequences $\big(v_i^{(j)}\big)_{i\in \N}$ for $j=1,\ldots,d$, dense in $[0,1]$, and so that for all $N \in \mathbb{N}$, there exists $n_0 \in \mathbb{N}$ such that for $n \geqslant n_0$, we have
$ \big\{(v_{i_1}^{(1)},\ldots,v_{i_d}^{(d)}),\; 1\leqslant i_1,\ldots,i_d\leqslant N \big\}\subset 
(x_i)_{i =1,\ldots, n}$. \\
\; \\
In our opinion, Condition-Grid is reasonable and natural. Its purpose is to guarantee that the partial derivatives of $Y$ are consistently estimable from $y$ everywhere on $[0,1]^d$ (see, for instance, the proof of Theorems \ref{th:cond_after_var} and \ref{th:cond_before_var} for $\kappa = 1$ in the appendix).
We believe that, for the results for which Condition-Grid is assumed, one could replace it by a milder condition and prove similar results. Then the proofs would be based on essentially the same ideas as the current ones, but could be more cumbersome.

In some other cases, we only need to assume that the observation points constitute a sequence.\\

\tbf{Condition-Sequence}. For all $n \in \mathbb{N}$ and $i \leqslant n$, we have $x^{(n)}_i = x^{(i)}_i$.
\\
\; \\
Condition-Sequence implies that sequences of conditional expectations with respect to the observations are martingales. This condition is necessary in some of the proofs (for instance, that of Theorem \ref{th:cond_before_var}) where convergence results for martingales are used.

\section{Variance parameter estimation}\label{sec:var}

\subsection{Model and assumptions}\label{ssec:model_var}

In this section, we focus on the estimation of a single variance parameter when the correlation function is known. Hence, we let $p=1$, $\theta = \sigma^2$, and for $x \in \R^d$,
\begin{align}\label{def:var}
k_{\sigma^2}(x)
=\sigma^2 k_1(x),
\end{align}
where $k_1$ is a fixed known function such that $\widetilde k_1$ defined by $\widetilde k_1(u,v)=k_1(u-v)$ is a correlation function on $[0,1]^d\times [0,1]^d$. \\

We define the Fourier transform of a function $h\colon \R^d \to \R$ by 
\[
\widehat h(\omega)=\frac{1}{(2\pi)^d} \int_{\R^d} h(t)e^{-\i \omega ^\top t}dt,
\]
where $\i^2=-1$ and we make the following assumption. \\

\tbf{Condition-Var}. Let  $\kappa$ be fixed in $\{0,1,2\}$.
\begin{itemize}
	\item[-] If $\kappa=0$, $k_1$ is $\alpha$-H\"older, which means that there exist non-negative constants $c$ and $\alpha$ such that
	\[
	|k_1(t)-k_1(t')|\leqslant c\norme{t-t'}^{\alpha },
	\]
	for all $t$ and $t'$ in $\R^d$, where $\norme{.}$ is the Euclidean norm.  Furthermore, the Fourier transform $\widehat k_{1}$ of $k_{1}$ satisfies, for some fixed $P<\infty$,
	\begin{align}\label{eq:fourier_transform}
	\widehat k_{1}(\omega) \norme{\omega}^P \underset{\norme{\omega} \to \infty}{\longrightarrow} \infty.
	\end{align}
	
	\item[-] If $\kappa=1$, the Gaussian process $Y$ is differentiable in quadratic mean. For $i=1,\ldots, d$, let $k_{1,i}= - \partial^2 k_1/ \partial x_i^2$. Remark that the covariance function of $\partial Y/\partial x_i$ is given by $\widetilde k_{1,i}$ defined by $\widetilde k_{1,i}(u,v)=k_{1,i}(u-v)$. Then $k_{1,i}$ is $\alpha$-H\"older for a fixed $\alpha >0$. Also, \eqref{eq:fourier_transform} holds with $\widehat k_1$ replaced by the Fourier transform $\widehat k_{1,i}$ of $k_{1,i}$ for $i=1,\ldots,d$.
	\item[-] If $\kappa=2$, the Gaussian process $Y$ is twice differentiable in quadratic mean. For $i,j=1,\ldots, d$, let $k_{2,i,j}=\partial^4 k_1/( \partial x_i^2\partial x_j^2)$. Remark that the covariance function of $\partial^2 Y/(\partial x_i\partial x_j)$ is given by $\widetilde k_{2,i,j}$ defined by $\widetilde k_{2,i,j}(u,v)=k_{2,i,j}(u-v)$.
	Then $k_{2,i,j}$ is $\alpha$-H\"older for a fixed $\alpha >0$. Also, \eqref{eq:fourier_transform} holds with $\widehat k_1$ replaced by the Fourier transform $\widehat k_{2,i,j}$ of $k_{2,i,j}$ for $i,j=1,\ldots,d$.
\end{itemize}

These assumptions make the conditioning by $\{Y\in \mathcal E_{\kappa}\}$ valid for $\kappa=0,1,2$ as established in the following lemma.

\begin{lemma} \label{lem:lower:bounded:proba:bounded:GP_var}
	Assume that Condition-Var holds.   Then for all $\kappa \in \{0,1,2\}$ and for any compact $K$ in $(0,+\infty)$, we have
	\[
	\underset{ \sigma^2 \in K }{\inf}
	\P_{\sigma^2 } 
	\left( Y\in \mathcal E_{\kappa} \right)
	>0.
	\]
\end{lemma}

\tbf{Proof of Lemma \ref{lem:lower:bounded:proba:bounded:GP_var}}. It suffices to follow the same lines as in the proof of  \citep[Lemma A.6]{LLBDR17} noticing that Condition-Var implies the conditions of \citep[Lemma A.6]{LLBDR17} \citep[see the discussion in][]{LLBDR17}.
\hfill $\square$

\subsection{Asymptotic conditional distribution of the maximum likelihood estimator}\label{ssec:after_var}

The log-likelihood function in \eqref{eq:log:lik} for $\sigma^2$  can be written as
\begin{align}\label{def:log_likeli_var}
\mathcal L_n(\sigma^2)=- \frac{n}{2} \ln( 2 \pi )-\frac n2  \ln (\sigma^2)-\frac 12\ln (|R_1|)-\frac{1}{2\sigma^2} y^{\top} R_1^{-1} y,
\end{align}
where $R_1 = [ k_1 (  x_i- x_j )   ]_{1 \leqslant i,j \leqslant n}$. Then the standard MLE is given by
\begin{align}\label{def:MLE_var}
\bar{\sigma}_{n}^2 \in \underset{\sigma^2>0}{\argmax \; } \mathcal L_{n}(\sigma^2).
\end{align}

Now we show that, for $\kappa=0,1,2$, $\sqrt n \left(\bar{\sigma}_n^2-\sigma_0^2 \right)$ is asymptotically Gaussian distributed conditionally to $\{Y\in \mathcal E_{\kappa}\}$.

\begin{theorem}\label{th:cond_after_var} For $\kappa=1,2$, we assume that Condition-Grid holds. For $\kappa=0,1,2$, under Condition-Var, the MLE $\bar{\sigma}_n^2$ of $\sigma_0^2$ defined by \eqref{def:MLE_var} conditioned on $\{ Y\in \mathcal E_{\kappa} \}$ is asymptotically Gaussian distributed. More precisely, 
	\[
	\sqrt n
	\left(\bar{\sigma}_n^2-\sigma_0^2 \right)
	\xrightarrow[n\to+\infty]{\mathcal{L} | Y\in \mathcal E_{\kappa}} \mathcal N(0,2\sigma_0^4).
	\]
\end{theorem}

It is well known that $\sqrt n \left(\bar{\sigma}_n^2-\sigma_0^2 \right)$ converges (unconditionally) to the $\mathcal N(0,2\sigma_0^4)$ distribution. Hence, conditioning by $\{ Y\in \mathcal E_{\kappa} \}$ has no impact on the asymptotic distribution of the MLE.

\subsection{Asymptotic conditional distribution of the constrained maximum likelihood estimator}\label{ssec:before_var}

Here, we assume that the compact set  $\Theta$ is $[\sigma_l^2,\sigma_u^2]$ with $0< \sigma_l^2 < \sigma_0^2 < \sigma_u^2<+\infty$ and we consider  the cMLE $\widehat \sigma_{n,c}^2$ of $\sigma_0^2$ derived by maximizing on the compact set $\Theta$ the constrained log-likelihood in \eqref{eq:const:log:lik}:
\begin{align}\label{def:MLE_cond}
\widehat \sigma_{n,c}^2 \in \underset{\sigma^2 \in \Theta} {\argmax \; }\mathcal L_{n,c}(\sigma^2).
\end{align}

Now we show that the conditional asymptotic distribution of the cMLE is the same as the asymptotic distribution of the MLE.

\begin{theorem}\label{th:cond_before_var}
	For $\kappa=1,2$, we assume that Condition-Grid holds. For $\kappa=0,1,2$, under Condition-Var and Condition-Sequence, the cMLE $\widehat \sigma_{n,c}^2$ of $\sigma_0^2$ defined in \eqref{def:MLE_cond} is asymptotically Gaussian distributed. More precisely,
	\[
	\sqrt n \left(\widehat \sigma_{n,c}^2-\sigma_0^2 \right) \xrightarrow[n\to+\infty]{\mathcal{L} |  Y\in \mathcal E_{\kappa}} \mathcal N(0,2\sigma_0^4).
	\]
\end{theorem}

\section{Microergodic parameter estimation for the isotropic Mat\'ern model}\label{sec:matern}

\subsection{Model and assumptions}\label{ssec:model_matern}

In this section, we let $d=1,2$ or $3$ and we consider the isotropic Mat\'ern family of covariance functions on $\R^d$. We refer to, e.g., \citep{stein99interpolation} for more details. 
Here $k_{\theta} = k_{\theta,\nu}$ is given by, for $x \in [0,1]^d$,
\[
k_{\theta,\nu}(x)
=\sigma^2 K_{\nu}\left(\frac{\norme{x}}{\rho}\right)=\frac{\sigma^2}{\Gamma(\nu)2^{\nu-1}}\left(\frac{\norme{x}}{\rho}\right)^{\nu}\kappa_{\nu}\left(\frac{\norme{x}}{\rho}\right). 
\]
The Mat\'ern covariance function is given by $\widetilde k_{\theta,\nu}(u,v)=k_{\theta,\nu}(u-v)$. 
The parameter
$\sigma^2>0$ is the variance of the process, $\rho>0$ is the  correlation length parameter that controls how fast the covariance function decays with the distance, and $\nu>0$ is the regularity parameter of the process. The function $\kappa_{\nu}$ is the modified Bessel function of the second kind of order $\nu$ \citep[see][]{AS64}. We assume in the sequel that the smoothness parameter $\nu$ is known. Then $\theta = (\sigma^2 , \rho)$ and $p=2$. \\

\tbf{Condition-$\nu$}. For $\kappa=0$ (respectively $\kappa=1$ and  $\kappa=2$), we assume that $\nu>0$ (resp. $\nu>1$ and $\nu>2$).\\

We remark that Condition-$\nu$ naturally implies Condition-Var so that the conditioning by $\{Y\in \mathcal E_{\kappa}\}$ is valid for any $\kappa=0,1,2$ as established in the next lemma. We refer to \citep{stein99interpolation} for a reference on the impact of $\nu$ on the smoothness of the Mat\'ern function $k_{\theta,\nu}$ and on its Fourier transform.

\begin{lemma} \label{lem:lower:bounded:proba:bounded:matern}
	Assume that Condition-$\nu$ holds. Then for all $\kappa\in\{0,1,2\}$ and for any compact $K$ of $(0,\infty)^2$, we have 
	\[
	\underset{ (\sigma^2 , \rho) \in K}{\inf}
	\P_{\sigma^2 , \rho} 
	\left(
	Y\in \mathcal E_{\kappa}
	\right)
	>0.
	\]
\end{lemma}

\tbf{Proof of Lemma \ref{lem:lower:bounded:proba:bounded:matern}}. This lemma is a special case of  \citep[Lemma A.6]{LLBDR17}.  
\hfill $\square$\\

\subsection{Asymptotic conditional distribution of the maximum likelihood estimator}\label{ssec:after_matern}

The log-likelihood function in \eqref{eq:log:lik} for $\sigma^2$ and $\rho$ under the Mat\'ern model with fixed parameter $\nu$ can be written as
\begin{align}\label{def:log_likeli_matern}
\mathcal L_n(\sigma^2,\rho)=-\frac n2 \ln (2\pi)-\frac n2  \ln (\sigma^2)-\frac 12\ln (|R_{\rho,\nu}|)-\frac{1}{2\sigma^2} y^{\top} R_{\rho,\nu}^{-1} y,
\end{align}
where $R_{\rho,\nu} = [ K_{\nu} ( \norme{ x_i - x_j } / \rho ) ]_{1 \leqslant i,j \leqslant n}$. Let $\Theta = [\sigma_l^2 , \sigma_u^2] \times [ \rho_l , \rho_u ]$ with fixed $0 < \sigma_l^2 < \sigma_u^2 < \infty$ and fixed $0 < \rho_l < \rho_u < \infty$. Moreover, assume that the true parameters are such that  $\sigma_l^2/(\rho_l^{2 \nu}) < \sigma_0^2/(\rho_0^{2 \nu}) < \sigma_u^2/(\rho_u^{2 \nu})$. Then the MLE is given by
\begin{align}\label{def:MLE_matern}
(\widehat \sigma_{n}^2, \widehat \rho_n) \in \underset{(\sigma^2,\rho) \in \Theta}{\argmax \; } \mathcal L_{n}(\sigma^2,\rho).
\end{align}

It has been shown in \citep{zhang04inconsistent} that the parameters $\sigma_0^2$ and $\rho_0$ can not be estimated consistently but that the microergodic parameter $\sigma_0^2 / \rho_0^{2 \nu}$ can. Furthermore, it is shown in \citep{ShaKau2013} that $\sqrt n
\left( \widehat \sigma_{n}^2 / \widehat \rho_{n}^{2\nu} 
-\sigma_0^2 / \rho_0^{2\nu}
\right)$ converges to a $\mathcal N\big(0,2\left( \sigma_0^2 / \rho_0^{2\nu} \right)^2\big)$ distribution. In the next theorem, we show that this asymptotic normality also holds conditionally to $\{Y\in \mathcal E_{\kappa}\}$.

\begin{theorem}\label{th:cond_after_matern}
	For $\kappa=1,2$, we assume that Condition-Grid holds. For $\kappa=0,1,2$, under Condition-$\nu$, the estimator $\widehat \sigma_n^2/ \widehat \rho_n^{2\nu}$ of the microergodic parameter $\sigma_0^2/\rho_0^{2\nu}$ defined by \eqref{def:MLE_matern} and conditioned on $\{ Y \in \mathcal E_{\kappa} \}$ is asymptotically Gaussian distributed. More precisely, 
	\[
	\sqrt n
	\bigg(\frac{\widehat \sigma_{n}^2}{\widehat \rho_{n}^{2\nu}}-\frac{ \sigma_0^2}{ \rho_0^{2\nu}}
	\bigg)
	\xrightarrow[n\to+\infty]{\mathcal{L}| Y \in \mathcal E_{\kappa}} \mathcal N\bigg(0,2\bigg(\frac{\sigma_0^2}{\rho_0^{2\nu}}\bigg)^2\bigg).
	\]
\end{theorem}

\subsection{Asymptotic conditional distribution of the constrained maximum likelihood estimator}\label{ssec:before_matern}

We turn to the constrained log-likelihood and its maximizer. We consider two types of estimation settings obtained by maximizing the constrained log-likelihood \eqref{eq:const:log:lik} under the Mat\'ern model. In the first setting, $\rho=\rho_1$ is fixed and \eqref{eq:const:log:lik} is maximized over $\sigma^2$ (in the case $\rho_1 = \rho_0$ this setting is already covered by Theorem \ref{th:cond_before_var}). In the second setting, \eqref{eq:const:log:lik} is maximized over both $\sigma^2$ and $\rho$. Under the two settings, we show that the cMLE has the same asymptotic distribution as the MLE, conditionally to $\{ Y \in \mathcal E_{\kappa} \}$.

\begin{theorem}[Fixed correlation length parameter $\rho_1$]\label{th:cond_before_matern_i} 
	For $\kappa=1,2$, we assume that Condition-Grid holds. Assume that Condition-$\nu$ and Condition-Sequence hold. 
	Let for $\rho \in [ \rho_l , \rho_u ]$, 
	\begin{align} \label{def:MLE_cond_matern}
	\widehat{\sigma}_{n,c}^2 (\rho)
	\in \underset{ \sigma^2 \in [ \sigma_l^2 , \sigma_u^2 ] }{\argmax \; }
	\mathcal{L}_{n,c}  ( \sigma^2 , \rho ).
	\end{align}
	Let $\rho_1 \in [ \rho_l , \rho_u ] $ be fixed. 
	Then $\widehat \sigma_{n,c}^2( \rho_1 )$ is asymptotically Gaussian distributed for $\kappa=0,1,2$. More precisely,
	\begin{align*}
	\sqrt n \bigg( \frac{\widehat \sigma_{n,c}^2( \rho_1 )}{\rho_1^{2\nu}}-\frac{\sigma_0^2}{\rho_0^{2\nu}} \bigg) \xrightarrow[n\to+\infty]{\mathcal{L} | Y \in \mathcal E_{\kappa}} \mathcal N\bigg(0,2\bigg(\frac{\sigma_0^2}{\rho_0^{2\nu}}\bigg)^2\bigg).
	\end{align*}
	
\end{theorem}

\begin{theorem}[Estimated correlation length parameter]\label{th:cond_before_matern_ii} 
	For $\kappa=1,2$, we assume that Condition-Grid holds.
	Assume that Condition-$\nu$ holds. 
	Let $\widehat \sigma_{n,c}^2( \rho )$ be defined as in \eqref{def:MLE_cond_matern} and let $(\widehat{\sigma}_{n,c}^2,\widehat \rho_{n,c})$ be defined by
	\begin{align*}
	(\widehat \sigma_{n,c}^2,\widehat \rho_{n,c}) \in \underset{(\sigma^2,\rho) \in \Theta}{\argmax \; } \mathcal L_{n,c}(\sigma^2,\rho).
	\end{align*}
	Notice that $\widehat \sigma_{n,c}^2 = \widehat \sigma_{n,c}^2( \widehat \rho_{n,c} )$.
	\begin{itemize}
		\item[(i)] For $\kappa=0$, assume that one of the following two conditions hold.
		\begin{itemize}
			\item[a)] We have $\nu > 1$, $d=1$ and ${\max}_{x \in [0,1]} {\min}_{i=1,\ldots,n} \ |x - x_i| = o( 1/ \sqrt n )$.
			\item[b)] We have $\nu > 2$ and there exists a sequence $(a_n)_{n \in \mathbb{N}}$ with $a_n = o(1/n^{1/4})$ as $n \to \infty$, so that, for all $x \in [0,1]^d$, there exists $d+1$ points $v_1,\ldots,v_{d+1}$ with $\{ v_1,\ldots,v_{d+1}\} \subset \{ x_1,\ldots,x_n \}$, so that $x$ belongs to the convex hull of $v_1,\ldots,v_{d+1}$ and ${\max}_{j=1,\ldots,d+1} \norme{x - v_j} \leqslant a_n$. 
		\end{itemize}  
		\item[(ii)] For $\kappa=1,2$, assume that one of the following two conditions hold.
		\begin{itemize}
			\item[a)] We have $\nu > \kappa+1$, $d=1$ and ${\max}_{x \in [0,1]} {\min}_{i=1,\ldots,n} \ |x - x_i| = o( 1/ \sqrt n )$.
			\item[b)] We have $\nu > \kappa+2$ and the observation points $\{ x_1,\ldots,x_n \}$ are so that, for all $n \geqslant 2^d$, with $N = \lfloor n^{1/d} \rfloor$,
			\[
			\{x_1,\ldots,x_{n}\} 
			\supset
			\bigg\{
			\bigg(\frac{i_1}{N-1},\ldots,\frac{i_d}{N-1}\bigg),\; 0\leqslant i_1,\ldots,i_d\leqslant N-1
			\bigg\}.
			\]
		\end{itemize}  
	\end{itemize}
	Then $\widehat \sigma_{n,c}^2 /  \widehat \rho_{n,c}^{2 \nu}$ is asymptotically Gaussian distributed for $\kappa=0,1,2$. More precisely,  
	\begin{align*}
	\sqrt n \bigg(
	\frac{ 
		\widehat \sigma_{n,c}^2 
	}{
	\widehat \rho_{n,c}^{2\nu} 
}
-\frac{\sigma_0^2}{\rho_0^{2\nu}} 
\bigg) 
\xrightarrow[n\to+\infty]{\mathcal{L} | Y \in \mathcal E_{\kappa} } \mathcal N\bigg(0,2\bigg(\frac{\sigma_0^2}{\rho_0^{2\nu}}\bigg)^2\bigg).
\end{align*}
\end{theorem}

In Theorem \ref{th:cond_before_matern_ii}, we assume that $\nu$ is larger than in Condition-$\nu$, and we assume that the observation points have specific quantitative space filling properties. The condition (i) b) also implies that a portion of the observation points are located in the corners and borders of $[0,1]^d$. Furthermore, the condition (ii) b) implies that the majority of the observation points are located on regular grids. We believe that these two last conditions could be replaced by milder ones, at the cost of similar proofs but more cumbersome than the present ones.

We make stronger assumptions in Theorem \ref{th:cond_before_matern_ii} than in Theorem \ref{th:cond_before_matern_i} because the former is more challenging than the latter. Indeed, since $\rho = \rho_1$ is fixed in Theorem \ref{th:cond_before_matern_i}, we can use the equivalence of two fixed Gaussian measures in order to obtain asymptotic properties of the conditional mean function of $Y$ under $k_{1,\rho_1,\nu}$ (see the developments following \eqref{eq:cv_Bn} in the proofs).
This is not possible anymore when considering the conditional mean function of $Y$ under $k_{1,\widehat \rho_{n,c},\nu}$, where $\widehat \rho_{n,c}$ is random. Hence, we use other proof techniques, based on reproducing kernel Hilbert spaces, for studying this conditional mean function, for which the above additional conditions are needed. We refer for instance to the developments following \eqref{eq:to:show:limsup:mn}  in the appendix for more details.

\section{Numerical results}\label{sec:num}
In this section, we illustrate numerically the conditional asymptotic normality of the MLE and the cMLE of the microergodic parameter for the Mat\'ern 5/2 covariance function. The numerical experiments were implemented using the R package ``LineqGPR'' \citep{LopezLopera2017LineqGP}.


\subsection{Experimental settings}
\label{subsec:numsettings}

We let $d=1$ in the rest of the section.
Since the event $\{Y \in \mathcal E_{\kappa}\}$ can not be simulated exactly in practice, we consider the piecewise affine interpolation $Y_m$ of $Y$ at $t_1, \ldots , t_m \in [0,1]$, with $m > n$ \citep{Maatouk2017GPineqconst,LLBDR17}. Then, the event $\{Y \in \mathcal E_{\kappa}\}$ is approximated by the event $\{Y_m \in \mathcal E'_{\kappa}\}$, where $\mathcal E'_{0}$ (respectively $\mathcal E'_{1}$, $\mathcal E'_{2}$) is the set of continuous bounded between $\ell$ and $u$ (resp. increasing, convex) functions. We can simulate efficiently $Y_m$ conditionally to $\{ Y_m \in \mathcal E'_{\kappa} \}$ by using Markov Chain Monte Carlo procedures \citep[see, for instance,][]{Pakman2014Hamiltonian}.

In Section \ref{sec:num}, we consider the Mat\'ern 5/2  function defined by
\begin{equation*}
k_{\theta,5/2}(x) = \sigma^2 \bigg(1 + \frac{\sqrt{5} |x|}{\rho} + \frac{5}{3} \frac{x^2}{\rho^2}\bigg) \exp\bigg\{-\frac{\sqrt{5}|x|}{\rho} \bigg\},
\end{equation*}
for $x \in \R$ and with $\theta = (\sigma^2,\rho)$.
Remark that $k_{\theta,5/2}$ is obtained by the parametrization of \citep{roustant2012dicekriging,LopezLopera2017LineqGP} rather than that of Section \ref{ssec:model_matern}. For an easy reading, we keep the same notation.

\subsection{Numerical results when $\rho_0$ is known}
\label{ss:resultsCase1}

We let $m = 300$ and $x_1,...,x_n$ be equispaced in $[0,1]$ in the rest of Section \ref{sec:num}.
For $\kappa = 0,1$, we generate $N = 1,000$ trajectories of $Y_m$ given $\{ Y_m \in \mathcal E'_{\kappa} \}$. For each of these trajectories, we compute the MLE $\bar{\sigma}_{m,n}^2 ( \rho_0 )$ with $\bar{\sigma}_{m,n}^2 ( \rho ) = y^\top R_{\rho}^{-1} y / n = \argmax_{\sigma^2 \in (0 , \infty)} \mathcal{L}_n(\sigma^2 , \rho)$, where $R_{\rho} = \big[ \widetilde k_{1,\rho,5/2} ( x_i , x_j )  \big]_{i,j=1,\ldots,d}$ for $\rho \in (0 , \infty)$. Then we evaluate the cMLE as follows. We let $m_{m,n,\rho,y}$ be the conditional mean function of $Y_m$ given $y$ under covariance function $\widetilde k_{1,\rho,5/2}$. We simulate $n_t = 1,000$ trajectories $Z_1,...,Z_{n_t}$ of a Gaussian process with zero mean function and covariance function $\widetilde k_{m,n,1,\rho_0,5/2}$, where $\widetilde k_{m,n,1,\rho,5/2}$ is the covariance function of $Y_m$ given $y$ under covariance function $\widetilde k_{1,\rho,5/2}$. Then we let $B_n(\sigma^2 , \rho_0)$ be approximated by $\ln\big((1 / n_t) \sum_{i=1}^{n_t} \mathds{1}_{ m_{m,n,\rho_0,y} + \sigma Z_i \in \mathcal{E}'_{\kappa} } \big)$. The term $A_n(\sigma^2,\rho_0)$ can be easily approximated as it does not depend on the trajectory of $Y_m$ under consideration. We maximize the resulting approximation of $\mathcal{L}_{n,c}(\sigma^2 , \rho_0)$ on $1,000$ equispaced values of $\sigma^2$ between $(1 - 4\sqrt{2/n}) \sigma_0^2$ and $(1 + 4\sqrt{2/n}) \sigma_0^2$, yielding the approximated cMLE estimator $\widehat{\sigma}_{m,n,c}^2( \rho_0 )$. 

\begin{figure}[t!]
	\centering
	\subfigure[\label{subfig:BoundednessCase1n20} $n = 20$]{\includegraphics[width = 0.33\textwidth]{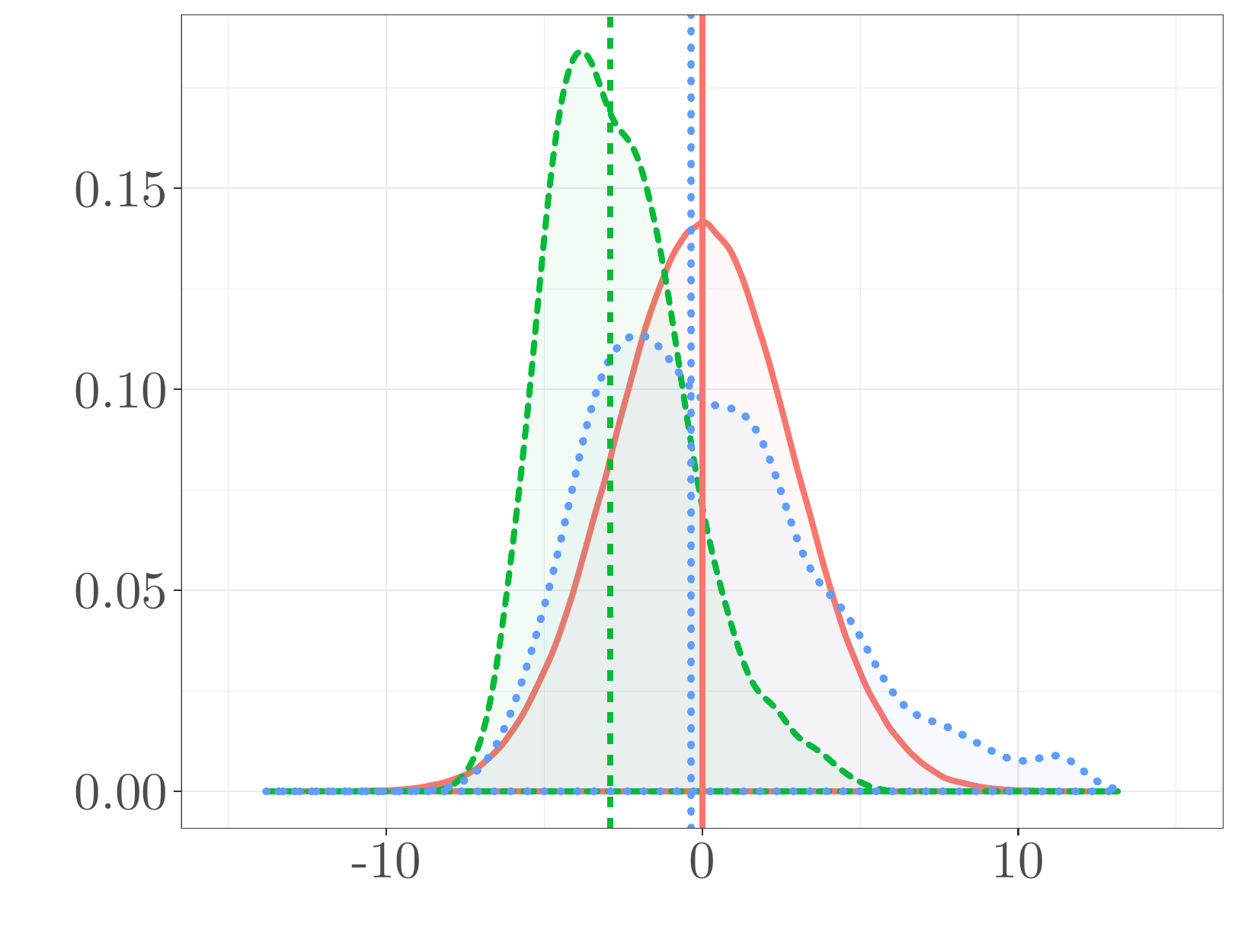}}
	\subfigure[\label{subfig:BoundednessCase1n50} $n = 50$]{\includegraphics[width = 0.33\textwidth]{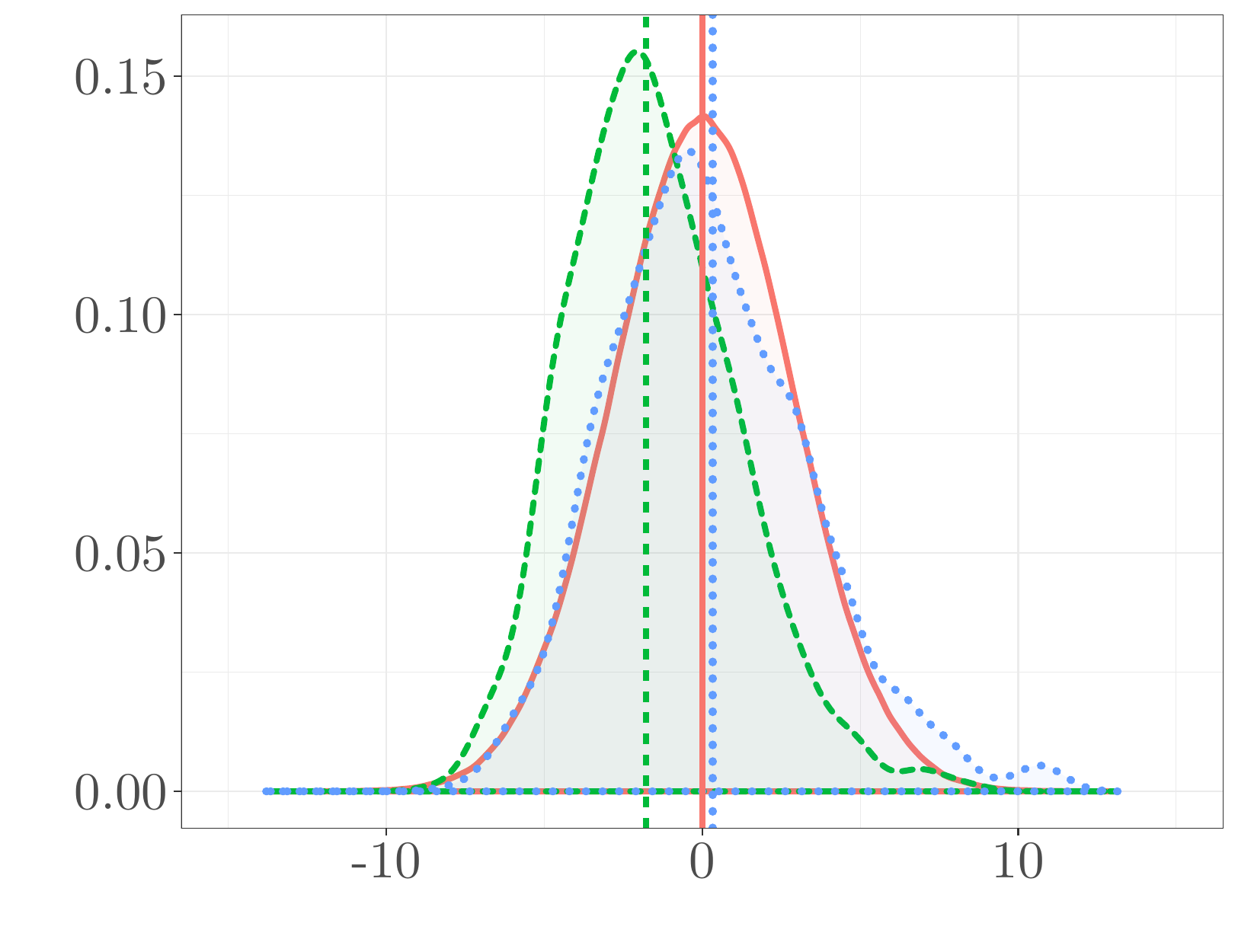}}
	\subfigure[\label{subfig:BoundednessCase1n80} $n = 80$]{\includegraphics[width = 0.33\textwidth]{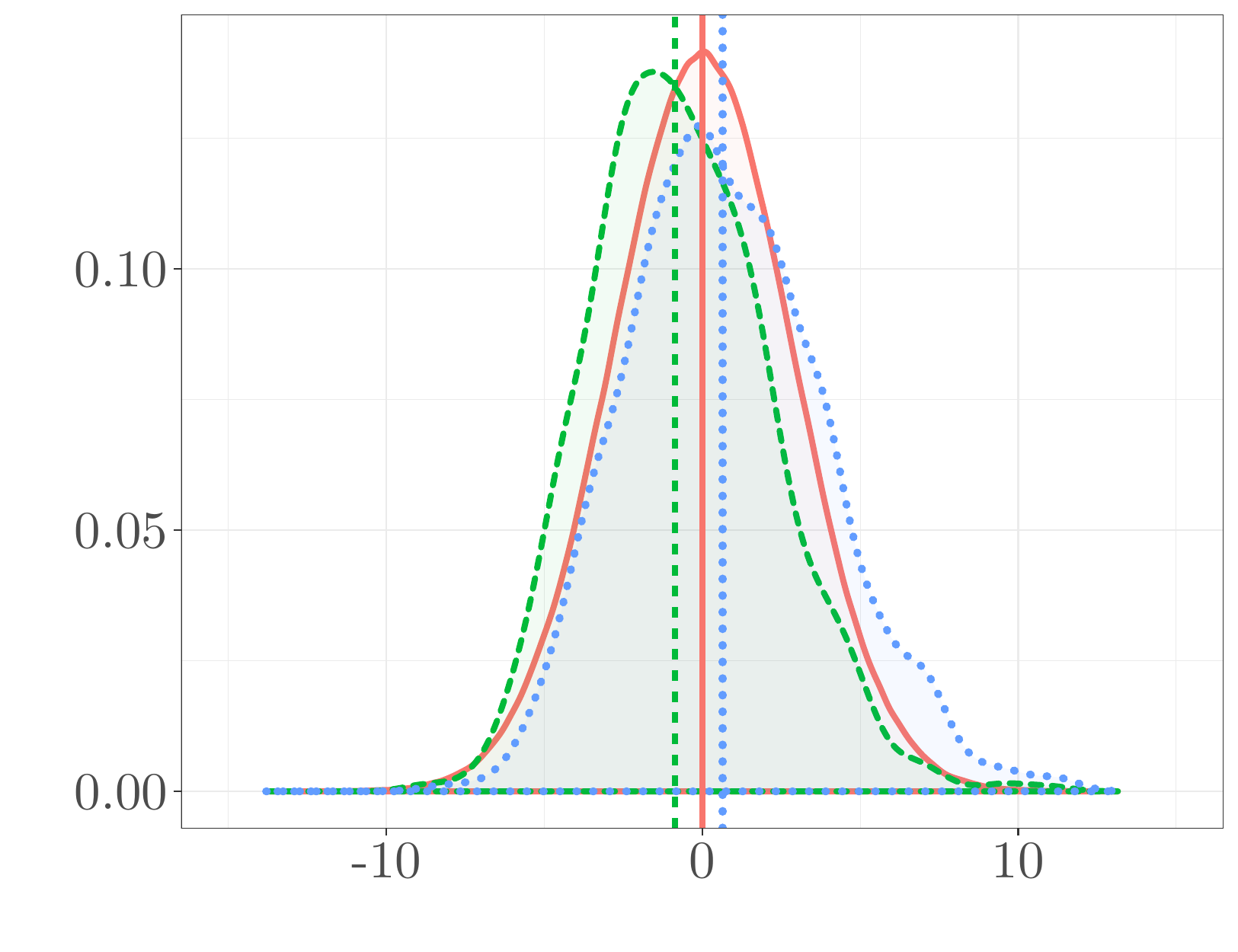}}
	\caption{Asymptotic conditional distributions of the variance parameter estimators under boundedness constraints. Here $(\sigma_0^2,\rho_0)=(2,0.2)$. Each panel shows: the limit distribution $\mathcal{N}(0, 2\sigma_0^4)$ (red, solid lines), the conditional distribution of the MLE (green, dashed lines), and the conditional distribution of the cMLE (blue, dotted lines). The vertical lines represent the median values of the distributions. Each sub-caption shows the number of observations $n$ used for the estimations.}
	\label{fig:AsympNormalityBoundednessCase1}
\end{figure}

In Figure \ref{fig:AsympNormalityBoundednessCase1}, we report the results for $\kappa = 0$ (boundedness constraints) with $(\sigma_0^2,\rho_0) = (2 , 0.2)$ and $n = 20, 50, 80$.
We show the probability density functions obtained from the samples $\left\{ n^{1/2} ( \bar{\sigma}_{m,n}^2 ( \rho_0 )_i - \sigma_0^2) \right\}_{i=1,...,N} $ and $\left\{ n^{1/2} ( \widehat{\sigma}_{m,n,c}^2 ( \rho_0 )_i - \sigma_0^2  ) \right\}_{i=1,...,N} $ obtained as discussed above. We also plot the probability density function of the limit $\mathcal N(0 , 2\sigma_0^4)$ distribution.
We observe that for a small number of observations, e.g. $n=20$, the distribution of the cMLE is closer to the limit distribution than that of the MLE in terms of median value. We also observe that, as $n$ increases, both distributions become more similar to the limit one. Nevertheless, the cMLE exhibits faster convergence.

\begin{figure}
	\centering
	\subfigure[\label{subfig:MonotonicityCase1n5} $n = 5$]{\includegraphics[width = 0.33\textwidth]{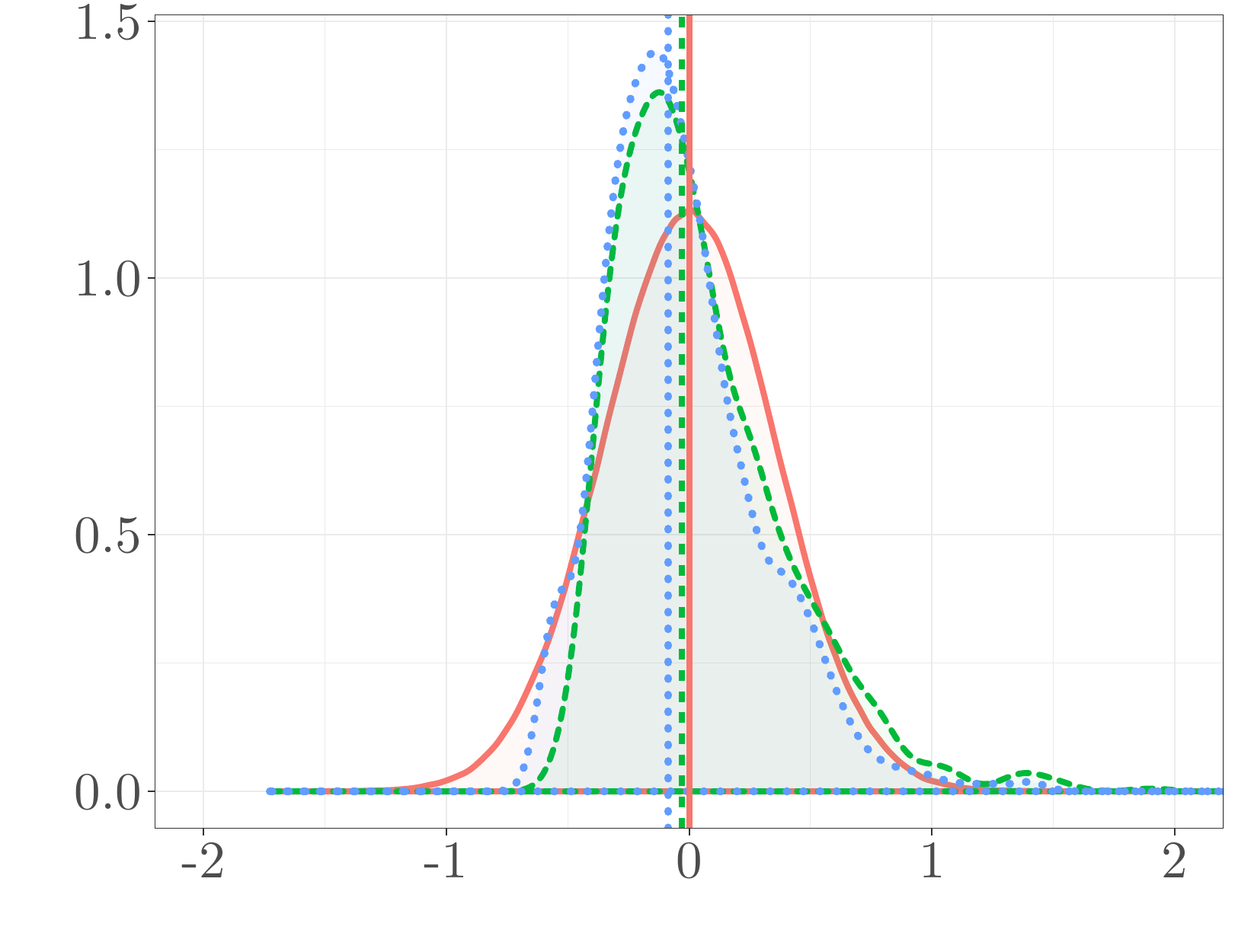}}
	\subfigure[\label{subfig:MonotonicityCase1n20} $n = 20$]{\includegraphics[width = 0.33\textwidth]{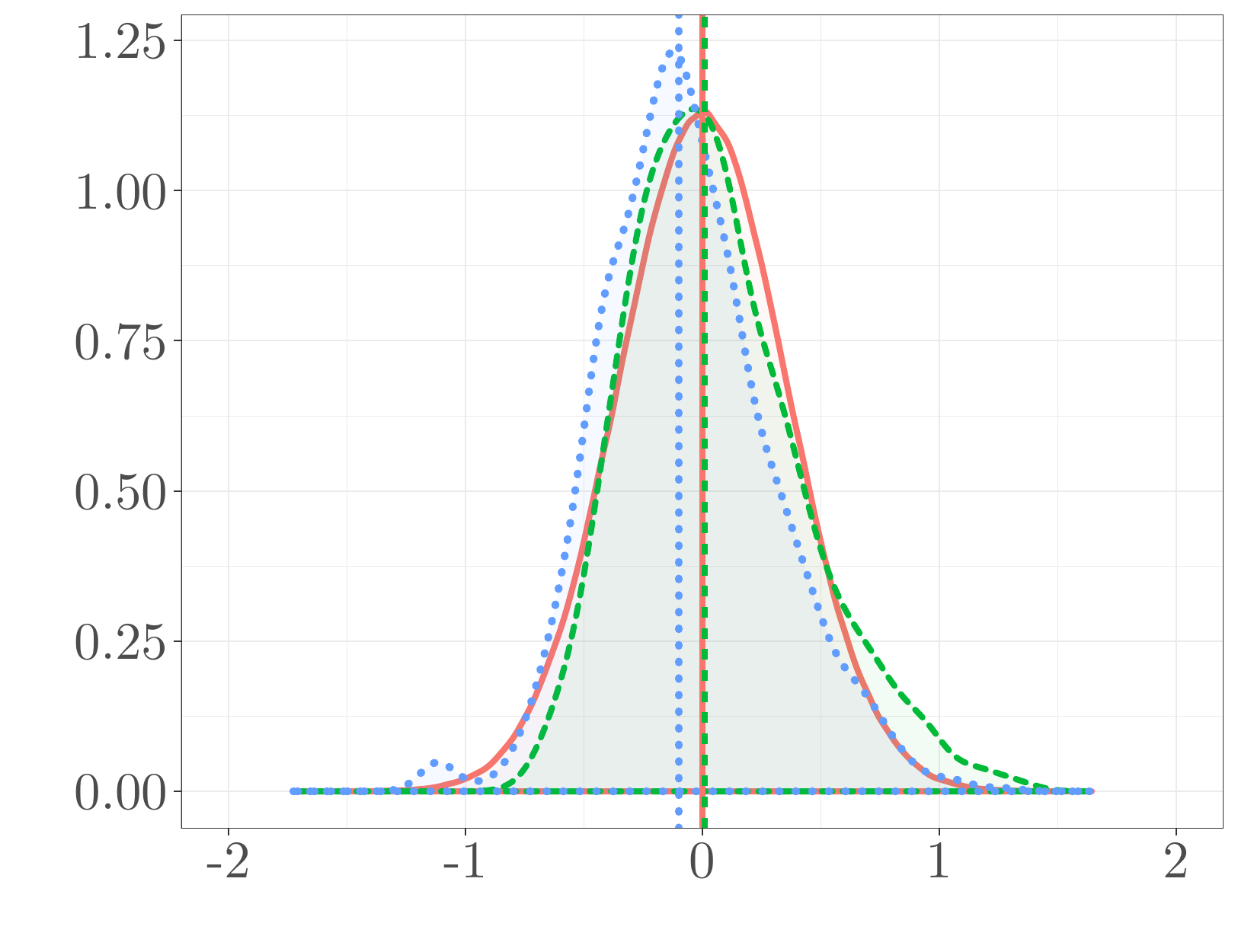}}
	\caption{Asymptotic conditional distributions of the variance parameter estimators under monotonicity constraints. Here $(\sigma_0^2,\rho_0)=(0.5^2,1)$. Panel description is the same as in Figure \ref{fig:AsympNormalityBoundednessCase1}.}
	\label{fig:AsympNormalityMonotonicityCase1}
\end{figure}

In Figure \ref{fig:AsympNormalityMonotonicityCase1}, we report the same quantities for $\kappa = 1$ (monotonicity constraints) and for $(\sigma_0^2,\rho_0) = (0.5^2,1)$. 
In this case, we observe that the distributions of both the MLE and the cMLE are close to the limit one even for small values of $n$ ($n=5, 20$).

\subsection{Numerical results when $\rho_0$ is unknown}
\label{ss:resultsCase3}

We let $\kappa = 0$, $(\sigma_0^2,\rho_0) = ( 2, 0.2)$ and $n = 20, 50, 80$. 
We proceed similarly as in the case where $\rho_0$ is known. To compute the MLE 
$( \bar{\sigma}_{m,n}^2( \widehat{\rho}_{m,n}) , \widehat{\rho}_{m,n} )$
of
$(\sigma_0^2 , \rho_0)$,
we maximize $ \mathcal{L}_n( \bar{\sigma}_{m,n}^2  (\rho ) , \rho)$  over a finite grid of values for $\rho$. To compute the cMLE
$( \widehat{\sigma}_{m,n,c}^2( \widehat{\rho}_{m,n,c}) , \widehat{\rho}_{m,n,c} )$
of
$(\sigma_0^2 , \rho_0)$, we evaluate $\ln\big((1 / n_t) \sum_{i=1}^{n_t} \mathds{1}_{ m_{m,n,\rho,y} + \sigma Z_{\rho,i} \in \mathcal{E}'_{\kappa} }\big) $ over $100^2$ pairs $(\sigma_{i,j}^2 , \rho_{i})_{i,j = 1,\ldots,100}$. Here $Z_{\rho,i}$ is generated as in Section \ref{ss:resultsCase1} but with $\rho_0$ replaced by $\rho$,  for $i=1,\ldots,100$. Then $\rho_1,\ldots,\rho_{100}$ are equispaced in $[0.1,0.3]$ and  for $i=1,\ldots,100$, $\sigma^2_{i,1},\ldots,\sigma^2_{i,100}$ are equispaced in 
\[
\rho_i^{5} \left[  \frac{\sigma_0^2}{\rho_0^{5}} - \frac{4\sqrt{2}}{\sqrt{n}}   \frac{\sigma_0^2}{\rho_0^{5}}
,
\frac{\sigma_0^2}{\rho_0^{5}} + \frac{4\sqrt{2}}{\sqrt{n}}   \frac{\sigma_0^2}{\rho_0^{5}}
\right].
\]
Hence, the estimator of the microergodic parameter $\sigma_0^2/\rho_0^{5}$ is restricted to be at distance less than $4$ times the asymptotic standard deviation of the microergodic parameter.

\begin{figure}[t!]
	\centering
	\subfigure[\label{subfig:BoundednessCase3n20b} $n = 20$]{\includegraphics[width = 0.33\textwidth]{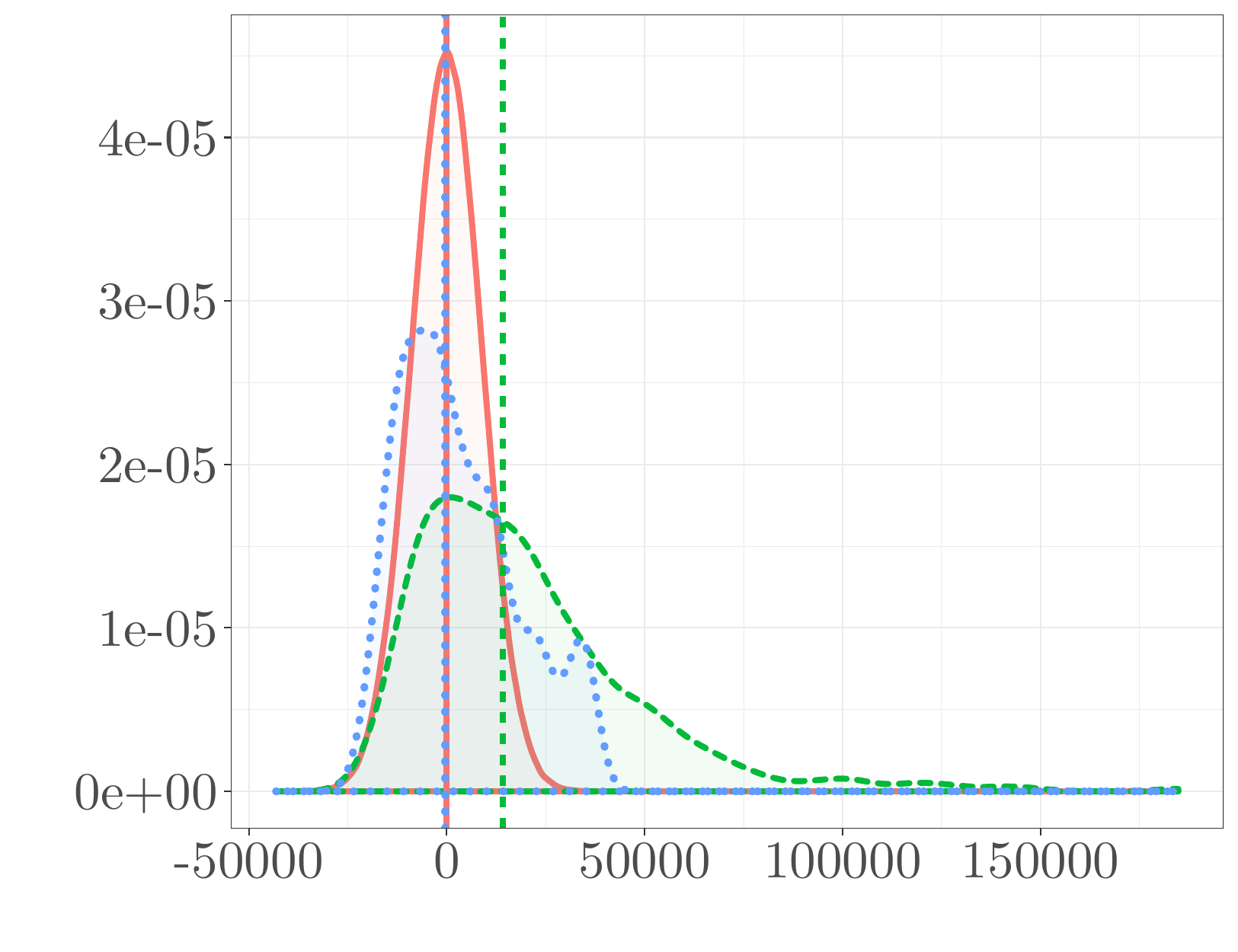}}
	\subfigure[\label{subfig:BoundednessCase3n50b} $n = 50$]{\includegraphics[width = 0.33\textwidth]{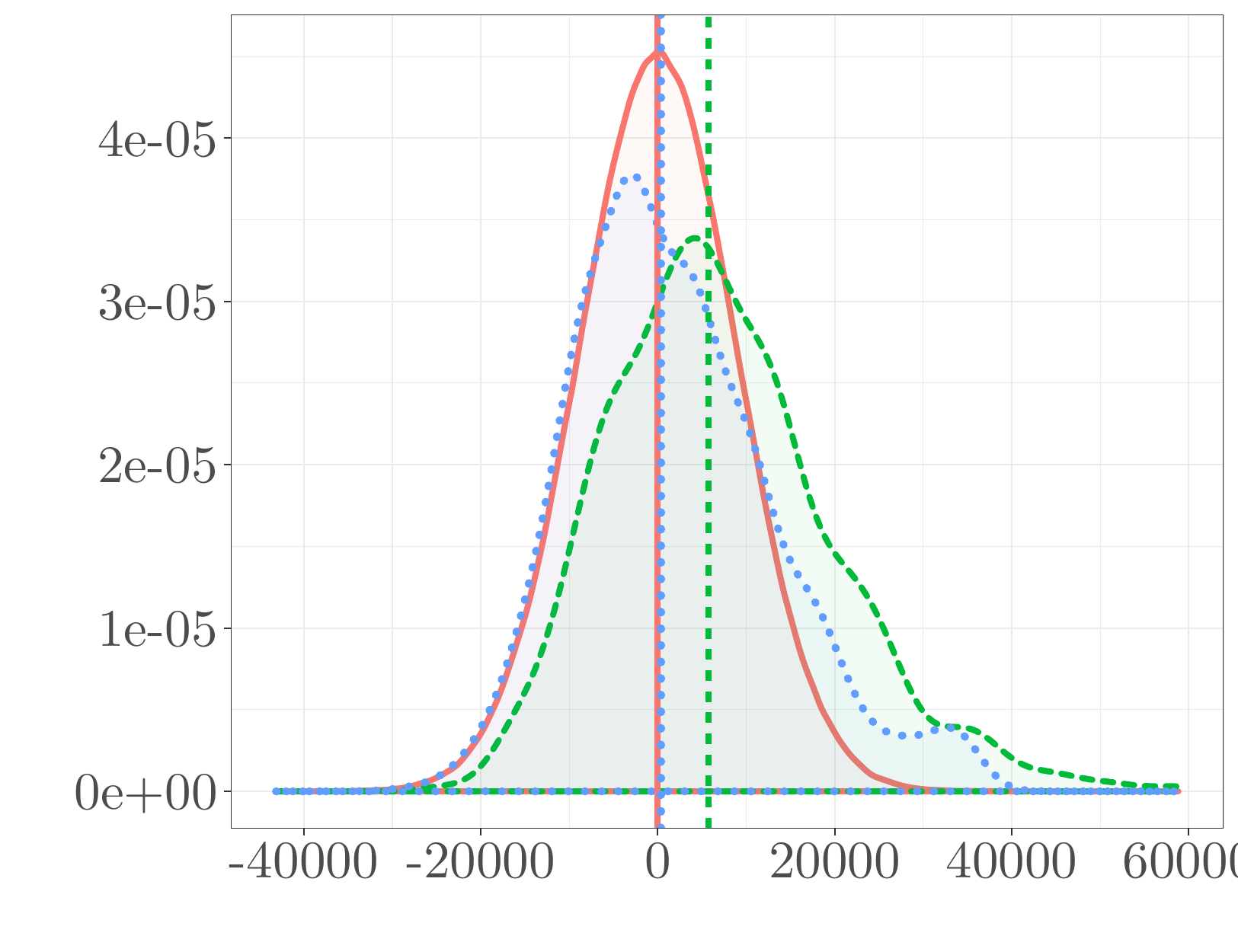}}
	\subfigure[\label{subfig:BoundednessCase3n80b} $n = 80$]{\includegraphics[width = 0.33\textwidth]{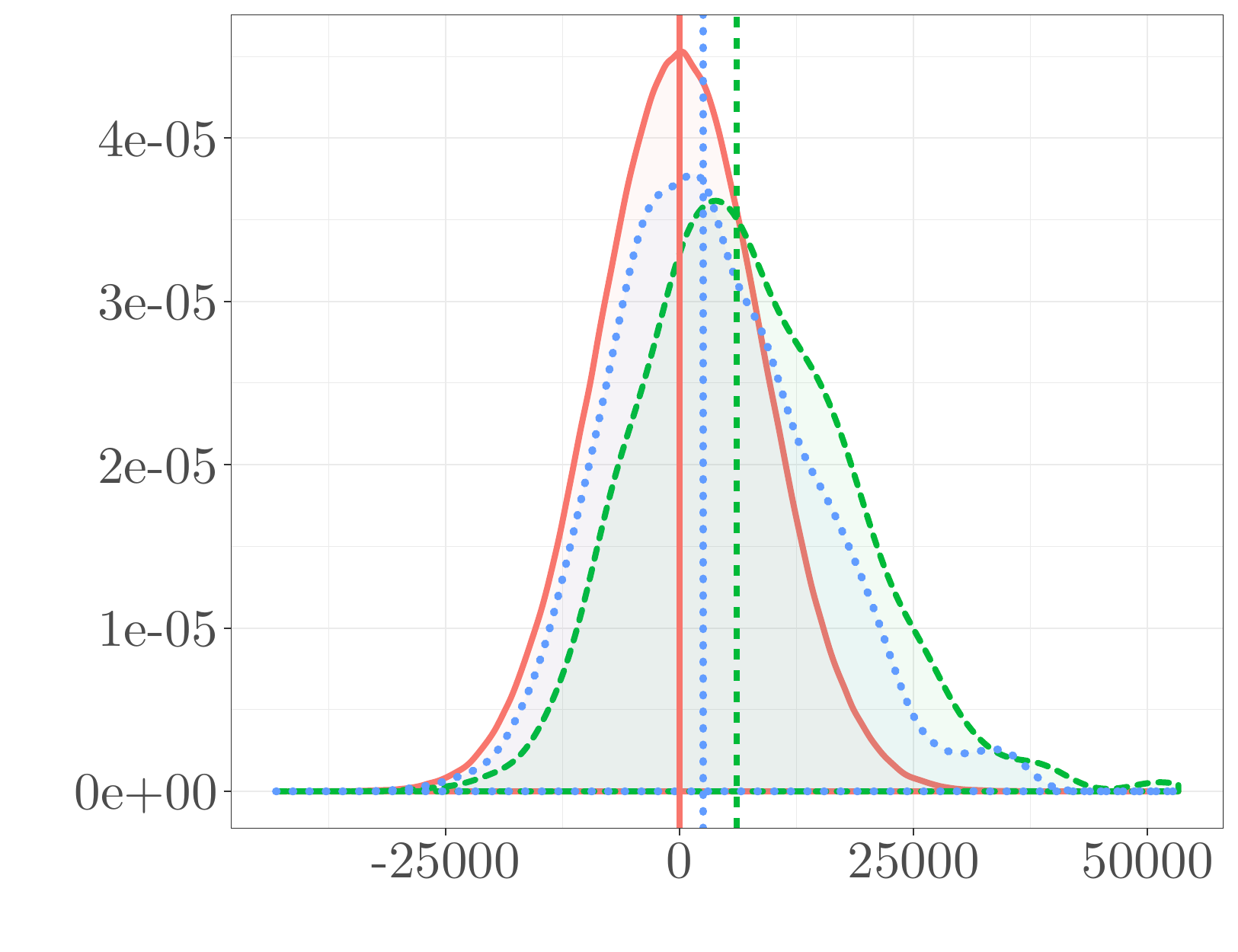}}
	\caption{Asymptotic conditional distributions of the microergodic parameter estimators for the isotropic $\nu =5/2$ Mat\'ern model under boundedness constraints. Here $(\sigma_0^2,\rho_0)=(2,0.2)$. Panel description is the same as in Figure \ref{fig:AsympNormalityBoundednessCase1}, with $\mathcal{N}(0, 2\sigma_0^4)$ replaced by $\mathcal{N}(0, 2(\sigma_0^2/\rho_0^5)^2)$.}
	\label{fig:AsympNormalityBoundednessCase3}
\end{figure}

In Figure \ref{fig:AsympNormalityBoundednessCase3},
we show the probability density functions obtained from the samples $\big\{ n^{1/2} ( \bar{\sigma}_{m,n}^2 ( \widehat{\rho}_{m,n} )_i / \widehat{\rho}_{m,n,i}^5
-
\sigma_0^2 / \rho_0^5 
) 
\big\}_{i=1,...,N}$
and
$\big\{ n^{1/2} ( \widehat{\sigma}_{m,n,c}^2 ( \widehat{\rho}_{m,n,c} )_i  / \widehat{\rho}_{m,n,c,i}^5
-
\sigma_0^2 / \rho_0^5 
) 
\big\}_{i=1,...,N} $, with $N=1,000$. 
Similarly to Section \ref{ss:resultsCase1}, we observe that the distribution of the cMLE tends to be closer to the limit one, than that of the MLE. Moreover, the convergence with the cMLE is faster than with the MLE in terms of median value. 

\section{Extensions} \label{sec:extensions}

\subsection{Results on prediction}

In the next proposition, we show that, conditionally to the inequality constraints, the predictions obtained when taking the constraints into account are asymptotically equal to the standard (unconstrained) Kriging predictions. Furthermore, the same is true when comparing the conditional variances obtained with and without accounting for the constraints.

\begin{proposition} \label{proposition:prediction}
	Let $\kappa = 0,1,2$ be fixed.
	Consider a Gaussian process $Y$ on $[0,1]^d$ with mean function zero and covariance function $\widetilde{k}$ of the form $\widetilde{k}(u,v) = k(u-v)$ for $u,v \in [0,1]^d$, where $k : \mathbb{R}^d \to \mathbb{R}$ satisfies Condition-Var (with $k_1$ replaced by $k$).
	Assume that Condition-Sequence holds. Recall that $y = (Y(x_1),\ldots,Y(x_n))^{\top}$ is the observation vector. For $x_0 \in [0,1]^d$, let $\widehat{Y}(x_0) = \mathbb{E}[Y(x_0) | y]$, $\widehat{\sigma}(x_0)^2 = \Var(Y(x_0) | y)$, $\widehat{Y}_c(x_0) = \mathbb{E}[Y(x_0) | y , Y \in \mathcal{E}_{\kappa} ]$, and $\widehat{\sigma}_c(x_0)^2 = \Var(Y(x_0) | y, Y \in \mathcal{E}_{\kappa})$. Then when $x_0 \not \in \{x_i \}_{i \in \mathbb{N}}$, we have
	\begin{equation} \label{eq:pred:means}
	\frac{ \widehat{Y}(x_0) - \widehat{Y}_c(x_0) }{ \widehat{\sigma}(x_0) }
	=o_{\P | Y \in \mathcal{E}_{\kappa}}(1),
	\end{equation}
	and
	\begin{equation} \label{eq:pred:variance}
	\frac{
		\widehat{\sigma}(x_0)^2-
		\widehat{\sigma}_c(x_0)^2}{
		\widehat{\sigma}(x_0)^2
	}
	=o_{\P | Y \in \mathcal{E}_{\kappa}}(1).
	\end{equation}
	
\end{proposition}

In Proposition \ref{proposition:prediction}, when taking the constraints into account or not, the predictions converge to the true values and the conditional variances converge to zero. Thus, the results in Proposition \ref{proposition:prediction} are given on a relative scale, by dividing the difference of predictions by the conditional standard deviation (without constraints), and by dividing the difference of conditional variances by the conditional variance (without constraints).

Similarly as for estimation in Sections \ref{sec:var} and \ref{sec:matern}, the conclusion of Proposition \ref{proposition:prediction} is that the constraints do not have an asymptotic impact on prediction.  

When there is no constraints, significant results on using misspecified covariance functions that are asymptotically equivalent to the true one have been obtained in \citep{AEPRFMCF,BELPUICF,UAOLPRFUISOS}.
Let $\widetilde{k}$, $\widehat{\sigma}(x_0)$, $\widehat{Y}(x_0)$ and $\widehat{Y}_c(x_0)$ be as in Proposition \ref{proposition:prediction}. Let $k_1$ satisfy Condition-Var and let $\widetilde{k}_1$ be defined from $k_1$ as in Proposition \ref{proposition:prediction}. Let the Gaussian measures of the Gaussian processes with mean functions zero and covariance functions $\widetilde{k}$ and $\widetilde{k}_1$ on $[0,1]^d$ be equivalent \citep[see][]{stein99interpolation}. Let $\widehat{\sigma}_1(x_0)$, $\widehat{Y}_1(x_0)$, $\widehat{Y}_{c,1}(x_0)$, and $\widehat{\sigma}_{c,1}(x_0)$ be defined as $\widehat{\sigma}(x_0)$, $\widehat{Y}(x_0)$, $\widehat{Y}_c(x_0)$ and $\widehat{\sigma}_c(x_0)$, when taking the conditional expectations with respect to $\widetilde{k}_1$ rather than $\widetilde{k}$. Then it is shown in \citep{AEPRFMCF,BELPUICF,UAOLPRFUISOS} \citep[see also][Chapter 4, Theorem 8]{stein99interpolation} that, when $x_0 \not \in \{x_i \}_{i \in \mathbb{N}}$,
\begin{equation} \label{eq:equivalence:prediction:pred}
\frac{\widehat{Y}_1(x_0) - \widehat{Y}(x_0)}{ \widehat{\sigma}(x_0) }
= o_{\P}(1),
\end{equation}
and
\begin{equation} \label{eq:equivalence:prediction:var}
\frac{\widehat{\sigma}_1(x_0)^2 - \widehat{\sigma}(x_0)^2}{ \widehat{\sigma}(x_0)^2 }
= o(1).
\end{equation}
Both expressions above mean that the predictions and conditional variances obtained from equivalent Gaussian measures are asymptotically equivalent.
A corollary of our Proposition \ref{proposition:prediction} is that this equivalence remains true when the predictions and conditional variances are calculated accounting for the inequality constraints.

\begin{corollary} \label{corollary:equivalence:pred:constraints}
	Let $\kappa = 0,1,2$ be fixed.
	Consider a Gaussian process $Y$ on $[0,1]^d$ with mean function zero and covariance function $\widetilde{k}$ of the form $\widetilde{k}(u,v) = k(u-v)$ for $u,v \in [0,1]^d$, where $k : \mathbb{R}^d \to \mathbb{R}$ satisfies Condition-Var.
	Assume that Condition-Sequence holds. Consider a covariance function $\widetilde{k}_1$ of the form $\widetilde{k}_1(u,v) = k_1(u-v)$ for $u,v \in [0,1]^d$, where $k_1 : \mathbb{R}^d \to \mathbb{R}$ satisfies Condition-Var. 
	
	Let the Gaussian measures of Gaussian processes with mean functions zero and covariance functions $\widetilde{k}$ and $\widetilde{k}_1$ on $[0,1]^d$ be equivalent. 
	Then when $x_0 \not \in \{x_i \}_{i \in \mathbb{N}}$, we have
	\begin{equation*} 
	\frac{\widehat{Y}_{c,1}(x_0) - \widehat{Y}_{c}(x_0)}{
		\widehat{\sigma}_c(x_0)
	}
	= o_{\P|Y\in \mathcal E_{\kappa}}(1), \quad \mbox{and} \quad
	\frac{
		\widehat{\sigma}_{c,1}(x_0)^2-
		\widehat{\sigma}_c(x_0)^2
	}
	{
		\widehat{\sigma}_c(x_0)^2
	}
	= o_{\P|Y\in \mathcal E_{\kappa}}(1).
	\end{equation*}
	
\end{corollary}

Finally, an important question for Gaussian processes is to assess the asymptotic accuracy of predictions obtained from (possibly consistently) estimated covariance parameters. In this section, we have restricted the asymptotic analysis of prediction to fixed (potentially misspecified) covariance parameters.

When no constraints are considered, and under increasing-domain asymptotics, predictions obtained from consistent estimators of covariance parameters are generally asymptotically optimal \citep{bachoc14asymptotic,bachoc2018asymptotic}. 
Under fixed-domain asymptotics, without considering constraints, the predictions obtained from estimators of the covariance parameters can be asymptotically equal to those obtained from the true covariance parameters \citep{putter01effect}. It would be interesting, in future work, to extend the results given in \citep{putter01effect}, to the case of inequality constraints. This could be carried out by making Proposition \ref{proposition:prediction} uniform over subspaces of covariance parameters, and by following a similar approach as for proving Corollary \ref{corollary:equivalence:pred:constraints}.

\subsection{Microergodic parameter estimation for the isotropic Wendland model}

In this section, we let $d=1,2$ or $3$ and extend the results for the Mat\'ern covariance functions of Section \ref{sec:matern} to the isotropic Wendland family of covariance functions on $[0,1]^d$ \citep{gneiting2002compactly,bevilacqua2019estimation}. Here $k_{\theta} = k_{\theta,s,\mu}$, with $\theta=(\sigma^2,\rho)$, is given by
\[
k_{\theta,s,\mu}(x) 
=
\sigma^2 \phi_{s,\mu} \left(  \frac{ \norme{x} }{ \rho } \right),
\]
for $x \in \mathbb{R}^d$ with, for $t \geqslant 0$,
\[
\phi_{s,\mu}(t)
=
\begin{cases}
\frac{1}{B(2s,\mu+1)}
\int_{\norme{x}}^1
u(u^2 - \norme{x}^2)^{s-1}  (1-u)^\mu du
& ~ ~ \mbox{if} ~ ~ \norme{x} < 1, \\
0  & ~ ~ \mbox{else.}
\end{cases}
\]
The parameters $s >0$ and $\mu \geqslant (d+1)/2 +s$ are considered to be fixed and known.
The Wendland covariance function is given by $\widetilde k_{\theta,s,\mu}(u,v)=k_{\theta,s,\mu}(u-v)$. 
The parameter $s$ drives the smoothness of the Wendland covariance function, similarly as for the Mat\'ern covariance function \citep{bevilacqua2019estimation}. The parameters $\sigma^2 > 0$ and $\rho > 0$ are interpreted similarly as for the Mat\'ern covariance functions and are to be estimated. 
We remark that, for appropriate equality conditions on $\nu$ (see Section \ref{sec:matern}), $s$ and $\mu$, the Gaussian measures obtained from the Wendland and Mat\'ern covariance functions are equivalent \citep{bevilacqua2019estimation}. The Wendland covariance function is compactly supported, which is a computational benefit \citep{bevilacqua2019estimation}.

Let us define the MLE $(\widehat{\sigma}_n^2 , \widehat{\rho}_n)$ in the exact same way as in Section \ref{ssec:after_matern} but for the Wendland covariance functions, with $\Theta = [\sigma_l^2 , \sigma_u^2] \times [ \rho_l , \rho_u ]$ fixed as in Section \ref{ssec:after_matern} and with $\sigma_l^2/(\rho_l^{1+2 s}) < \sigma_0^2/(\rho_0^{1+2 s}) < \sigma_u^2/(\rho_u^{1+2 s})$. 

It is shown in \citep{bevilacqua2019estimation} that the parameters $\sigma_0^2$ and $\rho_0$ cannot be estimated consistently but that the parameter $\sigma_0^2 / \rho_0^{1+2s}$ can. Furthermore, $\sqrt n
\left( \widehat \sigma_{n}^2 / \widehat \rho_{n}^{1+2s} 
-\sigma_0^2 / \rho_0^{1+2s}
\right)$ converges to a $\mathcal N\big(0,2\left( \sigma_0^2 / \rho_0^{1+2s} \right)^2\big)$ distribution. Then, we can  extend Theorem \ref{th:cond_after_matern}, providing the asymptotic conditional distribution of the MLE of the microergodic parameter for the Mat\'ern model, to the Wendland model. \\

\tbf{Condition-$s,\mu$}.
We assume that $\mu \geqslant d/2 + 7/2 +s$ and for $\kappa=1$ (respectively $\kappa=2$), we assume that $s + 1/2>1$ (resp. $s+1/2>2$).\\

\begin{theorem}\label{th:cond_after_wendland}
	For $\kappa=1,2$, we assume that Condition-Grid holds. For $\kappa=0,1,2$, under Condition-$s,\mu$, the MLE $\widehat \sigma_n^2/ \widehat \rho_n^{1+2s}$ of the microergodic parameter $\sigma_0^2/\rho_0^{1+2s}$, conditioned on $\{ Y \in \mathcal E_{\kappa} \}$, is asymptotically Gaussian distributed. More precisely,  
	\[
	\sqrt n
	\bigg(\frac{\widehat \sigma_{n}^2}{\widehat \rho_{n}^{1+2s}}-\frac{ \sigma_0^2}{ \rho_0^{1+2s}}
	\bigg)
	\xrightarrow[n\to+\infty]{\mathcal{L}| Y \in \mathcal E_{\kappa}} \mathcal N\bigg(0,2\bigg(\frac{\sigma_0^2}{\rho_0^{1+2s}}\bigg)^2\bigg).
	\]
\end{theorem}

Now we define the cMLEs $\widehat{\sigma}_{n,c}^2(\rho)$ and $(\widehat \sigma_{n,c}^2,\widehat \rho_{n,c})$ as in Section \ref{ssec:before_matern}, but for the Wendland covariance functions. Then we can extend Theorems \ref{th:cond_before_matern_i}  and  \ref{th:cond_before_matern_ii} to the Wendland model. 

\begin{theorem}[Fixed correlation length parameter $\rho_1$]\label{th:cond_before_wendland_i} 
	For $\kappa=1,2$, we assume that Condition-Grid holds. Assume that Condition-$s,\mu$ and Condition-Sequence hold.
	Let $\rho_1 \in [ \rho_l , \rho_u ] $ be fixed.
	Then $\widehat \sigma_{n,c}^2( \rho_1 )$ is asymptotically Gaussian distributed for $\kappa=0,1,2$. More precisely, 
	\begin{align*}
	\sqrt n \bigg( \frac{\widehat \sigma_{n,c}^2( \rho_1 )}{\rho_1^{1+2s}}-\frac{\sigma_0^2}{\rho_0^{1+2s}} \bigg) \xrightarrow[n\to+\infty]{\mathcal{L} | Y \in \mathcal E_{\kappa}} \mathcal N\bigg(0,2\bigg(\frac{\sigma_0^2}{\rho_0^{1+2s}}\bigg)^2\bigg).
	\end{align*}
	
\end{theorem}

\begin{theorem}[Estimated correlation length parameter]\label{th:cond_before_wendland_ii} 
	For $\kappa=1,2$, we assume that Condition-Grid holds.
	Assume that Condition-$s,\mu$ holds. Assume the same two conditions (i) and (ii) as in Theorem \ref{th:cond_before_matern_ii}, but with $\nu$ replaced by $s + 1/2$.
	Then $\widehat \sigma_{n,c}^2 /  \widehat \rho_{n,c}^{1 + 2s}$ is asymptotically Gaussian distributed for $\kappa=0,1,2$. More precisely,  
	\begin{align*}
	\sqrt n \bigg(
	\frac{ 
		\widehat \sigma_{n,c}^2 
	}{
	\widehat \rho_{n,c}^{1+2s} 
}
-\frac{\sigma_0^2}{\rho_0^{1+2s}} 
\bigg) 
\xrightarrow[n\to+\infty]{\mathcal{L} | Y \in \mathcal E_{\kappa} } \mathcal N\bigg(0,2\bigg(\frac{\sigma_0^2}{\rho_0^{1+2s}}\bigg)^2\bigg).
\end{align*}
\end{theorem}

\subsection{Noisy observations} \label{section:noisy}

The results above hold for a continuous Gaussian process that is observed exactly. It is thus natural to ask whether similar results hold for discontinuous Gaussian processes or for Gaussian processes observed with errors. In the next proposition, we show that the standard model of discontinuous Gaussian process with a nugget effect yields a zero probability to satisfy bound constraints. Hence, it does not seem possible to define, in a meaningful way, a discontinuous Gaussian process conditioned by bound constraints. 

\begin{proposition} \label{proposition:nugget}
	Let $\mathcal{E}_0$ be defined as in Section \ref{ssec:framework} with $- \infty < \ell $ or $u < + \infty$. Let $Y$ be a Gaussian process on $[0,1]^d$ of the form
	\[
	Y = Y_c + Y_\delta,
	\]
	where $Y_c$ is a continuous Gaussian process on $[0,1]^d$ and $Y_{\delta}$ is a Gaussian process on $[0,1]^d$ with mean function zero and covariance function $\widetilde{k}_{\delta}$ given by
	\[
	\widetilde{k}_{\delta}( u , v) = \delta \mathds{1}_{ \{ u = v \} },
	\]
	for $u,v \in [0,1]^d$. In addition, assume that $Y_c$ and $Y_{\delta}$ are independent. Then
	\[
	\mathbb{P}  ( Y \in \mathcal{E}_0 ) = 0.
	\]
\end{proposition} 

This proposition can be extended to monotonicity and convexity constraints. Hence, in the rest of this section, we consider constrained continuous Gaussian processes observed with noise.

In the case of noisy observations, obtaining fixed-domain asymptotic results on the (unconstrained) MLE of the covariance parameters and the noise variance is challenging, even more so than in the noise-free context. To the best of our knowledge, the only covariance models that have been investigated theoretically, under fixed-domain asymptotics with measurement errors, are the Brownian motion \citep{stein90comparison} and the exponential model \citep{CheSimYin2000,chang2017mixed}. 

In the case of the exponential model, we let $\theta = (\sigma^2 , \rho)$, $\Theta =  [\sigma_l^2 , \sigma_u^2] \times [ \rho_l , \rho_u ]$ with fixed $0 < \sigma_l^2 < \sigma_u^2 < \infty$ and fixed $0 < \rho_l < \rho_u < \infty$. We let $\Delta = [ \delta_l , \delta_u ]$ with $0 < \delta_l < \delta_u < \infty$ being fixed. We consider the set $\{ k_{\sigma^2,\rho} ; (\sigma^2,\rho) \in \Theta \}$ defined by $k_{\sigma^2,\rho}(t) = \sigma^2 e^{ - |t|/\rho }$ for $(\sigma^2,\rho) \in \Theta$ and $t \in \mathbb{R}$. We let $\widetilde{k}_{\sigma^2,\rho}(u,v) = k_{\sigma^2,\rho}(u-v)$ for $u,v \in [0,1]$. We let $Y$ be a Gaussian process on $[0,1]$ with mean function zero and covariance function $k_{\sigma^2_0,\rho_0}$ with $(\sigma^2_0,\rho_0) \in \Theta$. We consider the triangular array of observation points defined by 
\begin{equation*} 
(x_1,\ldots,x_n) = (0,1/(n-1),\ldots,1), 
\end{equation*}
for $n\geqslant 2$. We consider that the $n$ observations are given by 
\[
y_i = Y(x_i) +\epsilon_i,
\]
for $i=1,\ldots,n$ where $\epsilon_1,\ldots,\epsilon_n$ are independent, independent of $Y$, and follow the $\mathcal{N}(0,\delta_0^2)$ distribution. Then the log-likelihood is 
\begin{align*}
\mathcal L_n(\sigma^2,\rho,\delta)=-\frac n2 \ln (2\pi)-\frac 12\ln (|R_{\sigma^2,\rho,\delta}|)-\frac{1}{2} y^{\top} R_{\sigma^2,\rho,\delta}^{-1} y,
\end{align*}
with $(\sigma^2,\rho)  \in \Theta$, $\delta \in \Delta$, $y = (y_1,\ldots,y_n)^{\top}$ and $R_{\sigma^2,\rho,\delta} = [\widetilde{k}_{\sigma^2,\rho}(x_i,x_j)]_{1 \leqslant i,j \leqslant n} + \delta^2 I_n$.
The MLE is given by
\begin{align*}
(\widehat \sigma_{n}^2, \widehat \rho_n , \widehat \delta_n^2) \in \underset{(\sigma^2,\rho) \in \Theta, \delta \in \Delta}{\argmax \; } \mathcal L_{n}(\sigma^2,\rho,\delta).
\end{align*}

In \citep{CheSimYin2000}, it is shown that the MLE $\widehat \sigma_{n}^2/\widehat \rho_n$ of the microergodic parameter and the MLE $\widehat \delta_n^2$ of the noise variance jointly satisfy the central limit theorem
\begin{equation} \label{eq:TCL:noisy}
\begin{pmatrix}
n^{1/4}
\left(
\widehat \sigma_{n}^2/\widehat \rho_n
-
\sigma_0^2 / \rho_0
\right)
\\
n^{1/2}
\left(
\widehat \delta_n^2
-
\delta_0^2
\right)
\end{pmatrix}
\xrightarrow[n\to+\infty]{\mathcal{L}  }
\mathcal{N} 
\left(
\begin{pmatrix}
0 \\ 0
\end{pmatrix}
,
\begin{pmatrix}
4 \sqrt{2} \delta_0 ( \sigma_0^2 / \rho_0 )^{3/2}
& 
0 \\
0 
&
2 \delta_0^4
\end{pmatrix}
\right).
\end{equation}
Hence, the rate of convergence of the MLE of the microergodic parameter is decreased from $n^{1/2}$ to $n^{1/4}$, because of the measurement errors. The rate of convergence of the MLE of the noise variance is $n^{1/2}$. 

In the next proposition, we show that these rates are unchanged when conditioning by the boundedness event $\{ Y \in \mathcal{E}_0 \}$.

\begin{proposition} \label{prop:MLE:noisy}
	Consider the setting defined above, with $\theta_0,\delta_0$ in the interior of $\Theta \times \Delta$. Then, as $n \to \infty$,
	\[
	n^{1/4}
	\left(
	\widehat \sigma_{n}^2/\widehat \rho_n
	-
	\sigma_0^2 / \rho_0
	\right)
	= O_{\P | Y \in \mathcal{E}_0}
	(1),
	~ ~ ~ ~
	n^{1/4}
	\left(
	\widehat \sigma_{n}^2/\widehat \rho_n
	-
	\sigma_0^2 / \rho_0
	\right)
	\neq o_{\P | Y \in \mathcal{E}_0}
	(1),
	\]
	and
	\[
	n^{1/2}
	\left(
	\widehat \delta_n^2
	-
	\delta_0^2
	\right)
	= O_{\P | Y \in \mathcal{E}_0}
	(1),
	~ ~ ~ ~
	n^{1/2}
	\left(
	\widehat \delta_n^2
	-
	\delta_0^2
	\right)
	\neq o_{\P | Y \in \mathcal{E}_0}
	(1).
	\]
\end{proposition}

It would be interesting to see whether the central limit theorem in \eqref{eq:TCL:noisy} still holds conditionally to $\{ Y \in \mathcal{E}_0 \}$. This would be an extension to the noisy case of Theorems \ref{th:cond_after_var} and \ref{th:cond_after_matern}. Nevertheless, to prove Theorem \ref{th:cond_after_var}, we have observed that, in the noiseless case, the MLE of $\sigma_0^2 $ is a normalized sum of the independent variables $W_{n,1}^2 ,\ldots,W_{n,n}^2 $, with 
\[
W_{n,i}  \defeq 
\frac{y_i-\E[y_i|y_1,\ldots,y_{i-1}]}{ \sqrt{ \Var(y_i|y_1,\ldots,y_{i-1}) }},
\]
for $i=1 , \ldots ,n$. We have taken advantage of the fact that conditioning by $W_{n,1} , \ldots, W_{n,k}$ enables to condition by $Y(x_1),\ldots,Y(x_k)$ and to approximately condition by the event $\{ Y \in \mathcal{E}_0 \} $, while leaving the distribution of $W_{n,k+1} , \ldots, W_{n,n}$ unchanged. We refer to the proof of Theorem \ref{th:cond_after_var} for more details.

In contrast, in the noisy case, the authors of \citep{CheSimYin2000} show that  the MLE of $\sigma_0^2 / \rho_0 $ is also a normalized sum of the independent variables $W_{n,1}^2 ,\ldots,W_{n,n}^2 $, but each $W_{n,i}$ depends on the observation vector $y=(y_1,\ldots,y_n)$ entirely \citep[see][Equations (3.40) and (3.42)]{CheSimYin2000}. Hence, it appears significantly more challenging to address the asymptotic normality of the MLE of $\sigma_0^2 / \rho_0 $ and $\delta_0$, conditionally to $\{ Y \in \mathcal{E}_0 \} $. We leave this question open to future work.

The constrained likelihood and the cMLE can be naturally extended to the noisy case. Nevertheless, the asymptotic study of the cMLE, in the context of the exponential covariance function as in Proposition \ref{prop:MLE:noisy}, seems to require substantial additional work. Indeed, to analyze the cMLE in the noiseless case for the Mat\'ern covariance functions, we have relied on the results of \citep{ShaKau2013} and \citep{WanLoh2011}, that are specific to the noiseless case. Furthermore, the martingale arguments, used for instance in the point \tbf{4)} in the proof of Theorem \ref{th:cond_before_var}, require the observation points to be taken from a sequence. Hence, these martingale arguments are not available in the framework of this section, in which the observation points are taken on regular grids. Finally, the RKHS arguments, used for instance in the point \tbf{4)} in the proof of Theorem \ref{th:cond_before_matern_ii}, require to work with covariance functions that are at least twice differentiable, which is not the case with the exponential covariance functions. Hence, we leave the asymptotic study of the cMLE, in the noisy case, open to future research.

\section{Concluding remarks} \label{sec:ccl}

We have shown that the MLE and the cMLE are asymptotically Gaussian distributed, conditionally to the fact that the Gaussian process satisfies either boundedness, monotonicity or convexity constraints. Their asymptotic distributions are identical to the unconditional asymptotic distribution of the MLE. In simulations, we confirm that the MLE and the cMLE have very similar performances when the number $n$ of observation points becomes large enough. We also observe that the cMLE is more accurate for small or moderate values of $n$.

Hence, since the computation of the cMLE is more challenging than that of the MLE, we recommend to use the MLE for large data sets and the cMLE for smaller ones. In the proofs of the asymptotic behavior of the cMLE, one of the main steps is to show that $\mathbb{P}_{\theta} ( Y \in \mathcal E_{\kappa} | y)$ converges to one as $n$ goes to infinity. Hence, in practice, one may evaluate $\mathbb{P}_{\theta} ( Y \in \mathcal E_{\kappa} | y)$, for some values of $\theta$, in order to gauge whether this conditional probability is not too close to $1$ so that it is worth using the cMLE despite the additional computational cost. Similarly, Proposition \ref{proposition:prediction} (and its proof) show that if $\mathbb{P}_{\theta} ( Y \in \mathcal E_{\kappa} | y)$ is close to one, then it is approximately identical to predict new values of $Y$ with accounting for the constraints or not. The latter option is then preferable, as it is computationally less costly.

Our theoretical results could be extended in different ways. First, we remark that the techniques we have used to show that $\mathbb{P}_{\theta} ( Y \in \mathcal E_{\kappa} )$ and $\mathbb{P}_{\theta} ( Y \in \mathcal E_{\kappa} | y )$ are asymptotically negligible (see \eqref{ass:terme1} and \eqref{ass:terme2}) can be used for more general families of covariance functions. Hence, other results on the (unconditional) asymptotic distribution of the MLE could be extended to the case of constrained Gaussian processes in future work. These types of results exist for instance for the product exponential covariance function \citep{ying93maximum}.

Also, in practice, computing the cMLE requires a discretization of the constraints, for instance using a piecewise affine interpolation as in Section \ref{sec:num}, or a finite set of constrained points \citep{DaVeiga2012GPineqconst}. Thus it would be interesting to extend our results by taking this discretization into account. 

Finally, in this paper, we have focused on Gaussian processes that are either observed directly or with an additive Gaussian noise. These contexts are relevant in practice when applying Gaussian processes to computer experiments \citep{santner03design} and to regression problems in machine learning \citep{rasmussen06gaussian}. Nowadays, it has also become standard to study other more complex models of latent Gaussian processes, for instance in Gaussian process classification \citep{rasmussen06gaussian,nickisch2008approximations}. Some authors have also considered latent Gaussian processes subjected to inequality constraints for modelling point processes \citep{lopez2019gaussian}. It would be interesting to obtain asymptotic results similar to those in our article, for latent Gaussian processes. This could be a challenging problem, as few asymptotic results are available even for unconstrained latent Gaussian process models.
We remark that some of the techniques we have used in this paper could be useful when considering latent Gaussian processes under constraints. These techniques are, in particular, Lemmas \ref{lem:discret} and \ref{lem:c_to_hat} and their applications.

\noi\\
\textbf{Acknowledgements}

The authors are indebted to the anonymous reviewers for their helpful comments and suggestions, that lead to an improvement of the manuscript. This research was conducted within the frame of the Chair in Applied Mathematics OQUAIDO, gathering partners in technological research (BRGM, CEA, IFPEN, IRSN, Safran, Storengy) and academia (CNRS, Ecole Centrale de Lyon, Mines Saint-\'Etienne, University of Grenoble, University of Nice, University of Toulouse) around advanced methods for Computer Experiments.\\

\appendix

\section{Additional notation and intermediate results}\label{sec:intermediate}

For $a > 0$, let $f_a : (0,\infty) \to \mathbb{R}$ be defined by $f_a(t) = - \ln(t) - a/t$. We will repeatedly use the fact that $f_a$ has a unique global maximum at $a$ and  $f_a''(t) = 1/t^2 - 2a / t^3$. In addition, let $\xi_{*} = {\inf}_{x \in [0,1]^d } \ \xi(x)$, $\xi^{*} = {\sup}_{x \in [0,1]^d } \ \xi(x) $, and $\xi^{**} = {\sup}_{x \in [0,1]^d } | \xi(x) |$ for any stochastic process $\xi : [0,1]^d \to \mathbb{R}$. 

Now we establish three lemmas that will be useful in the sequel.

\begin{lemma}\label{lem:discret}
	Let $(X_n)_{n\in \N}$ be a sequence of  r.v.'s and $(m_{k,n})_{n, k\in \N, \; k \leqslant n}$ and $(M_{k,n})_{n, k\in \N, \; k \leqslant n}$ be two triangular arrays of  r.v.'s. We consider a random vector $(m,M)^{\top}$ such that $m\leqslant m_{k,n}\leqslant M_{k,n} \leqslant M$ for all $k\leqslant n$. We assume that $\P(m=\ell)=\P(M=u)=0$ and $\P(\ell \leqslant m\leqslant M\leqslant u)>0$ for some fixed $\ell$ and $u \in \mathbb{R}$. Moreover, we consider a sequence $(k_n)_{n \in \mathbb{N}}$ so that, $k_n \leqslant n$, $k_n \to_{n \to \infty} \infty$ and 
	\begin{align}\label{ass:discret}
	(m_{k_n,n},M_{k_n,n})^{\top} \xrightarrow[n \to +\infty]{\textrm{a.s.}} (m,M)^{\top}.
	\end{align}
	Then for any $a\in \R$,
	\begin{align}\label{eq:discret}
	\lim_{n \to+\infty} 
	\Big|\P(X_n\leqslant a | \ell \leqslant m_{k_n,n}\leqslant  M_{k_n,n}\leqslant u)- \P(X_n\leqslant a | \ell \leqslant m \leqslant M \leqslant u)  \Big|=0.
	\end{align}
\end{lemma}

\tbf{Proof of Lemma \ref{lem:discret}}. For the sake of simplicity, we denote by $E_{k,n}$ (respectively $E$) the event $\{\ell \leqslant m_{k,n} \leqslant M_{k,n}\leqslant u\}$ (resp. $\{\ell\leqslant  m\leqslant M\leqslant u\}$). Then 
\begin{align}\label{ineq:discret}
\left|\P(X_n\leqslant a | E_{k_n,n})- \P(X_n\leqslant a | E)  \right| 
\leqslant & \frac{ \left| \P(X_n\leqslant a,\; E_{k_n,n})- \P(X_n\leqslant a,\; E)  \right|}{\P(E_{k_n,n})} \notag \\
& + \left|\frac{1}{\P(E_{k_n,n})} - \frac{1}{\P(E)}\right| \P(X_n\leqslant a,\; E).
\end{align} 

(i) By \eqref{ass:discret}, $\P(E_{k_n,n})$ goes to $\P(E)=\P(\ell \leqslant m \leqslant M\leqslant u)>0$ as $n$ goes to $+\infty$. Thus $1/\P(E_{k_n,n})$ is well-defined for large values of $n$ and bounded as $n \to \infty$. Moreover, by trivial arguments of set theory, one gets 
\begin{align*}
\left| \P(X_n\leqslant a,\; E_{k_n,n})- \P(X_n\leqslant a,\; E)  \right| & 
\leqslant  \P(E_{k_n,n} \Delta E) 
=  \P(E_{k_n,n} \setminus E),
\end{align*} 
since $\P(E \setminus E_{k_n,n})=0$. Now let $\varepsilon>0$. One has  
\begin{align*}
\P&(E_{k_n,n} \setminus E)  = \P(\ell \leqslant m_{k_n,n}\leqslant M_{k_n,n} \leqslant u,\; (m,M) \notin [\ell,u]^2)\\
&\leqslant  \P(\ell \leqslant m_{k_n,n}\leqslant M_{k_n,n} \leqslant u,\; m<\ell)
+ \P(\ell \leqslant m_{k_n,n}\leqslant M_{k_n,n} \leqslant u,\; M>u)\\
&\leqslant  \P(\ell \leqslant m_{k_n,n},\; m<\ell)
+ \P(M_{k_n,n} \leqslant u,\; M>u).
\end{align*}
One may decompose $\P(\ell \leqslant m_{k_n,n},\; m<\ell)$ into
\begin{align*}
\P(\ell+\varepsilon \leqslant m_{k_n,n},\; m<l ) + \P(\ell \leqslant m_{k_n,n} \leqslant \ell+\varepsilon,\; m<l ) 
\leqslant \P(|m_{k_n,n}-m|> \varepsilon) + \P(\ell \leqslant m_{k_n,n} \leqslant \ell+\varepsilon).
\end{align*}
The first term in the right hand-side goes to $0$ as $n$ goes to infinity. 
By Portemanteau's lemma and \eqref{ass:discret},
\[
\displaystyle\limsup_{n \to +\infty }  \P(\ell \leqslant m_{k_n,n} \leqslant \ell+\varepsilon)\leqslant \P(\ell \leqslant m \leqslant \ell +\varepsilon) \underset{\varepsilon\to 0}{\longrightarrow} 0.
\]
We handle similarly the term $\P(M_{k_n,n} \leqslant u,\; M>u)$.
Hence, in the r.h.s. of \eqref{ineq:discret}, the first term goes to $0$ as $n \to \infty$. 

(ii) Now we turn to the control of the second term in \eqref{ineq:discret}. Upper bounding $\P(X_n\leqslant a,\; E)$ by 1, it remains to control $\Big|\frac{1}{\P(E_{k_n,n})} - \frac{1}{\P(E)}\Big|$ which is immediate by the convergence in distribution of $(m_{k_n,n},M_{k_n,n})^{\top}$ as $n$ goes to infinity (implied by the a.s. convergence) and the fact that $\P(E)>0$ and $\P(m=\ell)=\P(M=u)=0$. 
The proof is now complete.
\hfill $\square$

\begin{lemma}\label{lem:c_to_hat}
	Consider three sequences of random functions $f_n,g_n,h_n : [x_{inf}, x_{sup}] \to \mathbb{R}$, with $0 < x_{inf} < x_{sup} < \infty$ fixed.
	Consider that for all $x \in [x_{inf}, x_{sup}]$, $f_n(x)$, $g_n(x)$, and $h_n(x)$ are functions of $Y$ and $x$ only. 
	Let 
	\[
	\widehat{x}_n \in \underset{x \in [x_{inf}, x_{sup}] }{\argmax \; }
	f_n(x).
	\]
	Assume the following properties.
	\begin{enumerate}
		\item[(i)] There exists $A>0$, $B>0$ and $\delta>0$ such that
		\begin{align}\label{ass:vois}
		\displaystyle  f_n(x) - f_n( \widehat{x}_n ) \leqslant -An(x - \widehat{x}_n)^2,\ \forall x \in  [x_{inf}, x_{sup}];\  |x - \widehat{x}_n|\leqslant \delta,
		\end{align}
		and
		\begin{align}\label{ass:loin}
		\displaystyle \underset{ \substack{ | x - \widehat{x}_n | > \delta  \\ x \in  [x_{inf}, x_{sup}] }}{\sup}   f_{n}(x) - f_n  ( \widehat{x}_n)  \leqslant -Bn,
		\end{align}
		with probability going to $1$ as $n \to \infty$. 
		
		\item[(ii)] There exists $C>0$ such that for all $x_1 ,x_2 \in  [ x_{inf},x_{sup} ]$
		\begin{align}\label{ass:terme1}
		\displaystyle \left\vert
		g_n(x_1) - g_n(x_2)
		\right\vert 
		\leqslant C|x_1 - x_2 |,
		\end{align}
		with probability going to $1$ as $n \to \infty$. 
		\item[(iii)] One has, for $\kappa = 0,1,2$,
		\begin{align}\label{ass:terme2}
		\displaystyle \underset{
			x_1 ,x_2 \in  [ x_{inf},x_{sup} ]
		}  {\sup}
		\left\vert 
		h_n(x_1) - h_n(x_2)
		\right\vert =o_{\P | Y\in \mathcal E_\kappa}(1).
		\end{align}
	\end{enumerate}
	Then, with
	\begin{align*}
	\widehat{\widehat{x}}_n \in \underset{x \in [x_{inf}, x_{sup}] }{\argmax \; }
	\{ f_n(x) + g_n(x) + h_n(x) \},
	\end{align*}
	we have
	\begin{align}\label{eq:c_to_hat}
	\sqrt n
	|
	\widehat{\widehat{x}}_n
	-
	\widehat{x}_n
	|=o_{\P | Y\in \mathcal E_\kappa}\left( 1 \right).
	\end{align}
\end{lemma}

\tbf{Proof of Lemma \ref{lem:c_to_hat}}. Let $\varepsilon >0$. First, we have, with probability (conditionally to $\{Y\in \mathcal E_\kappa\}$)   going to 1 as $n \to \infty$, from \eqref{ass:vois}, \eqref{ass:terme1} and \eqref{ass:terme2}
\begin{align*}
\underset{ 
	\substack{ 
		| x - \widehat{x}_n | \geqslant {\varepsilon}/{\sqrt{n}}
		\\
		| x - \widehat{x}_n | \leqslant {1}/{n^{1/4}}
	}
}{\sup} &
\left(
f_n( x )
+
g_n( x )
+
h_n( x )
-
f_n( \widehat{x}_n )
-
g_n( \widehat{x}_n )
-
h_n( \widehat{x}_n )
\right)
& \\
& \leqslant 
- A n 
\left( \frac{ \varepsilon }{ \sqrt{n} } \right)^2
+ \frac{C}{n^{1/4}}
+o_{\mathbb{P} Y\in \mathcal E_\kappa }(1) 
=
- A \varepsilon^2 + o_{\mathbb{P} | Y \in \mathcal E_\kappa }(1). 
\end{align*}

Second, from \eqref{ass:vois}, \eqref{ass:terme1} and \eqref{ass:terme2}, we have, with probability (conditionally to $\{Y\in \mathcal E_\kappa\}$)   going to 1 as $n \to \infty$,
\begin{align*}
\underset{ 
	\substack{ 
		| x - \widehat{x}_n | \geqslant 1/n^{1/4}
		\\
		| x - \widehat{x}_n | \leqslant \delta
	}
}{\sup}
&
\left(
f_n( x )
+
g_n( x )
+
h_n( x )
-
f_n( \widehat{x}_n )
-
g_n( \widehat{x}_n )
-
h_n( \widehat{x}_n )
\right)
& \\
& \leqslant 
- A n 
\left( \frac{ 1 }{ n^{1/4} } \right)^2
+ C \delta
+o_{\mathbb{P}  | Y \in \mathcal E_\kappa}(1) 
\underset{n\to\infty}{\longrightarrow}  - \infty. 
\end{align*}

Third, from \eqref{ass:loin}, \eqref{ass:terme1} and \eqref{ass:terme2}, we have, with probability (conditionally to $\{Y\in \mathcal E_\kappa\}$)   going to 1 as $n \to \infty$,
\begin{align*}
\underset{ 
	| x - \widehat{x}_n | \geqslant  \delta
}{\sup}
&
\left(
f_n( x )
+
g_n( x )
+
h_n( x )
-
f_n( \widehat{x}_n )
-
g_n( \widehat{x}_n )
-
h_n( \widehat{x}_n )
\right)
& \\
& \leqslant 
- B n 
+ C ( x_{sup} - x_{inf} )
+o_{\mathbb{P}  | Y \in \mathcal E_\kappa}(1)
\underset{n\to\infty}{\longrightarrow}  - \infty. 
\end{align*}

Finally, for all $\varepsilon >0$ there exists $c >0$ so that, with probability (conditionally to $\{Y\in \mathcal E_\kappa\}$)    going to $1$ as $n \to \infty$,
\[
\underset{ 
	| x - \widehat{x}_n | \geqslant \varepsilon/\sqrt{n}
}{\sup}
\left(
f_n( x )
+
g_n( x )
+
h_n( x )
-
f_n( \widehat{x}_n )
-
g_n( \widehat{x}_n )
-
h_n( \widehat{x}_n )
\right)
\leqslant - c.
\]
Hence, we have, by definition of $\widehat{\widehat{x}}_n$ 
\begin{align*}
\sqrt n
|
\widehat{\widehat{x}}_n
-
\widehat{x}_n
|=o_{\P  | Y \in \mathcal E_\kappa}\left( 1 \right).
\end{align*}
\hfill $\square$

\begin{lemma} \label{lem:mean:value:sup:GP}
	Let $\{ k_{\theta} ; \theta \in \Theta \}$ be the set of functions in Section \ref{sec:framework} where $\Theta$ is compact. Assume that $k_{\theta}$ satisfies Condition-Var in the case $\kappa = 0$, where $c$ and $\alpha$ can be chosen independently of $\theta$.
	Let $Z_{n,\theta}$ be a Gaussian process with mean function zero and covariance function $(x_1,x_2) \mapsto \Cov_{\theta} ( Y(x_1) , Y(x_2) | y) $.
	Then, we have
	\[
	\underset{ \theta \in \Theta}{\sup}
	\;
	\mathbb{E} \bigg[
	\underset{x \in [0,1]^d}{\sup}
	| Z_{n,\theta} (x)|
	\bigg]
	\underset{n\to\infty}{\to} 0.
	\]
\end{lemma}

\tbf{Proof of Lemma \ref{lem:mean:value:sup:GP}}. This result is proved as an intermediate result in the proof of \citep[Lemma A.3]{LLBDR17}. There, the result was for fixed $\theta$, but it can be made uniform over 
$\theta\in \Theta$ with no additional difficulties. 
\hfill $\square$\\

\section{Proofs for Sections \ref{sec:var} and \ref{sec:matern} - Boundedness}\label{sec:proof:bounded}

We let $\kappa = 0$ throughout Section \ref{sec:proof:bounded}.

\subsection{Estimation of the variance parameter}\label{ssec:proof:bounded_var}

\tbf{Proof of Theorem \ref{th:cond_after_var} under boundedness constraints}.

\tbf{1)} Let $m_{k,n} = {\min}_{i=1,\ldots,k} \ y_i$, $M_{k,n} = {\max}_{i=1,\ldots,k} \ y_i$, and $(m,M)^{\top} = (Y_*,Y^*)^{\top}$, where  $Y_*$ and $Y^*$ have been defined in Appendix \ref{sec:intermediate}. 
We clearly have $m\leqslant m_{k_n,n}\leqslant M_{k_n,n} \leqslant M$.
Since $(x_i)_{i \in \mathbb{N}}$ is dense, for any sequence $(k_n)_{n \in \mathbb{N}}$ so that $k_n \to \infty$ as $n \to \infty$ and $k_n \leqslant n$, we have $(m_{k_n,n}, M_{k_n,n})^{\top} \to (m,M)^{\top}$ a.s. as $n \to \infty$ (up to re-indexing $x_1, \ldots , x_n$).\\

\tbf{2)}  Let $k \in \mathbb{N}$ be fixed. We have
\begin{align*}
\sqrt n\left(\bar{\sigma}_{n}^2-\sigma_0^2\right) 
= & \frac{1}{\sqrt n} \left(y^{\top} R_1^{-1}y-n\sigma_0^2\right).
\end{align*}
Writing the Gaussian probability density function of $y$ as the product of the conditional probability density functions of $y_i$  given $y_1,\ldots,y_{i-1}$ leads to
\begin{align*}
\frac{1}{\sqrt n}\left(y^{\top} R_1^{-1}y-n\sigma_0^2\right)
&= \frac{\sigma_0^2}{\sqrt n} \sum_{i=1}^n \left(\frac{(y_i-\E[y_i|y_1,\ldots,y_{i-1}])^2}{\Var(y_i|y_1,\ldots,y_{i-1})}-1\right).
\end{align*}
The terms in the sum above are independent. Indeed,
\[
\Cov(y_l, y_i-\E[y_i|y_1,\ldots,y_{i-1}])=0, \quad \textrm{for any} \quad l\leqslant i-1
\]
and the Gaussian distribution then leads to independence. Therefore,
\begin{align*}
\frac{1}{\sqrt n}\left(y^{\top} R_1^{-1}y-n\sigma_0^2\right)
= \frac{\sigma_0^2}{\sqrt n} \sum_{i=1}^k \left(\frac{(y_i-\E[y_i|y_1,\ldots,y_{i-1}])^2}{\Var(y_i|y_1,\ldots,y_{i-1})}-1\right)
+\frac{\sigma_0^2}{\sqrt n} \sum_{i=k+1}^n \left(\frac{(y_i-\E[y_i|y_1,\ldots,y_{i-1}])^2}{\Var(y_i|y_1,\ldots,y_{i-1})}-1\right).
\end{align*}
The first term is $o_{\P}(1)$ being the sum of $k$ r.v.'s (whose variances are all equal to $2$) divided by the square root of $n$. Because $\P_{\sigma^2}\Big(\ell \leqslant \underset{i=1,\ldots,k}{\min} y_i\leqslant\underset{i=1,\ldots,k}{\max} y_i \leqslant u\Big)>0$, the first term is also $o_{\P}(1)$ conditionally to $\Big\{\ell \leqslant \underset{i=1,\ldots,k}{\min} y_i\leqslant\underset{i=1,\ldots,k}{\max} y_i \leqslant u\Big\}$.
The second term is equal to $\sigma_0^2 / \sqrt{n}$ times the sum of $n-k$   independent variables with zero mean  and variance 2 and is also independent of $y_1,\ldots,y_k$. Hence, from the central limit theorem and Slutsky's lemma \citep[Lemma 2.8]{van2000asymptotic}, we obtain that 
\[
\frac{1}{\sqrt n}\left(y^{\top} R_1^{-1}y-n\sigma_0^2\right)
\xrightarrow[n \to \infty]{ \mathcal{L} | y \in \mathcal E_{0,k} } \mathcal{N}( 0 , 2 \sigma_0^4 ),
\]
where $\mathcal E_{0,k}\defeq \Big\{y:\; \ell \leqslant \underset{i=1,\ldots,k}{\min} y_i \leqslant \underset{i=1,\ldots,k}{\max} y_i \leqslant u\Big\}$ and $\xrightarrow[n \to \infty]{ \mathcal{L} | y\in \mathcal E_{0,k} }$ is defined similarly as $\xrightarrow[n \to \infty]{ \mathcal{L} |  Y\in \mathcal E_0 }$. \\

\tbf{3)} Hence, for $x \in \mathbb{R}$, there exists a sequence $\tau_n \underset{n \to \infty}{\longrightarrow} \infty$ satisfying $\tau_n = o(n)$ as $n \to \infty$ so that:
\[
\mathbb{P}
\left(
\Big.
\sqrt{n}
\left(
\bar{\sigma}_n^2 - \sigma_0^2
\right)
\leqslant x
\Big|
\ell \leqslant \underset{i=1,\ldots,\tau_n}{\min} y_i\leqslant \underset{i=1,\ldots,\tau_n}{\max} y_i \leqslant u
\right)
\underset{n \to \infty}{\longrightarrow}
\mathbb{P}
\left(
V \leqslant x
\right),
\]
with $V \sim \mathcal{N}( 0, 2 \sigma_0^4 )$. The above display naturally holds. Indeed, if $(S_{\tau,n})_{n \in \mathbb{N},\tau=1,\ldots,n}$ is a triangular array of numbers so that, for any fixed $\tau$, $S_{\tau,n} \to S$ as $n \to \infty$, where $S$ does not depend on $\tau$, then there exists a sequence $\tau_n \to \infty$ so that $S_{\tau_n,n} \to S$ as $n \to \infty$.

Therefore, from Lemma \ref{lem:discret},
\[
\mathbb{P}
\left(
\left.
\sqrt{n}
\left(
\bar{\sigma}_n^2 - \sigma_0^2
\right)
\leqslant x
\right|
\ell \leqslant Y(x) \leqslant u,\; \forall x\in [0,1]^d
\right)
\underset{n \to \infty}{\longrightarrow}
\mathbb{P}
\left(
V \leqslant x
\right).
\]
This concludes the proof.
\hfill $\square$ \\


\tbf{Proof of Theorem \ref{th:cond_before_var}  under boundedness constraints}. 
We apply Lemma \ref{lem:c_to_hat} to the sequences of functions $f_n$, $g_n$ and $h_n$ defined by $f_n( \sigma^2 ) = \mathcal{L}_n( \sigma^2  )$, $g_n(x) = A_n( \sigma^2)$, and $h_n( \sigma^2 ) = B_n( \sigma^2 )$. Here we recall that for $\sigma^2\in \Theta$,
\[
A_n( \sigma^2  )
= -\ln \P_{\sigma^2}
\left( Y\in\mathcal E_0 \right) 
\;
\mbox{and}
\;
B_n( \sigma^2 )
= \ln \P_{\sigma^2 }
\left( \left. Y\in\mathcal E_0 \right| y \right) .
\]

In order to apply Lemma \ref{lem:c_to_hat}, we need to check that the conditions \eqref{ass:vois} to \eqref{ass:terme2} hold.\\

\tbf{1)} By \eqref{def:log_likeli_var}, one has
\[
\mathcal L_n(\sigma^2)=-\frac n2 \ln 2\pi-\frac n2  \ln (\sigma^2)-\frac 12\ln (|R_1|)-\frac{1}{2\sigma^2} y^{\top} R_1^{-1} y.
\]

Now $y^{\top} R_1^{-1} y$ is the square of the norm of a Gaussian vector with variance-covariance matrix $\sigma_0^2I_n$, where  $I_n$ stands for the identity matrix of dimension $n$. Thus one can write $y^{\top} R_1^{-1} y$ as 
the sum of the squares of $n$ independent and identically distributed r.v.'s $\varepsilon_i$, where $\varepsilon_i$ is Gaussian distributed with mean 0 and variance $\sigma_0^2$.
We prove that \eqref{ass:vois} is satisfied.
One may rewrite $\mathcal{L}_{n}(\sigma^2)$ as
\begin{equation} \label{eq:link:lik:fa:brownian}
\mathcal{L}_{n}(\sigma^2)
=
- \frac{n}{2} \ln( 2 \pi)
- \frac{1}{2} \ln(|R_1|)
+ \frac{n}{2} f_{\sigma_0^2 + o_{\mathbb{P}}(1)} (\sigma^2),
\end{equation}
where the $o_{\mathbb{P}}(1)$ above does not depend on $\sigma^2$ and $f_a$ has been introduced in Appendix \ref{sec:intermediate}. 
By a Taylor expansion and the definition of $\bar \sigma_n^2$, we have, with probability going to 1 as $n \to \infty$,
\begin{align*}
\mathcal{L}_{n}(\sigma^2) -  \mathcal{L}_{n}(\bar{\sigma}_n^2) & = (\sigma^2-\bar{\sigma}_n^2)   \mathcal{L}_{n}'(\bar{\sigma}_n^2) 
+\frac 12(\sigma^2-\bar{\sigma}_n^2)^2 \mathcal{L}_{n}''(\widetilde{\sigma}^2)\\
&= 
\frac{n}{4} f''_{\sigma_0^2 + o_{\mathbb{P}}(1)}(\widetilde{\sigma}^2)
(\sigma^2-\bar{\sigma}_n^2)^2  \\
&= 
\frac{n}{4} 
\left(
\frac{1}{ \widetilde{\sigma}^4 }
- 
2
\frac{ \sigma_0^2 + o_{\mathbb{P}}(1) }{ \widetilde{\sigma}^6 }
\right)
(\sigma^2-\bar{\sigma}_n^2)^2, 
\end{align*}
with $\widetilde{\sigma}^2$ in the interval with endpoints $\sigma^2$ and $\bar{\sigma}_n^2$. Hence, non-random constants $A >0 $ and $\delta>0$ exist for which \eqref{ass:vois} is satisfied. \\

\tbf{2)} Second, let us prove that \eqref{ass:loin} holds with the previous $\delta >0$ and for some $B >0$. From \eqref{eq:link:lik:fa:brownian}, $2\mathcal L_n/n+\ln(2\pi)+(1/n) \ln(|R_1|)$ converges uniformly on $[\sigma_l^2,\sigma_u^2]$ as $n$ goes to infinity to $f_{\sigma_0^2}$. The function $f_{\sigma_0^2}$ attains its unique maximum at $\sigma_0^2$, which implies the result since $\bar{\sigma}_n^2$ converges to $\sigma_0^2$ in probability. Hence \eqref{ass:loin} holds.\\

\tbf{3)} Now we consider \eqref{ass:terme1}.
Let us introduce the Gaussian process $Y_r$ with mean function zero and covariance function $\widetilde k_1$.
Let $\sigma_1^2 \leqslant \sigma_2^2$. Then, one has:
\begin{align*}
\Big|\exp\left\{ -A_n(\sigma_1^2 ) \right\}
& -
\exp\left\{ -A_n(\sigma_2^2 ) \right\}\Big|
= \Big|\P (\sigma_1 Y_r \in \mathcal E_0)- \P (\sigma_2 Y_r \in \mathcal E_0)\Big|\\
&\leqslant  \P \left(\frac{u}{\sigma_2}\leqslant  Y_r(x) \leqslant \frac{u}{\sigma_1},\; \forall x\in [0,1]^d\right) +  \P \left(\frac{\ell}{\sigma_2}\leqslant  Y_r(x) \leqslant \frac{\ell}{\sigma_1},\; \forall x\in [0,1]^d\right) \\
& \leqslant c\left|\frac{1}{\sigma_1}-\frac{1}{\sigma_2}\right|
\leqslant c |\sigma_2^2-\sigma_1^2|,
\end{align*}
by Tsirelson theorem in \citep{AW09}. 
Then, from Lemma \ref{lem:lower:bounded:proba:bounded:GP_var}, \eqref{ass:terme1} holds.\\

\tbf{4)} We turn to 
\[
B_n(\sigma^2) = \ln \P_{\sigma^2}(Y\in \mathcal E_0| y)=\ln \P_{\sigma^2}(\ell \leqslant Y(x)\leqslant u, \; \forall x\in[0,1]^d| y).
\] 
Let $m_{n,y}$ and $\sigma^2 \widetilde k_{n}$ be the conditional mean and covariance functions of $Y$ given $y$, under the probability measure $\P_{\sigma^2}$. Using Borell-TIS inequality \citep{Adler07}, with $Z_{n , \sigma^2}$ a Gaussian process with mean function zero and covariance function $\sigma^2 \widetilde k_{n}$, we obtain
\begin{align}
\P_{\sigma^2}(Y^*> u| y)
&\leqslant
\P_{\sigma^2}\bigg( Z_{n,\sigma^2}^* > u-\underset{x\in [0,1]^d}{\sup} m_{n,y}(x) | y\bigg)\nonumber\\
&\leqslant 
\P_{\sigma^2}\bigg(Z_{n,\sigma^2}^{**}> u-\underset{x\in [0,1]^d}{\sup} m_{n,y}(x) | y\bigg)\nonumber\\
&\leqslant
\exp\left\{
-\frac{
	\bigg( \bigg(u-\underset{x\in [0,1]^d}{\sup} m_{n,y}(x)-\E[Z_{n,\sigma^2}^{**}] \bigg)_+\bigg)^2
}{
2\underset{x\in [0,1]^d}{\sup} \E[Z_{n,\sigma^2}(x)^2]
}
\right\}.\label{ineq:BTIS}
\end{align}

But by Lemma \ref{lem:mean:value:sup:GP}, ${\sup}_{\sigma^2 \in [ \sigma_l^2 , \sigma_u^2 ]} \ \E[Z_{n,\sigma^2}^{**}]\to 0$ as $n\to+\infty$. Additionally, one can simply show that ${\sup}_{x\in [0,1]^d} \ \E[Z_{n,\sigma^2}(x)^2] =
{\sup}_{x\in [0,1]^d} \ \sigma^2 \widetilde k_{n }(x,x)$ goes to zero uniformly in $\sigma^2 \in [ \sigma_l^2 , \sigma_u^2 ]$ as $n \to \infty$. By \citep[Proposition 2.8]{BBG16} and because the sequence of observation points is dense,
\[
\displaystyle \underset{x\in [0,1]^d}{\sup} |m_{n,y}(x)-Y(x)|\xrightarrow[n\to +\infty]{\textrm{a.s.}} 0,
\]
from which we deduce that on $\{Y^* < u-\delta\}$, a.s.
\[
\displaystyle\limsup_{n\to+\infty } \bigg(u-\underset{x\in [0,1]^d}{\sup} m_{n,y}(x)\bigg) \geqslant \delta.
\]
Consequently, \eqref{ineq:BTIS} leads to 
\begin{align}\label{limit:BTIS_sup}
\mathds{1}_{ \{ Y^* < u-\delta \} }
\underset{\sigma^2 \in [ \sigma_l^2 , \sigma_u^2 ]}{\sup} \P_{\sigma^2} (Y^*> u| y) \xrightarrow[n\to +\infty]{\textrm{a.s.}} 0.
\end{align}
Similarly, taking $-Y$ instead of $Y$, one may prove easily that 
\begin{align}\label{limit:BTIS_inf}
\mathds{1}_{ \{ Y_* > l+\delta \} }
\underset{\sigma^2 \in [ \sigma_l^2 , \sigma_u^2 ]}{\sup} \P_{\sigma^2} (Y_*< l| y) \xrightarrow[n\to+\infty]{\textrm{a.s.}} 0.
\end{align}
Then, we deduce that
\begin{align}\label{limit:BTIS}
\mathds{1}_{ \{ \ell+\delta <  Y(x)  < u-\delta,\; \forall x\in [0,1]^d \} }
\underset{\sigma^2 \in [ \sigma_l^2 , \sigma_u^2 ]}{\sup} \P_{\sigma^2} (Y^*> u \; \textrm{or} \; Y_*<\ell| y) \xrightarrow[n\to+\infty]{\textrm{a.s.}} 0.
\end{align}
Now let $\varepsilon>0$, $\varepsilon'=2|\ln(1-\varepsilon)|$ and $\mathcal E_{0,\delta}\defeq \{ f\in \mathcal C([0,1]^d,\R) \quad \textrm{s.t.}   \; \ell+\delta \leqslant f(x) \leqslant u-\delta,\, \forall x\in [0,1]^d \}$. We have:
\begin{align*}
\P \bigg(\underset{\sigma^2 \in [ \sigma_l^2 , \sigma_u^2 ]}{\sup} \P_{\sigma^2} &\left(Y^*> u \; \textrm{or} \; Y_*<\ell| y\right)\geqslant \varepsilon,\; Y\in \mathcal E_{0,\delta}\bigg) \underset{n\to+\infty}{\longrightarrow} 0\\
&= \P\bigg(\underset{\sigma^2 \in [ \sigma_l^2 , \sigma_u^2 ]}{\inf} B_n(\sigma^2) \leqslant - \varepsilon'/2,\;Y\in \mathcal E_{0,\delta}\bigg)\\
&= \P\bigg(\underset{\sigma^2 \in [ \sigma_l^2 , \sigma_u^2 ]}{\sup} |B_n(\sigma^2)| \geqslant \varepsilon'/2,\; Y\in \mathcal E_{0,\delta}\bigg)\\
&\geqslant \P\bigg(\underset{\sigma_1^2,\sigma_2^2 \in [ \sigma_l^2 , \sigma_u^2 ]}{\sup} |B_n(\sigma_1^2)-B_n(\sigma_2^2)| \geqslant \varepsilon',\; Y\in \mathcal E_{0,\delta}\bigg)
\end{align*}
by the triangular inequality and \eqref{limit:BTIS}. Therefore, 
\begin{align}
\P&\bigg(\underset{\sigma_1^2,\sigma_2^2 \in [ \sigma_l^2 , \sigma_u^2 ]}{\sup} |B_n(\sigma_1^2)-B_n(\sigma_2^2)| \geqslant \varepsilon',\; Y\in \mathcal E_0\bigg)\notag\\
=& \P\bigg(\underset{\sigma_1^2,\sigma_2^2 \in [ \sigma_l^2 , \sigma_u^2 ]}{\sup} |B_n(\sigma_1^2)-B_n(\sigma_2^2)| \geqslant \varepsilon',\; Y\in \mathcal E_{0,\delta}\bigg) \label{eq:eq1} \\
& + \P\bigg(\underset{\sigma_1^2,\sigma_2^2 \in [ \sigma_l^2 , \sigma_u^2 ]}{\sup} |B_n(\sigma_1^2)-B_n(\sigma_2^2)| \geqslant \varepsilon',\; Y\in \mathcal E_{0} \setminus \mathcal E_{0,\delta}\bigg).\label{eq:eq2}
\end{align}
As already shown, the term \eqref{eq:eq1} converges to 0 as $n\to+\infty$ for any fixed $\delta>0$. For \eqref{eq:eq2}, we have
\[
\underset{
	\substack{ 
		t_1 , t_2 \in \mathbb{R}
		\\
		t_1 \neq t_2
	}
}{\sup}
\frac{ |  \mathbb{P}_{\sigma_0^2} ( Y^* \leqslant t_1 ) - \mathbb{P}_{\sigma_0^2} ( Y^* \leqslant t_2 )   | }{ | t_1 - t_2 | }
< + \infty.
\]
This follows from Tsirelson theorem in \citep{AW09}. Hence
for all $\varepsilon>0$, there exists $\delta^*>0$ such that, 
\begin{equation} \label{eq:Ystar:no:concentration:u}
\P_{\sigma_0^2}(Y^* \in [u-\delta^*,u])\leqslant \varepsilon.
\end{equation}

Similarly, for all $\varepsilon>0$, there exists $\delta_*>0$ such that, 
\begin{equation} \label{eq:Ystar:no:concentration:l}
\P_{\sigma_0^2}(Y_* \in [\ell+\delta_*,\ell])\leqslant \varepsilon.
\end{equation}
Taking $\delta=\min(\delta_*,\delta^*)$, we conclude the proof of \eqref{ass:terme2}.\\

\tbf{5)} Finally, we remark that with probability going to one as $n \to \infty$, $\bar{\sigma}_n^2={\argmax \; }_{ \sigma^2\in [\sigma_l^2,\sigma_u^2]} \ \mathcal L_n(\sigma^2)$. Hence, one may apply Lemma \ref{lem:c_to_hat} to obtain
\[
\sqrt n |\widehat \sigma_{n,c}^2-\bar{\sigma}_n^2 | =o_{\P | Y\in \mathcal E_0}\left( 1 \right).
\]

By Theorem \ref{th:cond_after_var} and Slutsky's lemma, we conclude the proof. 
\hfill $\square$

\subsection{Isotropic Mat\'ern process}\label{ssec:proof:bounded_matern}

Before proving Theorems \ref{th:cond_after_matern}, \ref{th:cond_before_matern_i} and \ref{th:cond_before_matern_ii}, we establish an intermediate result useful in the sequel.

\begin{lemma} \label{lem:mater:change:rho}
	For $\rho \in [ \rho_l , \rho_u ] $, let
	\[
	\bar{\sigma}_n^2 (\rho)
	\in \underset{ \sigma^2 \in ( 0 , \infty) }{\argmax \; }
	\mathcal{L}_n  ( \sigma^2 , \rho ),
	\quad \mbox{and} \quad
	\widehat{\sigma}_n^2 (\rho)
	\in \underset{ \sigma^2 \in [ \sigma_l^2 , \sigma_u^2 ] }{\argmax \; }
	\mathcal{L}_n  ( \sigma^2 , \rho ).
	\] 
	Then, we have
	\begin{equation} \label{eq:mater:change:rho:1}
	\underset{ \rho_1 , \rho_2 \in [ \rho_l , \rho_u ] }{\sup}
	\left|
	\frac{\bar{\sigma}_n^2 (\rho_1)}{ \rho_1^{2 \nu} }
	-
	\frac{\bar{\sigma}_n^2 (\rho_2)}{ \rho_2^{2 \nu} }
	\right|
	= o_\P( 1 / \sqrt{n} )
	\end{equation}
	and 
	\begin{equation} \label{eq:mater:change:rho:2}
	\underset{ \rho_1 , \rho_2  \in [ \rho_l , \rho_u ] }{\sup}
	\left|
	\frac{\widehat{\sigma}_n^2 (\rho_1)}{ \rho_1^{2 \nu} }
	-
	\frac{\widehat{\sigma}_n^2 (\rho_2)}{ \rho_2^{2 \nu} }
	\right|
	= o_\P( 1 / \sqrt{n} ). 
	\end{equation}
\end{lemma}

\tbf{Proof of Lemma \ref{lem:mater:change:rho}}. We have $\bar{\sigma}_n^2 (\rho_1) /  (\rho_1^{2 \nu}) = (1/n) (1/  (\rho_1^{2 \nu})) y^{\top} R_{\rho_1 , \nu}^{-1} y$ so that from \citep[Lemma 1]{ShaKau2013}, we get $\frac{\bar{\sigma}_n^2 (\rho_1)}{ \rho_1^{2 \nu} } \geqslant \frac{\bar{\sigma}_n^2 (\rho_2)}{ \rho_2^{2 \nu} }$ for $\rho_1 \leqslant \rho_2$. Thus,
\[
\underset{ \rho_1 , \rho_2 \in [ \rho_l , \rho_u ] }{\sup}
\left|
\frac{\bar{\sigma}_n^2 (\rho_1)}{ \rho_1^{2 \nu} }
-
\frac{\bar{\sigma}_n^2 (\rho_2)}{ \rho_2^{2 \nu} }
\right|
=
\left|
\frac{\bar{\sigma}_n^2 (\rho_l)}{ \rho_l^{2 \nu} }
-
\frac{\bar{\sigma}_n^2 (\rho_u)}{ \rho_u^{2 \nu} }
\right|.
\]
Then, it is shown in the proof of \citep[Theorem 3]{WanLoh2011} \citep[see, also,][for its proof]{Wang10} that 
\[
\left|
\frac{\bar{\sigma}_n^2 (\rho_l)}{ \rho_l^{2 \nu} }
-
\frac{\bar{\sigma}_n^2 (\rho_0)}{ \rho_0^{2 \nu} }
\right|
= 
o_\P( 1 / \sqrt{n} ),
\]
and similarly for $\rho_u$. Hence, \eqref{eq:mater:change:rho:1} follows. Also, let $p_{\sigma_l^2, \sigma_u^2}$ be the function from $( 0 , \infty )$ to $[\sigma_l^2, \sigma_u^2]$ defined by $p_{\sigma_l^2, \sigma_u^2}(t) = \min( \max( t,\sigma_l^2 ),\sigma_u^2 )$. Then, since $\sigma^2 \mapsto \mathcal{L}_n( \sigma^2 ,\rho)$ is first increasing and then decreasing, we have $\widehat{\sigma}_n^2 (\rho) = p_{\sigma_l^2, \sigma_u^2}( \bar{\sigma}_n^2(\rho))$. Notice that $p_{\sigma_l^2, \sigma_u^2}$ is continuous and bounded by $\sigma_u^2$. Hence, \eqref{eq:mater:change:rho:2} follows.
\hfill $\square$\\

%
%

\tbf{Proof of Theorem \ref{th:cond_after_matern}  under boundedness constraints}. 
Because $\sigma_l^2/(\rho_l^{2 \nu}) < \sigma_0^2/(\rho_0^{2 \nu}) < \sigma_u^2/(\rho_u^{2 \nu})$, we have $\widehat{\sigma}_n^2 = \widehat{\sigma}_n^2(\widehat{\rho}_n)$ with the notation of Lemma \ref{lem:mater:change:rho}. Also, with probability going to $1$ as $n \to \infty$, $ \bar{\sigma}_{n}^2( \widehat \rho_{n} ) = \widehat{\sigma}_{n}^2( \widehat \rho_{n} ) $ with the notation of Lemma \ref{lem:mater:change:rho}.
From Lemma  \ref{lem:mater:change:rho} and with the notation therein, we have
\begin{align*}
\sqrt n\left(\frac{\widehat \sigma_{n}^2(\widehat{\rho}_n)}{\widehat \rho_{n}^{2\nu}}-\frac{ \sigma_0^2}{ \rho_0^{2\nu}}\right) 
= \sqrt n\left(\frac{\widehat \sigma_{n}^2 ( \rho_0 ) }{ \rho_{0}^{2\nu}}-\frac{ \sigma_0^2}{ \rho_0^{2\nu}}\right)
+ \sqrt n\left(\frac{\widehat \sigma_{n}^2(\widehat{\rho}_n)}{\widehat \rho_{n}^{2\nu}}-\frac{\widehat \sigma_{n}^2( \rho_0 )}{ \rho_0^{2\nu}}\right)
= \frac{1}{\sqrt n}\frac{1}{\rho_0^{2\nu}} \left(y^{\top} R_{\rho_0,\nu}^{-1}y-n\sigma_0^2\right)
+ o_{\P}(1),
\end{align*}
where we have used that, with probability going to $1$ as $n \to \infty$, $\widehat \sigma_{n}^2( \rho_0 ) = \bar{\sigma}_n^2 (\rho_0) $. Then we conclude by applying Theorem \ref{th:cond_after_var} when $\kappa=0$ and with $K_0=K_{\nu}(\cdot{}/\rho_0)$.
\hfill $\square$\\

\tbf{Proof of Theorem \ref{th:cond_before_matern_i}  under boundedness constraints}. We apply Lemma \ref{lem:c_to_hat} to the sequences of functions $f_n$, $g_n$ and $h_n$ defined by $f_n( \sigma^2 ) = \mathcal{L}_n( \sigma^2 , \rho_1 )$, $g_n(x) = A_n( \sigma^2 , \rho_1 )$, and $h_n( \sigma^2 ) = B_n( \sigma^2 , \rho_1)$.\\

\tbf{1)} We have, with $\sigma_1^2$ so that $\sigma_1^2 / \rho_1^{2 \nu} = \sigma_0^2 / \rho_0^{2 \nu}$, 
\begin{align*}
\mathcal L_n(\sigma^2 , \rho_1)
= & 
- \frac{n}{2} \ln( 2 \pi )
-\frac n2  \ln (\sigma^2)-\frac 12 \ln (| R_{\rho_1,\nu}|) -\frac{1}{2\sigma^2} y^{\top} R_{\rho_1,\nu}^{-1} y \\
= &
- \frac{n}{2} \ln( 2 \pi )
-\frac n2  \ln (\sigma^2)
-\frac 12 \ln (| R_{\rho_1,\nu}|) 
-\frac{1}{2\sigma^2}
\left[
\frac{\rho_1^{2 \nu}}{\rho_0^{2 \nu}} y^{\top} R_{\rho_0,\nu}^{-1} y
+ o_\P( \sqrt{n} )
\right],
\end{align*}
from \eqref{eq:mater:change:rho:1} in Lemma \ref{lem:mater:change:rho}, observing that 
$\bar{\sigma}_n^2 (\rho_1) = (1/n) y^{\top} R_{\rho_1,\nu}^{-1} y$. Thus 
\begin{align*}
\mathcal L_n(\sigma^2 , \rho_1)
= 
- \frac{n}{2} \ln( 2 \pi )
-\frac n2 \ln (\sigma^2)
-\frac 12 \ln (| R_{\rho_1,\nu}|) 
-\frac{ \sigma_1^2}{2\sigma^2}
\left[
\frac{1}{\sigma_0^2} y^{\top} R_{\rho_0,\nu}^{-1} y
+ o_\P( \sqrt{n} )
\right],
\end{align*}
where $(1/\sigma_0^2) y^{\top} R_{\rho_0,\nu}^{-1} y$ is a sum of the squares of independent standard Gaussian variables. Hence, we show \eqref{ass:vois} and \eqref{ass:loin} exactly as for the proof of Theorem \ref{th:cond_before_var}  when $\kappa=0$. \\

\tbf{2)} Assumption \eqref{ass:terme1} is satisfied since it has been established in the proof of Theorem \ref{th:cond_before_var}  when $\kappa=0$ (for any $\rho_0 \in (0 , \infty))$ and does not involve $y$.\\

\tbf{3)} We turn to $B_n(\sigma^2 , \rho_1)=\ln \P_{(\sigma^2 , \rho_1)}(Y\in \mathcal E_0| y)$. 
Similarly to the proof of Theorem \ref{th:cond_before_var}  when $\kappa=0$, we have that, for any $\delta >0$:
\begin{align}\label{eq:cv_Bn}
\mathbb{P}_{(\sigma_1^2,\rho_1)}
\bigg(
\mathds{1}_{ \{ \ell+\delta\leqslant Y(x) \leqslant u - \delta,\; \forall x\in [0,1]^d \}}
\underset{\sigma^2\in [\sigma_l^2,\sigma_u^2]}{\sup}
\left|
B_n(\sigma^2,\rho_1)
\right|
\underset{n\to\infty}{\longrightarrow} 0
\bigg)
=1.
\end{align}
Now, for $(\sigma^2,\rho) \in (0 , \infty)^2$, we recall that $\P_{(\sigma^2,\rho)}$ is the measure on $\Omega$ for which $Y : ( \Omega , \mathcal{A} ,\P_{\sigma^2 , \rho} )  \to (\mathcal{C}( [0,1]^d ,\R), \mathcal{B} )$ has the distribution of a Gaussian process on $[0,1]^d$ with mean function zero and covariance function $\widetilde k_{\sigma^2, \rho,\nu}$. 
By \citep[Theorem 2]{zhang04inconsistent}, the measures $\P_{(\sigma_0^2,\rho_0)}$ and  $\P_{(\sigma_1^2,\rho_1)}$ are equivalent as soon as $\sigma_0^2/\rho_0^{2\nu}=\sigma_1^2/\rho_1^{2\nu}$  meaning 
that for any set $A \in \mathcal{A} $,
\[
\P_{(\sigma_0^2,\rho_0)}(A)=1 \Leftrightarrow \P_{(\sigma_1^2,\rho_1)}(A)=1.
\] 
Then, one gets
\[
\mathbb{P}_{(\sigma_0^2,\rho_0)}
\bigg(
\mathds{1}_{ \{ \ell+\delta\leqslant Y(x) \leqslant u - \delta,\; \forall x\in [0,1]^d \}}
\underset{\sigma^2\in [\sigma_l^2,\sigma_u^2]}{\sup}
\left|
B_n(\sigma^2,\rho_1)
\right|
\underset{n\to\infty}{\longrightarrow} 0
\bigg)
=1,
\]
which can also be written as
\[
\mathds{1}_{ \{ \ell+\delta\leqslant Y(x) \leqslant u - \delta,\; \forall x\in [0,1]^d \}}
\underset{\sigma^2\in [\sigma_l^2,\sigma_u^2]}{\sup}
\left|
B_n(\sigma^2,\rho_1)
\right|
\xrightarrow[n \to \infty] {\textrm{a.s.}} 0. 
\]
This implies \eqref{ass:terme2}, as for the proof of Theorem \ref{th:cond_before_var}  when $\kappa=0$.\\

\tbf{4)} From a special case of Theorem \ref{th:cond_after_matern} when $\kappa=0$ and with $\rho_l = \rho_u = \rho_1$, we have 
\[
\sqrt n \left(\widehat \sigma_{n}^2 ( \rho_1 ) -\sigma_1^2 \right) \xrightarrow[n\to+\infty]{\mathcal{L} | Y\in \mathcal E_0} \mathcal N(0,2\sigma_1^4).
\]

Therefore, by Lemma \ref{lem:c_to_hat} and Slutsky's lemma, we conclude the proof of Theorem \ref{th:cond_before_var}  when $\kappa=0$. \hfill $\square$\\

\tbf{Proof of Theorem \ref{th:cond_before_matern_ii}  under boundedness constraints}. 
Let $\kappa=0$ in this proof.
We apply Lemma \ref{lem:c_to_hat} to the sequences of functions $f_n$, $g_n$ and $h_n$ defined by $f_n( x ) = \mathcal{L}_n( x \widehat{\rho}_{n,c}^{2 \nu}  , \widehat{\rho}_{n,c} )$, $g_n(x) = A_n( x \widehat{\rho}_{n,c}^{2 \nu} , \widehat{\rho}_{n,c} )$, and $h_n( x ) = B_n( x \widehat{\rho}_{n,c}^{2 \nu} , \widehat{\rho}_{n,c} )$. \\

\tbf{1)}  We naturally have
\[
\frac{ \widehat{\sigma}_{n,c}^2 }{ \widehat{\rho}_{n,c}^{2 \nu} }
\in \underset{ x \in [ \sigma_l^2 / \widehat{\rho}_{n,c}^{2 \nu}  , \sigma_u^2 / \widehat{\rho}_{n,c}^{2 \nu} ] }{\argmax \; }
\mathcal{L}_{n,c} ( x \widehat{\rho}_{n,c}^{2 \nu}, \widehat{\rho}_{n,c}^{2 \nu} ). 
\]

Also, we have
\begin{align*}
\mathcal L_n( x \widehat{\rho}_{n,c}^{2 \nu} , \widehat{\rho}_{n,c})
= &
- \frac{n}{2} \ln( 2 \pi )
-\frac n2  \ln ( x \widehat{\rho}_{n,c}^{2 \nu}  )-\frac 12 \ln (| R_{ \widehat{\rho}_{n,c},\nu}|) -\frac{1}{2 x \widehat{\rho}_{n,c}^{2 \nu} } y^{\top} R_{  \widehat{\rho}_{n,c} ,\nu}^{-1} y \\
= &
- \frac{n}{2} \ln( 2 \pi )
-\frac n2  \ln ( \widehat{\rho}_{n,c}^{2 \nu} )
-\frac 12 \ln (| R_{ \widehat{\rho}_{n,c},\nu}|)
\frac n2  \ln ( x )
-\frac{1}{2 x}
\left[
\frac{
	y^{\top} R_{\rho_0,\nu}^{-1} y
}{
\rho_0^{2 \nu}
}
+ o_\P( \sqrt{n} )
\right],
\end{align*}
from \eqref{eq:mater:change:rho:1} in Lemma \ref{lem:mater:change:rho}, observing that 
$\bar{\sigma}_n^2 (\rho) = (1/n) y^{\top} R_{\rho,\nu}^{-1} y$. Thus
\begin{align} \label{eq:mat:const:case:three:fa}
\mathcal L_n( x \widehat{\rho}_{n,c}^{2 \nu} , \widehat{\rho}_{n,c})  = 
- \frac{n}{2} \ln( 2 \pi )
-\frac n2  \ln ( \widehat{\rho}_{n,c}^{2 \nu} )
 -\frac 12 \ln (| R_{ \widehat{\rho}_{n,c},\nu}|)
+ \frac{n}{2} f_{ (\sigma_0^2 / \rho_0^{2 \nu}) + V_n + o_\P( 1 / \sqrt n ) }(x),
\end{align}
where $\sqrt{n} V_n$ converges in distribution to a centered Gaussian distribution with variance $ 2 (  \sigma_0^2 / \rho_0^{2 \nu} )^2$. Here $V_n$ and the above $o_\P( 1 / \sqrt{n} ) $ do not depend on $x$. Hence, we show \eqref{ass:vois} and \eqref{ass:loin} exactly as for the proof of Theorem \ref{th:cond_before_var}  when $\kappa=0$. \\

\tbf{2)}  One can also see that $\delta >0$ can be chosen so that $ \sigma_l^2 / \rho_l^{ 2 \nu }+ \delta \leqslant \sigma_0^2 / \rho_0^{ 2 \nu }  \leqslant \sigma_u^2 / \rho_u^{ 2 \nu }- \delta$. 
Furthermore, let $x_{inf} = \sigma_{l}^2 / \rho_l^{2\nu}$ and $x_{sup} = \sigma_{u}^2 / \rho_u^{2\nu}$. Then, from \eqref{eq:mat:const:case:three:fa}, one can show that with 
\begin{align}\label{def:tilde_sigma}
\frac{ \widetilde{\sigma}_{n}^2  ( \widehat \rho_{n,c} )}{ \widehat{\rho}_{n,c}^{2 \nu} }
&\in \underset{ x \in [x_{inf} , x_{sup} ] }{\argmax \; }
\mathcal{L}_{n} ( x \widehat{\rho}_{n,c}^{2 \nu}, \widehat{\rho}_{n,c} ),
\end{align}
we have:
\begin{align}\label{eq:avant}
\frac{ \widetilde{\sigma}_{n}^2 ( \widehat \rho_{n,c} )}{ \widehat{\rho}_{n,c}^{2 \nu} }
=
\frac{ \widehat{\sigma}_{n}^2 (\widehat{\rho}_{n,c})}{ \widehat{\rho}_{n,c}^{2 \nu} },
\end{align}
with probability going to $1$ as $n \to \infty$. It is convenient to introduce $\widetilde{\sigma}_{n}^2 ( \widehat \rho_{n,c} )$ because this yields a non-random optimization domain in 
\eqref{def:tilde_sigma}. 
Hence, from Theorem \ref{th:cond_after_matern}  when $\kappa=0$ and from Lemma \ref{lem:mater:change:rho},
\begin{equation*} 
\sqrt n
\left(\frac{\widetilde{\sigma}_{n}^2 ( \widehat \rho_{n,c} )}{\widehat \rho_{n,c}^{2\nu}}-\frac{ \sigma_0^2}{ \rho_0^{2\nu}}
\right)
\xrightarrow[n\to+\infty]{\mathcal{L}| Y\in \mathcal E_0} \mathcal N\left(0,2\left(\frac{\sigma_0^2}{\rho_0^{2\nu}}\right)^2\right).
\end{equation*}
Let also 
\[
\frac{ \widetilde{\sigma}_{n,c}^2 ( \widehat \rho_{n,c} )}{ \widehat{\rho}_{n,c}^{2 \nu} }
\in \underset{ x \in [x_{inf} , x_{sup} ] }{\argmax \; }
\mathcal{L}_{n,c} ( x \widehat{\rho}_{n,c}^{2 \nu}, \widehat{\rho}_{n,c}^{2 \nu} ).
\]
Then, if we show \eqref{ass:terme1} and \eqref{ass:terme2}, we can show similarly as for \eqref{eq:avant} that: 
\[
\frac{ \widetilde{\sigma}_{n,c}^2 ( \widehat \rho_{n,c} ) }{ \widehat{\rho}_{n,c}^{2 \nu} }
=
\frac{ \widehat{\sigma}_{n,c}^2  }{ \widehat{\rho}_{n,c}^{2 \nu} },
\]
with probability going to $1$ as $n \to \infty$.
Hence, from Lemmas \ref{lem:c_to_hat} and \ref{lem:mater:change:rho} and Slutsky's lemma, we can obtain, if \eqref{ass:terme1} and \eqref{ass:terme2} hold,
\[
\sqrt n
\left(\frac{\widehat{\sigma}_{n,c}^2  }{\widehat \rho_{n,c}^{2\nu}}-\frac{ \sigma_0^2}{ \rho_0^{2\nu}}
\right)
\xrightarrow[n\to+\infty]{\mathcal{L}| Y\in \mathcal E_0} \mathcal N\left(0,2\left(\frac{\sigma_0^2}{\rho_0^{2\nu}}\right)^2\right).
\]
Therefore, in order to conclude the proof, it is sufficient to prove \eqref{ass:terme1} and \eqref{ass:terme2}. \\

\tbf{3)} We turn to \eqref{ass:terme1}. Let $\sigma Y_{\rho}$ be a Gaussian process with mean function zero and covariance function $\widetilde k_{\sigma^2,\rho,\nu}$. Then we have, from Lemma \ref{lem:lower:bounded:proba:bounded:matern},
\begin{align*}
\left|
g_n(x_1)
-
g_n(x_2)
\right|
& \leqslant 
c
\left|
\P ( x_1 \widehat{\rho}_{n,c}^{2 \nu} Y_{ \widehat{\rho}_{n,c} }\in \mathcal E_0 )
-
\P ( x_2 \widehat{\rho}_{n,c}^{2 \nu} Y_{ \widehat{\rho}_{n,c} }\in \mathcal E_0)
\right| \\
& 
\\
&
\leqslant 
c
\left|
x_1 - x_2
\right|
\underset{ \substack{
		\rho \in [\rho_l , \rho_u] \\
		t_1, t_2 \in [ l /  (x_{sup} \rho_u^{2\nu}) , u / (x_{inf} \rho_l^{2\nu}) ] 
		\\
		t_1 \neq  t_2 } 
}{\sup}
\frac{ \left|
	\P (  Y_{\rho}^* \leqslant t_1)
	-
	\P ( Y_{\rho}^* \leqslant t_2 )
	\right|}{ |t_1 - t_2| }
. &
\end{align*}
We introduce the following notation:
\[
F_{\rho}(t)=\P\bigg(\underset{x\in [0,1]^d}{\sup} Y_{\rho}(x)\leqslant t\bigg),
\]
and assume that 
\begin{align*} 
\underset{ \substack{
		\rho \in [\rho_l , \rho_u] \\
		t \in [ l /  (x_{sup} \rho_u^{2\nu}) , u / (x_{inf} \rho_l^{2\nu}) ] 
	} 
}{\sup}
\; \underset{\varepsilon >0}{\sup} \; \frac{F_{\rho}(t+\varepsilon)-F_{\rho}(t)}{\varepsilon}=+\infty.
\end{align*}
Therefore, there exists a sequence $(\rho_k,t_k,\varepsilon_k)_{k\in \N}$ such that 
\begin{align}\label{hyp:abs}
\frac{F_{\rho_k}(t_k+\varepsilon_k)-F_{\rho_k}(t_k)}{\varepsilon_k}\xrightarrow[k\to +\infty] {} +\infty.
\end{align}
We extract from $(\rho_k,t_k,\varepsilon_k)_{k\in \N}$ a subsequence (still denoted $(\rho_k,t_k,\varepsilon_k)_{k\in \N}$) such that $(\rho_k)_k$ is convergent and we denote by $\overline \rho$ its limit.
Let $\Phi$ be the cumulative distribution function of a standard Gaussian random variable. Then by the mean value theorem,
\begin{align*}
\frac{\Phi^{-1}\circ F_{\rho_k}(t_k+\varepsilon_k)-\Phi^{-1}\circ F_{\rho_k}(t_k)}{\varepsilon_k} \geqslant \frac{F_{\rho_k}(t_k+\varepsilon_k)-F_{\rho_k}(t_k)}{\varepsilon_k}\underset{p\in [0,1]}{\inf} \left(\Phi^{-1}\right)'(p)
\xrightarrow[k\to +\infty] {} +\infty 
\end{align*}
noticing that ${\inf}_{p\in [0,1]} \ \left(\Phi^{-1}\right)'(p)>0$ and using \eqref{hyp:abs}.\\

But, using the concavity of $\Phi^{-1}\circ F_{\rho}$ \citep[see][Theorem 10 in Section 11]{lifshits95}, one gets  
\begin{align*}
\frac{\Phi^{-1}\circ F_{\rho_k}(t_k+\varepsilon_k)-\Phi^{-1}\circ F_{\rho_k}(t_k)}{\varepsilon_k} 
& \leqslant 
\Phi^{-1}\circ F_{\rho_k}
\left( \frac{l}{x_{sup} \rho_u^{2\nu}} \right)
-\Phi^{-1}\circ F_{\rho_k}
\left( \frac{l}{x_{sup} \rho_u^{2\nu}} - 1 \right)
\\
& \underset{k \to \infty}{\longrightarrow} \Phi^{-1}\circ F_{\overline\rho}\left( \frac{l}{x_{sup} \rho_u^{2\nu}} \right)
-\Phi^{-1}\circ F_{\overline\rho}
\left( \frac{l}{x_{sup} \rho_u^{2\nu}} -1 \right).
\end{align*}
The convergence comes from the continuity of the function $\rho\mapsto F_{\rho}(t)$ for a fixed $t$ \citep[see the proof of][Lemma A.6]{LLBDR17}. From Lemma \ref{lem:lower:bounded:proba:bounded:matern}, the above limit is finite, which is contradictory with \eqref{hyp:abs}. Hence, \eqref{ass:terme1} is proved. \\

\tbf{4)} Finally, we turn to \eqref{ass:terme2}. We let $m_{n,\rho,y}$ and $\sigma^2 \widetilde k_{n,\rho}$ be the mean and covariance functions of $Y$ given $y$ under covariance function $\widetilde k_{\sigma^2,\rho,\nu}$. Our first aim is to show that, for any $\varepsilon >0$, with probability going to $1$ as $n \to \infty$,
\begin{equation} \label{eq:to:show:limsup:mn}
\underset{ \rho \in [\rho_l , \rho_u] }{\sup}
\underset{ x \in [0 ,1]^d}{\sup}
\left(
m_{n,\rho,y}(x) - Y^* 
\right)
\leqslant \varepsilon,
\end{equation}
and
\begin{equation} \label{eq:to:show:liminf:mn}
\underset{ \rho \in [\rho_l , \rho_u] }{\sup}
\underset{ x \in [0 ,1]^d} {\sup}
\left(
Y_*- m_{n,\rho,y}(x)  
\right)
\leqslant \varepsilon.
\end{equation}

Now we use tools from the theory of reproducing kernel Hilbert spaces (RKHSs) and refer to, e.g., \citep{wendland2004scattered} for the definitions and properties of RKHSs used in the rest of the proof. For $\rho \in [\rho_l , \rho_u] $, the function $m_{n,\rho,y}$ belongs to the RKHS of the covariance function $\widetilde k_{1,\rho,\nu}$. Its RKHS norm $\norme { m_{n,\rho,y} }_{\widetilde k_{1,\rho,\nu}} $ can be simply shown to satisfy
\[
\norme{m_{n,\rho,y}}^2_{\widetilde k_{1,\rho,\nu}}
=
y^{\top} R_{\rho , \nu}^{-1} y.
\]
Hence, from Lemma \ref{lem:mater:change:rho}, observing that 
$\bar{\sigma}_n^2 (\rho) = (1/n) y^{\top} R_{\rho,\nu}^{-1} y$, we have, with probability going to $1$  
as $n \to \infty$,
\begin{equation} \label{eq:controle:rkhs:norm}
\underset{ \rho \in [\rho_l , \rho_u] }{\sup}
\norme{ m_{n,\rho,y} }_{\widetilde k_{1,\rho,\nu}}
\leqslant 
c \sqrt{n}.
\end{equation}

\indent Consider the case a). Since $\nu >1$, the covariance function $\widetilde k_{1,\rho,\nu}$ is twice continuously differentiable on $[0,1]$. Hence, we have from  \citep[Theorem 1]{zhou2008derivative}, 
\[
\underset{ \rho \in [\rho_l , \rho_u] }{\sup}
\underset{x \in [0,1]}{\sup}
| m_{n,\rho,y}'(x) |
\leqslant 
c
\underset{ \rho \in [\rho_l , \rho_u] }{\sup}
\norme{ m_{n,\rho,y} }_{\widetilde k_{1,\rho,\nu}}
\leqslant 
c \sqrt{n},
\]
with probability going to $1$ as $n \to \infty$. Hence, since for $i=1,\ldots,n$ $m_{n,\rho,y}(x_i) = y_i \leqslant Y^*$, and from the assumption $\underset{x \in [0,1]}{\max} \underset{i=1,\ldots,n}{\min} |x - x_i| = o( 1/ \sqrt n )$, it follows that \eqref{eq:to:show:limsup:mn} holds. Similarly, one can show that, with probability going to $1$  
as $n \to \infty$, \eqref{eq:to:show:liminf:mn} holds.\\ 

\indent Consider the case b). 
Since $\nu >2$, the covariance function $\widetilde k_{1,\rho,\nu}$ is four times continuously differentiable on $[0,1]^d$. Hence, we have also from  \citep[Theorem 1]{zhou2008derivative}, 
\begin{equation} \label{eq:bound:second:der}
\underset{ \rho \in [\rho_l , \rho_u] }{\sup}
\underset{i,j = 1,\ldots,d}{\max}
\;
\underset{x \in [0,1]^d}{\sup}
\left|
\frac{\partial^2 m_{n,\rho,y}}{ \partial x_i \partial x_j } (x) 
\right|
\leqslant 
c
\underset{ \rho \in [\rho_l , \rho_u] }{\sup}
\norme{ m_{n,\rho,y} }_{\widetilde k_{1,\rho,\nu}}
\leqslant 
c \sqrt{n},
\end{equation}
with probability going to $1$ as $n \to \infty$. For $\varepsilon >0$,  consider the event 
\begin{equation} \label{eq:event:no:jump:mny}
\left\{\underset{ \rho \in [\rho_l , \rho_u] }{\sup}~~
\underset{ x \in [0 ,1]^d} {\sup}
\left(
m_{n,\rho,y}(x) - Y^* 
\right)
\geqslant \varepsilon\right\}.
\end{equation}
Then, there exists $\bar{\rho}\in [\rho_l , \rho_u]$ and  $ \bar{x} \in [0,1]^d$ for which $m_{n,\bar{\rho},y} ( \bar{x} ) \geqslant Y^* + \varepsilon$. Let $H_{x} m_{n,\bar{\rho},y}$ be the Hessian matrix of $m_{n,\bar{\rho},y}$ at $x$ and $\norme{ H_{x} m_{n,\bar{\rho},y} }$ be its largest singular value.\\

\hspace{0.5cm} $\bullet$ For $d=1$, we consider $\{ v_1,v_2 \} \subset \{ x_1,\ldots,x_n\}$ for which $v_1 \leqslant \bar{x} \leqslant v_2$ with $| v_1 - v_2 | \leqslant 2a_n$ (the existence is assumed in the case b)). Then, necessarily $v_1 < \bar{x} < v_2$.
Because, $m_{n,\bar{\rho},y}(v_1) \leqslant {\max}_{i=1,\ldots,n} \ y_i \leqslant Y^*$, it follows that there exists $w_1 \in [ v_1 , \bar{x} ]$ for which $m'_{n,\bar{\rho},y}(w_1) \geqslant \varepsilon / a_n$. Similarly there exists $w_2 \in [   \bar{x} , v_2]$ for which $m'_{n,\bar{\rho},y}(w_2) \leqslant -\varepsilon /a_n$. Hence, there exists $w_3 \in [v_1,v_2]$ for which $ m''_{n,\bar{\rho},y}(w_3)  \leqslant - \varepsilon / a_n^2$ so that ${\sup}_{x \in [0,1]} \norme{ H_{x} m_{n,\bar{\rho},y} } \geqslant c \varepsilon / a_n^2$. \\

\hspace{0.5cm} $\bullet$ For $d=2$, we consider $\{ v_1,v_2,v_3 \} \subset \{ x_1,\ldots,x_n\}$ for which $ \bar{x}$ belongs to the convex hull of $v_1,v_2,v_3$. Then, if $\bar{x}$ belongs to one of the three segments with end points $v_1,v_2$ or $v_2,v_3$ or $v_1,v_3$, from the previous step with $d=1$, it follows that ${\sup}_{x \in [0,1]^2} \norme{ H_{x} m_{n,\bar{\rho},y} } \geqslant c \varepsilon / a_n^2$. Consider now that $\bar{x}$ does not belong to one of these segments and consider the (unique) intersection point $r$ of the line with direction $ v_1 - \bar{x} $ and of the segment with endpoints $v_2$ and $v_3$. If $m_{n,\bar{\rho},y}(r) \leqslant Y^* + \varepsilon / 2$, by considering the triplet $(v_1 , \bar{x} , r)$, from the reasoning of the case $d=1$, it follows that ${\sup}_{x \in [0,1]^2} \norme{ H_{x} m_{n,\bar{\rho},y} } \geqslant c \varepsilon / a_n^2$. If $m_{n,\bar{\rho},y}(r) \geqslant Y^* + \varepsilon / 2$, by considering the triplet $(v_2 , r , v_3)$, it also follows that ${\sup}_{x \in [0,1]^2} \norme{ H_{x} m_{n,\bar{\rho},y} } \geqslant c \varepsilon / a_n^2$. \\

\hspace{0.5cm} $\bullet$ For $d=3$, we consider $\{ v_1,v_2,v_3,v_4 \} \subset \{ x_1,\ldots,x_n\}$ for which $ \bar{x}$ belongs to the convex hull of $v_1,v_2,v_3,v_4$.
Let $\mathrm{ch}( z_1,z_2,z_3 )$ be the convex hull of $z_1,z_2,z_3 \in [0,1]^d$ (a two-dimensional triangle). 
If $\bar{x}$ belongs to one of the four triangles $\mathrm{ch}( v_1,v_2,v_3 )$, $\mathrm{ch}( v_1,v_2,v_4 )$, $\mathrm{ch}( v_1,v_3,v_4 )$, $\mathrm{ch}( v_2,v_3,v_4 )$, then from the previous step with $d=2$, it follows that ${\sup}_{x \in [0,1]^3} \norme{ H_{x} m_{n,\bar{\rho},y} } \geqslant c \varepsilon / a_n^2$. Now if $\bar{x}$ does not belongs to one of these triangles, then there exists a plane $P_l$ containing $\bar{x}$, intersecting $\mathrm{ch}( v_1,v_2,v_3 )$, $\mathrm{ch}( v_1,v_2,v_4 )$, $\mathrm{ch}( v_1,v_3,v_4 )$, and being parallel to $\mathrm{ch}( v_2,v_3,v_4 )$. Let $E$ be the intersection of this plane $P_l$ and of $\mathrm{ch}( v_1,v_2,v_3 ) \cup \mathrm{ch}( v_1,v_2,v_4 ) \cup \mathrm{ch}( v_1,v_3,v_4 )$. If there exists $x \in E$ so that $m_{n,y,\bar{\rho}}(x) \geqslant Y^* + \varepsilon /2$, then from the previous step with $d=2$, it follows that ${\sup}_{x \in [0,1]^3} \norme{ H_{x} m_{n,\bar{\rho},y} } \geqslant c \varepsilon / a_n^2$. If for all $x \in E$, $m_{n,y,\bar{\rho}}(x) \leqslant Y^* + \varepsilon /2$, then there exists $z_1,z_2,z_3 \in E$ so that $\bar{x} \in \mathrm{ch}( z_1,z_2,z_3 )$ and hence we obtain ${\sup}_{x \in [0,1]^3} \norme{ H_{x} m_{n,\bar{\rho},y} } \geqslant c \varepsilon / a_n^2$.\\

Hence eventually, in all the configurations of the case b), we have ${\sup}_{x \in [0,1]^d} \norme{ H_{x} m_{n,\bar{\rho},y} }\geqslant c \varepsilon / a_n^2$, under the event \eqref{eq:event:no:jump:mny}. Hence, from \eqref{eq:bound:second:der}, \eqref{eq:to:show:limsup:mn} follows in the case b). Analogously, \eqref{eq:to:show:liminf:mn} holds in that case.\\

Similarly to the proof of \eqref{ass:terme2} in the proof of Theorem \ref{th:cond_before_var}  when $\kappa=0$, one can show that ${\sup}_{(\sigma^2 , \rho)  \in \bar{\Theta} } \ \E[Z_{n,\sigma^2 , \rho}^{**}]\to 0$ as $n\to+\infty$ and that ${\sup}_{x\in [0,1]^d} \ \E[Z_{n,\sigma^2,\rho}(x)^2] = \sigma^2
{\sup}_{x\in [0,1]^d} \ \widetilde k_{n , \rho}(x,x)$ goes to zero uniformly in $(\sigma^2 , \rho)  \in \bar{\Theta}$ as $n \to \infty$ for any compact $\bar{\Theta} \subset (0,\infty)^2$. Here $Z_{n,\sigma^2,\rho}$ is defined as in Lemma \ref{lem:mean:value:sup:GP}. We also have as in the proof of Theorem \ref{th:cond_before_var}  when $\kappa=0$ that
\begin{align*}
\P_{\sigma^2 , \rho}(Y^*> u| y)
&\leqslant 
\exp\left\{-\frac{\left(\left(u-\underset{x\in [0,1]^d}{\sup} m_{n,\rho,y}(x)-\E[Z_{n,\sigma^2, \rho}^{**}]\right)_+\right)^2}{2\underset{x\in [0,1]^d}{\sup} \E[Z_{n,\sigma^2, \rho}(x)^2]}\right\}.
\end{align*}
Similarly, we bound the probability of the event  $\{Y_*< \ell\}$ conditionally to $y$.
Hence, we conclude (as in the proof of Theorem \ref{th:cond_before_var}  when $\kappa=0$), also from \eqref{eq:Ystar:no:concentration:u} and \eqref{eq:Ystar:no:concentration:l}, that 
\[
\underset{ (\sigma^2 , \rho) \in \bar{\Theta} }{\sup}
| B_n( \sigma^2 , \rho ) |
=o_{ \P | Y\in \mathcal E_0 }(1).
\]
Consequently, \eqref{ass:terme2} follows and the proof is concluded.
\hfill $\square$\\

\section{Proofs for Sections \ref{sec:var} and \ref{sec:matern} - Monotonicity}\label{sec:proof:mono}

We let $\kappa = 1$ throughout Appendix \ref{sec:proof:mono}.

\subsection{Estimation of the variance parameter}\label{ssec:proof:bounded_var:mono}


\tbf{Proof of Theorem \ref{th:cond_after_var}  under monotonicity constraints}. 
The proof is similar to that of Theorem \ref{th:cond_after_var}   when $\kappa=0$ and is also divided into the three steps \tbf{1)}, \tbf{2)} and \tbf{3)}. \\

\tbf{1)} 
For $n\in \N$, let $N_n$ be the greatest integer such that Condition-Grid holds. Now we define 
\begin{align*}
m_{i,n}=\min \Bigg\{ &\frac{1}{v_{j_i}^{(i)}-v_{j_{i-1}}^{(i)}} \Bigg( 
Y\left(v_{j_1}^{(1)},\ldots,v_{j_{i-1}}^{(i-1)},v_{j_i}^{(i)},v_{j_{i+1}}^{(i+1)},\ldots, v_{j_d}^{(d)}\right)
\\
&-
Y\left(v_{j_1 }^{(1)},\ldots,v_{j_{i-1}}^{(i-1)},v_{j_i -1}^{(i)},v_{j_{i+1}}^{(i+1)},\ldots, v_{j_d}^{(d)}\right)\Bigg),
\\
&
j_k \in \{1,\ldots,N_n\}, \; \forall k\in \{1,\ldots,d\}\diagdown \{i\},\; 
j_i\in \{2,\ldots,N_n\}
\Bigg\},
\end{align*}
for all $i=1,\ldots,d$. Since $N_n \to \infty$ as $n\to \infty$, $m_{i,n} \xrightarrow[n\to +\infty]{\textrm{a.s.}} {\inf}_{x\in [0,1]^d} \ \partial Y(x)/\partial x_i$ since $Y$ is $\mathcal C^1$ a.s. by Condition-Var. Now we notice that 
$m_{i,n}=g_{i,n}(y_1,\ldots,y_n)$ and we define 
$m_{k,i,n}= g_{i,k}(y_1,\ldots,y_k)$. One can see that a slightly different version of Lemma \ref{lem:discret} can be shown (up to re-indexing $x_1,\ldots,x_n$) with $m_{k,n}=\min\{m_{k,1,n},\ldots,m_{k,d,n}\}$, $m=\inf_{x\in [0,1]^d} \min_{i=1,\ldots,d} \partial Y(x)\partial x_i$, and $M_{k,n}=M=u=+\infty$. After applying this different version, points \tbf{2)} and  \tbf{3)} in the proof of Theorem \ref{th:cond_after_var}   when $\kappa=0$ remain unchanged. This concludes the proof.
\hfill $\square$\\


\tbf{Proof of Theorem \ref{th:cond_before_var}  under monotonicity constraints}. 
The proof is similar to that of Theorem \ref{th:cond_before_var}  when $\kappa=0$ and is also divided into the five steps \tbf{1)} to \tbf{5)}. 
We apply Lemma \ref{lem:c_to_hat} to the sequences of functions $f_n$, $g_n$ and $h_n$ defined by $f_n( \sigma^2 ) = \mathcal{L}_n( \sigma^2  )$, $g_n(x) = A_n( \sigma^2)$ and $h_n( \sigma^2 ) = B_n( \sigma^2 )$. Here we recall that for $\sigma^2\in \Theta$,
\[
A_n( \sigma^2  )
= -\ln \P_{\sigma^2}
\left( Y\in\mathcal E_1 \right) 
\;
\mbox{and}
\;
B_n( \sigma^2 )
= \ln \P_{\sigma^2 }
\left( \left. Y\in\mathcal E_1 \right| y \right) .
\]

In order to apply Lemma \ref{lem:c_to_hat}, we need to check that the conditions \eqref{ass:vois} to \eqref{ass:terme2} hold.\\

\tbf{1)} and \tbf{2)} The proof that \eqref{ass:vois} and \eqref{ass:loin} are satisfied is identical to the proof for Theorem \ref{th:cond_before_var}  when $\kappa=0$, as \eqref{ass:vois} and \eqref{ass:loin} do not involve the event $\{ Y\in\mathcal E_1 \}$. \\

\tbf{3)} Let us introduce the Gaussian process $Y_r$ with mean function zero and covariance function $\widetilde{k}_1$. Then we have 
\[
A_n(\sigma^2) = -\ln  \mathbb{P} \left(
\forall x \in [0,1]^d , \forall i=1,\ldots,d, \sigma \frac{\partial}{\partial x_i} Y_r \geqslant 0  \right).
\]
Hence $A_n(\sigma^2)$ does not depend on $\sigma^2$ so that \eqref{ass:terme1} holds. \\

\tbf{4)} We turn to 
\[
B_n(\sigma^2) = \ln \P_{\sigma^2}
\left(
\left.
\forall x \in [0,1]^d , \forall i=1,\ldots,d,  \frac{\partial}{\partial x_i} Y \geqslant 0
\right|
y
\right).
\] 
For $i=1, \ldots ,d$,
let $m^{(1,i)}_{n,y}$ and $\sigma^2\widetilde k^{(1,i)}_{n}$ be the conditional mean and covariance function of $ \partial Y/ \partial x_i$ given $y$, under the probability measure $\P_{\sigma^2}$. We obtain
using Borell-TIS inequality \citep{Adler07} and a union bound, with $Z^{(1,i)}_{n , \sigma^2}$ a Gaussian process with mean function zero and covariance function $\sigma^2 \widetilde k^{(1,i)}_{n}$, 
\begin{align}
\P_{\sigma^2} \nonumber
&
\left(
\left.
\exists x \in [0,1]^d , \exists i=1,\ldots,d,  \frac{\partial}{\partial x_i} Y(x) \leqslant 0
\right|
y
\right)
\\
&\leqslant 
\sum_{i=1}^d
\P_{\sigma^2}
\left( \left.
\sup_{x \in [0,1]^d} \left(- Z^{(1,i)}_{n,\sigma^2}(x) \right)  \geqslant \inf_{x \in [0,1]^d} m_{n,y}^{(1,i)}(x) 
\right| y
\right) \nonumber\\
&\leqslant 
\sum_{i=1}^d
\P_{\sigma^2}
\left( \left.
\sup_{x \in [0,1]^d} \left| Z^{(1,i)}_{n,\sigma^2}(x) \right|  \geqslant \inf_{x \in [0,1]^d} m_{n,y}^{(1,i)}(x) 
\right| y
\right) \nonumber\\
&\leqslant 
\sum_{i=1}^d
\exp\left\{
-\frac{
	\left(\left(
	\underset{x\in [0,1]^d}{\inf} m_{n,y}^{(1,i)}(x)-
	\E \left[
	\sup_{x \in [0,1]^d} \left| Z^{(1,i)}_{n,\sigma^2}(x) \right|
	\right]
	\right)_+\right)^2
}{
2\underset{x\in [0,1]^d}{\sup} \E[Z^{(1,i)}_{n,\sigma^2}(x)^2]}
\right\}.\label{ineq:BTIS:mono}
\end{align}

One can see that Lemma \ref{lem:mean:value:sup:GP} can also be shown when $Z_{n,\theta}$ is replaced by $Z^{(1,i)}_{n,\theta}$ (here $\theta = \sigma^2$). Hence 
$
\sup_{\sigma^2 \in [\sigma_l^2 , \sigma_u^2 ]} \
\E \big[
\sup_{x \in [0,1]^d} \big| Z^{(1,i)}_{n,\sigma^2}(x) \big|
\big]$ goes to $0$ as $n \to \infty$. 
Additionally, one can simply show that ${\sup}_{x\in [0,1]^d} \ \E[Z^{(1,i)}_{n,\sigma^2}(x)^2] =
{\sup}_{x\in [0,1]^d} \ \sigma^2 \widetilde k^{(1,i)}_{n }(x,x)$ goes to zero uniformly in $\sigma^2 \in [ \sigma_l^2 , \sigma_u^2 ]$ as $n \to \infty$. 

One can see that the proof of \citep[Proposition 2.8]{BBG16} can be adapted to establish that, for $i = 1, \ldots , d$,
\[
\displaystyle \underset{x\in [0,1]^d}{\sup} 
\left|
m^{(1,i)}_{n,y}(x)- \frac{\partial }{\partial x_i} Y(x)
\right|
\xrightarrow[n\to +\infty]{\textrm{a.s.}} 0,
\]
from which we deduce that on the set $\{Y\in \mathcal E_{1, \delta}\}$, where
\[
\mathcal E_{1, \delta}
\defeq
\left\{ f\in \mathcal C^1([0,1]^d,\R) \quad \textrm{s.t.} \;  \partial f(x) / \partial x_i \geqslant \delta, \; \forall x\in [0,1]^d, i \in \{ 1,\ldots,d\}\right\},
\]
we have a.s., for $i=1,\ldots,d$,
\[
\displaystyle\liminf_{n\to+\infty } 
\left(
\underset{x\in [0,1]^d}{\inf} m_{n,y}^{(1,i)}(x)
\right)
\geqslant \delta.
\]
Consequently, from \eqref{ineq:BTIS:mono}, on $\{Y\in \mathcal E_{1, \delta}\}$, we have:
\begin{align}\label{limit:BTIS_sup:mono}
\underset{\sigma^2\in [\sigma_l^2, \sigma_u^2]}{\sup}\P_{\sigma^2} 
\left(
\left.
\exists x \in [0,1]^d , \exists i=1,\ldots,d,  \frac{\partial}{\partial x_i} Y(x) \leqslant 0
\right|
y
\right)
\xrightarrow[n\to +\infty]{\textrm{a.s.}} 0.
\end{align}

Similarly as in the proof of Theorem \ref{th:cond_before_var}  when $\kappa=0$, we can show, by applying Tsirelson theorem in \citep{AW09} to the processes $\partial Y / \partial x_i$, that 
\[
\mathbb{P}_{\sigma_0^2} \left(
Y \in \mathcal{E}_1 \cap \mathcal E_{1,\delta}^c
\right)
\underset{\delta \to 0}{\longrightarrow}  0.
\]
Hence we conclude the proof of \eqref{ass:terme2} as for Theorem \ref{th:cond_before_var}  when $\kappa=0$. \\

\tbf{5)} 
We conclude the proof as in \tbf{5)} for Theorem \ref{th:cond_before_var}  when $\kappa=0$.
\hfill $\square$

\subsection{Isotropic Mat\'ern process}\label{ssec:proof:bounded_matern:mono}

\tbf{Proof of Theorem \ref{th:cond_after_matern}  under monotonicity constraints}. The proof is the same as for Theorem \ref{th:cond_after_matern} when $\kappa=0$ and is concluded by applying Theorem \ref{th:cond_after_var} when $\kappa=1$.
\hfill $\square$\\

\tbf{Proof of Theorem \ref{th:cond_before_matern_i} under monotonicity constraints}. The proof is similar to that of Theorem \ref{th:cond_before_matern_i} when $\kappa=0$ and is also divided into the four steps \tbf{1)} to \tbf{4)}.  We apply Lemma \ref{lem:c_to_hat} to the sequences of functions $f_n$, $g_n$ and $h_n$ defined by $f_n( \sigma^2 ) = \mathcal{L}_n( \sigma^2 , \rho_1 )$, $g_n(x) = A_n( \sigma^2 , \rho_1 )$, and $h_n( \sigma^2 ) = B_n( \sigma^2 , \rho_1)$.\\

\tbf{1)} The proof that \eqref{ass:vois} and \eqref{ass:loin} are satisfied is identical to the proof of Theorem \ref{th:cond_before_matern_i} when $\kappa=0$, as 
\eqref{ass:vois} and \eqref{ass:loin} do not involve the event $\{Y\in \mathcal E_1\}$.\\

\tbf{2)} Let us introduce the Gaussian process $Y_r$ with mean function zero and covariance function $\widetilde k_{1,\rho_1,\nu}$. Then we have 
\[
A_n(\sigma^2, \rho_1) = - \ln  \mathbb{P} \left(
\forall x \in [0,1]^d , \forall i=1,\ldots,d, \sigma \frac{\partial}{\partial x_i} Y_r(x) \geqslant 0  \right).
\]
Hence $A_n(\sigma^2, \rho_1)$ does not depend on $\sigma^2$ so that \eqref{ass:terme1} holds. \\

\tbf{3)} We turn to $B_n(\sigma^2 , \rho_1)$. We conclude to \eqref{ass:terme2} following the same lines as in the proof of Theorem \ref{th:cond_before_matern_i}  when $\kappa=0$ and using the equivalence of measures.\\

\tbf{4)} We conclude the proof of Theorem \ref{th:cond_before_matern_i}  when $\kappa=1$ similarly as in the proof of Theorem \ref{th:cond_before_matern_i}  when $\kappa=0$ using Theorem \ref{th:cond_after_matern}  when $\kappa=1$. \hfill $\square$\\

\tbf{Proof of Theorem \ref{th:cond_before_matern_ii} under monotonicity constraints}. 
The proof follows the similar four steps of the proof Theorem \ref{th:cond_before_matern_ii}  when $\kappa=0$. We apply Lemma \ref{lem:c_to_hat} to the sequences of functions $f_n$, $g_n$ and $h_n$ defined by $f_n( x ) = \mathcal{L}_n( x \widehat{\rho}_{n,c}^{2 \nu}  , \widehat{\rho}_{n,c} )$, $g_n(x) = A_n( x \widehat{\rho}_{n,c}^{2 \nu} , \widehat{\rho}_{n,c} )$, and $h_n( x ) = B_n( x \widehat{\rho}_{n,c}^{2 \nu} , \widehat{\rho}_{n,c} )$.  \\

\tbf{1)} The proof that \eqref{ass:vois} and \eqref{ass:loin} are satisfied is identical to the proof of Theorem \ref{th:cond_before_matern_ii}  when $\kappa=0$, as 
\eqref{ass:vois} and \eqref{ass:loin} do not involve the event $\{Y\in \mathcal E_1\}$.\\

\tbf{2)} Similarly as in the proof of Theorem \ref{th:cond_before_matern_ii} when $\kappa=0$, we show that Theorem \ref{th:cond_before_matern_ii} when $\kappa=1$ holds if \eqref{ass:terme1} and \eqref{ass:terme2} are satisfied.\\

\tbf{3)} Similarly as in the proof of Theorem \ref{th:cond_before_matern_i} when $\kappa=1$, we show that \eqref{ass:terme1} holds. \\

\tbf{4)} Finally, we turn to \eqref{ass:terme2}. First, consider the case a). Recall the notation $Y'_*= {\inf}_{ x \in [0 ,1]} \ Y'(x)$ from Appendix \ref{sec:intermediate}. We proceed similarly as in the proof of Theorem \ref{th:cond_before_matern_ii} when $\kappa=0$. Since $\E_{\rho}[Y'(x)|y]=m_{n,\rho,y}'(x)$, it is then sufficient to show that, for all $i=1,\ldots, d$, for any $\varepsilon >0$, with probability going to $1$ as $n \to \infty$,
\begin{equation} \label{eq:to:show:liminf:mn:mono:a}
\underset{ \rho \in [\rho_l , \rho_u] }{\sup}
\;
\underset{ x \in [0 ,1]} {\sup}
\left(
Y'_*-  m'_{n,\rho,y}(x)  
\right)
\leqslant \varepsilon,
\end{equation}
in order to prove \eqref{ass:terme2} as in the proof of Theorem \ref{th:cond_before_matern_ii}  when $\kappa=0$. Analogously, since $\nu>2$, \eqref{eq:bound:second:der} holds. Consider $\bar x$, $\bar\rho$ so that  
$m'_{n,\bar \rho,y}(\bar x)\leqslant  Y'_* - \varepsilon$. There exists $i, j \in \{1,\ldots,n\}$, $x_i<x_j$ such that $\abs{x_i-\bar x}\leqslant a_n$ and   $\abs{x_j-\bar x}\leqslant a_n$. We have $m_{n,\bar \rho,y}(x_j)-m_{n,\bar \rho,y}(x_i)=Y(x_j)-Y(x_i)\geqslant (x_j-x_i)Y'_*$. Thus from the mean value theorem, there exists $w\in [0,1]$ so that $\abs{\bar x-w}\leqslant 2a_n$ and $\abs{m'_{n,\bar \rho,y}(w)-m'_{n,\bar \rho,y}(\bar x)}\geqslant \varepsilon$. Hence there exists $z\in [0,1]$ so that $m''_{n,\bar \rho,y}(z)\geqslant \varepsilon /(2a_n)$. Hence \eqref{eq:to:show:liminf:mn:mono:a} holds from \eqref{eq:bound:second:der}.\\

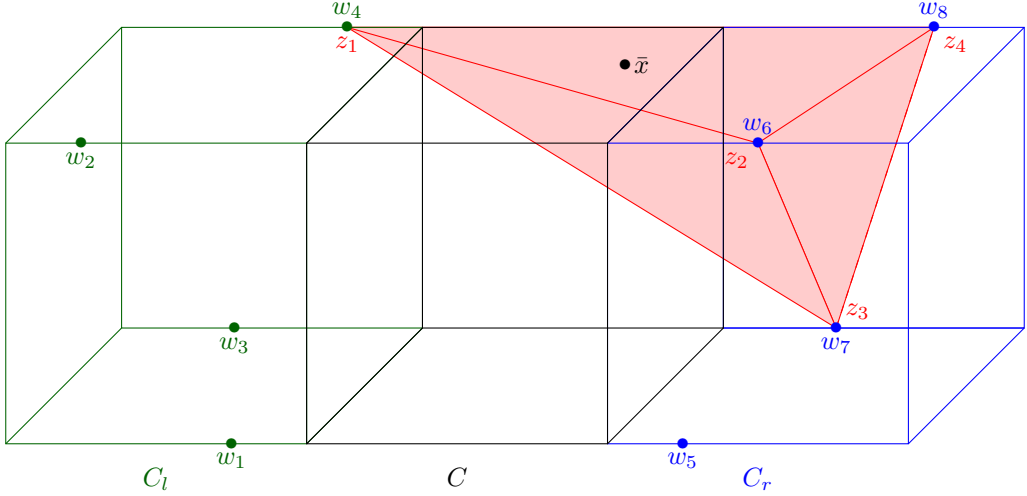
\begin{figure}
\centering
\begin{tikzpicture}

\draw[red,fill=red!20] (-5,4,0)--(2,4,4)--(1.5,0,0)--(2.8,4,0)--cycle;
\draw[red,fill=red!20] (-5,4,0)--(1.5,0,0)--(2,4,4)--(2.8,4,0)--cycle;

\draw (-1,3.8,0.8) node[right]{$\bar x$};
\draw (-1,3.8,0.8) node {$\bullet$};

\draw[red,dashed] (1.5,0,0)--(2.8,4,0);
\draw[red] (-5,4,0)--(2,4,4);

\draw[greenmoi] (-8,0,0)--(-4,0,0)--(-4,4,0)--(-8,4,0)--cycle; 
\draw[greenmoi] (-8,0,4)--(-4,0,4)--(-4,4,4)--(-8,4,4)--cycle; 
\draw[greenmoi] (-8,0,0) -- (-8,0,4); 
\draw[greenmoi] (-4,0,0) -- (-4,0,4); 
\draw[greenmoi] (-4,4,0) -- (-4,4,4); 
\draw[greenmoi] (-8,4,0) -- (-8,4,4); 

\draw[greenmoi] (-5,0,4) node[below]{$w_1$};
\draw[greenmoi] (-5,0,4) node {$\bullet$} ;

\draw[greenmoi] (-7,4,4) node[below]{$w_2$};
\draw[greenmoi] (-7,4,4) node {$\bullet$} ;

\draw[greenmoi] (-6.5,0,0) node[below]{$w_3$};
\draw[greenmoi] (-6.5,0,0) node {$\bullet$} ;

\draw[greenmoi] (-5,4,0) node[above]{$w_4$};
\draw[red] (-5,4,0) node[below]{$z_1$};
\draw[greenmoi] (-5,4,0) node {$\bullet$} ;

\draw[blue] (0,0,0)--(4,0,0)--(4,4,0)--(0,4,0)--cycle; 
\draw[blue] (0,0,4)--(4,0,4)--(4,4,4)--(0,4,4)--cycle; 
\draw[blue] (0,0,0) -- (4,0,0); 
\draw[blue] (4,0,0) -- (4,0,4); 
\draw[blue] (4,4,0) -- (4,4,4); 
\draw[blue] (0,4,0) -- (0,4,4); 

\draw[blue] (1,0,4) node[below]{$w_5$};
\draw[blue] (1,0,4) node {$\bullet$} ;

\draw[blue] (2,4,4) node[above]{$w_6$};
\draw[red] (2,4,4) node[below left]{$z_2$};
\draw[blue] (2,4,4) node {$\bullet$} ;

\draw[blue] (1.5,0,0) node[below]{$w_7$};
\draw[red] (1.5,0,0) node[above right]{$z_3$};
\draw[blue] (1.5,0,0) node {$\bullet$} ;

\draw[blue] (2.8,4,0) node[above]{$w_8$};
\draw[red] (2.8,4,0) node[below right]{$z_4$};
\draw[blue] (2.8,4,0) node {$\bullet$} ;

\draw (-4,0,0)--(0,0,0)--(0,4,0)--(-4,4,0)--cycle; 
\draw (-4,0,4)--(0,0,4)--(0,4,4)--(-4,4,4)--cycle; 
\draw (-4,0,0) -- (-4,0,4); 
\draw (0,0,0) -- (0,0,4); 
\draw (0,4,0) -- (0,4,4); 
\draw (-4,4,0) -- (-4,4,4); 

\draw (-2,-0.2,4) node[below]{$C$};
\draw[greenmoi] (-6,-0.2,4) node[below]{$C_l$};
\draw[blue] (2,-0.2,4) node[below]{$C_r$};

\end{tikzpicture}
\caption{The point $\bar x$ belongs to the  hypercube $C$ (in black, central) that lies between the hypercubes $C_l$ (in green, left) and $C_r$ (in blue, right). The mean value theorem ensures the existence of the 8 points $w_j$, $j=1,\ldots,8$, and $\bar{x}$ belongs to the convex hull $\mathrm{ch}( w_1,\ldots,w_8)$. Furthermore, $\bar{x}$ also belongs to $\mathrm{ch}( z_1,\ldots,z_4)$, with $\{z_1,\ldots,z_4\} \subset \{w_1,\ldots,w_8\}$ (in red).}\label{fig:mono}
\end{figure}

\indent Second, consider the case b). We shall only address the case $d=3$, the cases $d=1,2$ being treated similarly. For all $i=1,2,3$, let us define $Y^{(1,i)}=\partial Y/\partial x_i$ and $Y^{(1,i)}_*= {\inf}_{ x \in [0 ,1]^3} \ Y^{(1,i)}(x)$ following the notation of Appendix \ref{sec:intermediate}. Let also $m^{(1,i)}_{n,\rho,y}(x)=\E_{\rho}[Y^{(1,i)}(x)|y]=\partial m_{n,\rho,y}(x)/\partial x_i$. First, we want to show that, for all $i=1,2, 3$, for any $\varepsilon >0$, with probability going to $1$ as $n \to \infty$,
\begin{equation} \label{eq:to:show:liminf:mn:mono:b}
\underset{ \rho \in [\rho_l , \rho_u] }{\sup}
\;
\underset{ x \in [2/n^{1/3} ,1-2/n^{1/3}]^3} {\sup}
\left(Y^{(1,i)}_*-  m_{n,\rho,y}^{(1,i)}(x)  
\right)
\leqslant \varepsilon.
\end{equation}
Assume that there exists $i\in \{1,2,3\}$ and $ \bar{x} , \bar{\rho}$ for which $m^{(1,i)}_{n,\bar{\rho},y} ( \bar{x} ) \leqslant Y^{(1,i)}_* - \varepsilon$. There exists $\{i_1, \ldots, i_8\}\in \{1,\ldots,n\}$
so that, $\bar x$ belongs to the  hypercube $C$ with vertices $x_j$, $j\in \{i_1,\ldots, i_8\}$. We refer to Figure \ref{fig:mono} for an illustration. This hypercube lies between two adjacent hypercubes $C_l,C_r$ (with vertices in $\{x_1,\ldots, x_n\}$ and edge lengths $1/\left( \lfloor n^{1/3}\rfloor - 1 \right)$) which are obtained 
by translations (to the left and to the right) in the direction $i$.
Note that $C,C_l,C_r$ are disjoint and that the pairs $C,C_l$ and $C,C_r$ each have a common face which is orthogonal to the direction $i$. We now consider the $8$ vertices $v_1 , \ldots , v_8$ of $C_l$ and $C_r$ which are parallel to the direction $i$. For any $j =1,\ldots,8$, the endpoints of $v_j$ can be written as $x_a, x_b$ with $(x_a)_k = (x_b)_k$ for $k \neq i$ and with $(x_a)_i < (x_b)_i$. Then we have $m_{n,\bar \rho ,y}(x_b) - m_{n,\bar \rho ,y}(x_a) = Y(x_b) - Y(x_a) \geqslant Y^{(1,i)}_* ( (x_b)_i - (x_a)_i )$. Hence, from the mean value theorem, there exists $w_j \in v_j$ for which $m_{n,\bar \rho ,y}^{(1,i)}(w_j) \geqslant Y^{(1,i)}_* $. Also, it can be shown that $\bar{x}$ belongs to $\mathrm{ch}( w_1,\ldots,w_8)$. In addition, it can also be shown that $\bar{x}$ belongs to $\mathrm{ch}( z_1,\ldots,z_4)$, with $\{z_1,\ldots,z_4\} \subset \{w_1,\ldots,w_8\}$. Since $m_{n,\bar \rho ,y}^{(1,i)}(z_j) \geqslant Y^{(1,i)}_* $ for $j=1,\ldots,4$ and since $m^{(1,i)}_{n,\bar{\rho},y} ( \bar{x} ) \leqslant Y^{(1,i)}_* - \varepsilon$, we show \eqref{eq:to:show:liminf:mn:mono:b} as in the proof of  Theorem \ref{th:cond_before_matern_ii} b)  when $\kappa=0$.

We have, for $i =1,2,3$,
\begin{align*}
\P_{(\sigma^2,\rho)}
\left(
\left.
\underset{x \in [0,1]^3}{\inf} 
Y^{(1,i)}(x)
\leqslant 0
\right|
y
\right)
\leqslant & 
\P_{(\sigma^2,\rho)}
\left(
\left.
\underset{ x \in [2/n^{1/3} ,1-2/n^{1/3}]^3} {\inf}
Y^{(1,i)}(x)
\leqslant \delta/2
\right|
y
\right)
\\
& +
\P_{(\sigma^2,\rho)}
\left(
\left.
\underset{x \in [0,1]^3}{\inf} 
Y^{(1,i)}(x)
\leqslant 0,
\underset{ x \in [2/n^{1/3} ,1-2/n^{1/3}]^3} {\inf}
Y^{(1,i)}(x)
> \delta/2
\right|
y
\right)
\\
\eqdef & P_1 + P_2,
\end{align*}
say. As in the proof of Theorem \ref{th:cond_before_matern_ii} b)  when $\kappa=0$, we can show from \eqref{eq:to:show:liminf:mn:mono:b} that $ \mathds{1}_{ Y^{(1,i)}_* \geqslant \delta } {\sup}_{(\sigma^2,\rho) \in \bar{\Theta}} \ P_1 $ goes to zero in probability so that it is sufficient to show ${\sup}_{(\sigma^2,\rho) \in \bar{\Theta}} \ P_2 $ goes to zero in probability where $\bar{\Theta}$ is any compact set of $(0 , \infty)^2$. We have, with $Y^{(2,j,i)}(x) = \partial^2 Y(x)/ (\partial x_j \partial x_i) $, and with $m_{n,\rho,y}^{(2,j,i)}(x) = \partial^2 m_{n,\rho,y}(x)/ (\partial x_j \partial x_i)  = \mathbb{E} [ Y^{(2,j,i)}(x)  | y ]$,
\begin{align*}
P_2
\leqslant &
\P_{(\sigma^2,\rho)}
\left(
\left.
\max_{j=1,2,3}
\underset{x \in [0,1]^3}{\sup} 
\left| Y^{(2,j,i)}(x) \right|
\geqslant c n^{1/3}
\right|
y
\right)
\\
\leqslant &
\P_{(\sigma^2,\rho)}
\left(
\left.
\max_{j=1,2,3}
\underset{x \in [0,1]^3}{\sup} 
\left| m_{n,\rho,y}^{(2,j,i)}(x) \right|
\geqslant (c/2) n^{1/3}
\right|
y
\right)
\\
& + 
\P_{(\sigma^2,\rho)}
\left(
\left.
\max_{j=1,2,3}
\underset{x \in [0,1]^3}{\sup} 
\left| 
m_{n,\rho,y}^{(2,j,i)}(x) - Y^{(2,j,i)}(x)
\right|
\geqslant (c/2) n^{1/3}
\right|
y
\right).
\end{align*}
Using the Borell-TIS inequality as in the proof of \ref{th:cond_before_var}  when $\kappa=1$, we can show that 
\[
\underset{\rho,\sigma^2 \in \bar{\Theta}}{\sup}
~ ~
\P_{(\sigma^2,\rho)}
\left(
\left.
\max_{j=1,2.3}
\underset{x \in [0,1]^3}{\sup} 
\left| 
m_{n,\rho,y}^{(2,j,i)}(x) - Y^{(2,j,i)}(x)
\right|
\geqslant (c/2) n^{1/3}
\right|
y
\right)
=o_{\P}(1).
\]
Now consider that there exist $\bar{\rho} \in [ \rho_l , \rho_u ]$, $\bar{j}  \in \{1,2,3 \}$ and $\bar{x} \in [0,1]^3$ so that $| m_{n,\bar{\rho},y}^{(2,\bar j,i)}(\bar{x }) | \geqslant (c/2) n^{1/3}$. 
If $\bar{j} = i$, then by applications of the mean value theorem, by using that $m_{n,\rho,y}(x) = Y(x)$ for $x \in \{ x_1, \ldots , x_n\}$, we can show that there exists $\bar{w} \in [0,1]^3$ so that $\norme{ \bar{w} - \bar{x} } \leqslant c n^{-1/3}$ and $\abs{ m_{n,\bar{\rho},y}^{(2,\bar j,i)}(\bar{w }) } \leqslant  2 \max_{i,j \in \{ 1,2,3\}} Y^{(2,j,i)**}$ where with the notation of Appendix \ref{sec:intermediate}, $Y^{(2,j,i)**} \defeq {\sup}_{x \in [0,1]^3} |  Y^{(2,j,i)}(x)|$. Hence there exists $\bar{r} \in [0,1]^3$ so that ${\max}_{j,k \in \{ 1,2,3\}}  | m_{n,\rho,y}^{(3,k,j,i)}(\bar{r}) | \geqslant  c n^{2/3}$ where
$m_{n,\rho,y}^{(3,k,j,i)}(x) = \partial^3 / (\partial x_k \partial x_j \partial x_i) \ m_{n,\rho,y}(x)$.

If $\bar{j} \neq i$, we can consider  $\bar{z} \in \{ x_1, \ldots , x_n\}$ so that $ \norme{\bar{z} - \bar{x} } \leqslant c n^{-1/3}$. We also consider $9$ additional points $(\bar{z}_{k,l})_{k,l \in \{ 1,2,3\}} \subset \{ x_1, \ldots , x_n\}$ so that $\bar{z}_{k,l} = \bar{x} + l / ( \lfloor n^{1/3} \rfloor - 1 ) v_k$ where $v_1 = e_i$ with $e_i$ the $i$-th base column vector, $v_2 = e_{\bar j}$ and $v_3 = e_i + e_{\bar j}$. By applications of the mean value theorem, we can show that there exist $\bar{w}_1$ for which $|m_{n,\bar{\rho},y}^{(2, i ,i)}(\bar{w }_1) | \leqslant 3  Y^{(2,i,i)**}$, $\bar{w}_2$ for which $| m_{n,\bar{\rho},y}^{(2, \bar j , \bar j)}(\bar{w }_2) | \leqslant 3  Y^{(2, \bar j, \bar j)**}$ and $\bar{w}_3$ for which $| m_{n,\bar{\rho},y}^{(2, i , i)}(\bar{w }_3) + m_{n,\bar{\rho},y}^{(2, \bar j , \bar j)}(\bar{w }_3) + 2 m_{n,\bar{\rho},y}^{(2, \bar j , i)}(\bar{w }_3) | \leqslant 12  \max_{i,j \in \{ 1,2\}} Y^{(2,j,i)**}$. If $| m_{n,\bar{\rho},y}^{(2, \bar j , i)}(\bar{w }_3) | \leqslant (c/4) n^{1/3}$, then there exists $\bar{r} \in [0,1]^3$ so that ${\max}_{j,k \in \{ 1,2,3\}}  | m_{n,\rho,y}^{(3,k,j,i)}( \bar{r} ) | \geqslant  c n^{2/3}$. If $ |m_{n,\bar{\rho},y}^{(2, \bar j , i)}(\bar{w }_3) | \geqslant (c/4) n^{1/3}$, then $| m_{n,\bar{\rho},y}^{(2, i , i)}(\bar{w }_3) | \geqslant (c/4) n^{1/3} - 6  \max_{i,j \in \{ 1,2\}} Y^{(2,j,i)**}$ or 
$| m_{n,\bar{\rho},y}^{(2, \bar j , \bar j)}(\bar{w }_3) | \geqslant (c/4) n^{1/3} - 6  \max_{i,j \in \{ 1,2\}} Y^{(2,j,i)**}$. In all the cases, there exists $\bar{r} \in [0,1]^3$ so that ${\max}_{j,k \in \{ 1,2,3\}}   m_{n,\rho,y}^{(3,k,j,i)}( \bar{r} ) \geqslant  c n^{2/3}$.

We have also from \citep[Theorem 1]{zhou2008derivative}, and since $\nu >3$,
\begin{equation} \label{eq:bound:third:der}
\underset{ \rho \in [\rho_l , \rho_u] }{\sup}
\underset{i,j,k = 1,2,3}{\max}
\underset{x \in [0,1]^3}{\sup}
\left|
\frac{\partial^3}{\partial x_i  \partial x_j \partial x_k } m_{n,\rho,y}(x) 
\right|
\leqslant 
c
\underset{ \rho \in [\rho_l , \rho_u] }{\sup}
\norme{ m_{n,\rho,y} }_{k_{1,\rho,\nu}}
\leqslant 
c \sqrt{n},
\end{equation}
with probability going to $1$ as $n \to \infty$ from \eqref{eq:controle:rkhs:norm}.

Hence, we have that with probability going to $1$ as $n \to \infty$,
\begin{align*}
\underset{\sigma^2 , \rho}{\sup}
\P_{(\sigma^2,\rho)}
\left(
\left.
\max_{j=1,2,3}
\underset{x \in [0,1]^3}{\sup}  
| m_{n,\rho,y}^{(2,j,i)}(x) |
\geqslant (c/2) n^{1/3}
\right|
y
\right)
=
\mathds{1}_{ \{ \max_{j=1,2,3}
	\underset{x \in [0,1]^3}{\sup} 
	m_{n,\rho,y}^{(2,j,i)}(x)
	\geqslant (c/2) n^{1/3} \}}
= 0.
\end{align*}
This proves that \eqref{ass:terme2} holds so that the proof is complete.
\hfill $\square$\\

\section{Proofs for Sections \ref{sec:var} and \ref{sec:matern} - Convexity}\label{sec:proof:convex}

We let $\kappa = 2$ throughout Appendix \ref{sec:proof:convex}.

\subsection{Estimation of the variance parameter}\label{ssec:proof:bounded_var:convex}

\tbf{Proof of Theorem \ref{th:cond_after_var}  under convexity constraints}. 
The proof is similar to that of Theorem \ref{th:cond_after_var}   when $\kappa=1$, where the finite differences of order one are replaced by finite differences of order two. \\

\tbf{Proof of Theorem \ref{th:cond_before_var}  under convexity constraints}. 
The proof is similar to that of Theorem \ref{th:cond_before_var}  when $\kappa=0$ introducing the Gaussian process $V$ defined on $S_1\times \R^d$ by $V(v,x)=v^{\top} HY(x) v $ where $S_1=\{v\in \R^d,\; \norme{v}=1 \}$ and observing that
\[
\mathcal E_2= \Bigg\{ f\in \mathcal C([0,1]^d,\R),\; \forall x \in [0,1]^d , \underset{\substack{v\in \R^d\\ \norme{v}=1}}{\inf} v^{\top} Hf(x) v\leqslant 0 \Bigg\} ,
\]
where $Hf(x)$ represents the Hessian matrix of $f$ at $x$ which means that $Hf(x)_{i,j}=\partial^2  f(x)/(\partial x_i\partial x_j)$.
\hfill $\square$

\subsection{Isotropic Mat\'ern process}\label{ssec:proof:bounded_matern:convex}

\tbf{Proof of Theorem \ref{th:cond_after_matern}  under convexity constraints}. The proof is the same as for Theorem \ref{th:cond_after_matern} when $\kappa=0$ and is concluded by applying Theorem \ref{th:cond_before_var} when $\kappa=2$.
\hfill $\square$\\

\tbf{Proof of Theorem \ref{th:cond_before_matern_i} under convexity constraints}. The proof is similar to that of Theorem \ref{th:cond_before_matern_i} when $\kappa=0$.
\hfill $\square$\\

\tbf{Proof of Theorem \ref{th:cond_before_matern_ii} under convexity constraints}. 
The proof follows the similar four steps of the proof Theorem \ref{th:cond_before_matern_ii}  when $\kappa=0$. Points 
\tbf{1)} to 
\tbf{3)} are identical. Turning to \eqref{ass:terme2}, point 
\tbf{4)} can be treated similarly as in the proof of Theorem \ref{th:cond_before_matern_ii}  when $\kappa=1$ but with more cumbersome notation and arguments. In order to ease the reading of the paper, we omit this technical proof.
\hfill $\square$

\section{Proofs for Section \ref{sec:extensions}}

\tbf{Proof of Proposition \ref{proposition:prediction}}.
We have 
\begin{align*}
\frac{ \widehat{Y}(x_0) - \widehat{Y}_c(x_0) }{ \widehat{\sigma}(x_0) }
& =
\mathbb{E}
\left[
\left.
\frac{ Y(x_0) - \widehat{Y}(x_0) }{ \widehat{\sigma}(x_0) }
\right|
y
\right]
-
\frac{
	\mathbb{E}
	\left[
	\left.
	\frac{ Y(x_0) - \widehat{Y}(x_0) }{ \widehat{\sigma}(x_0) }
	\mathds{1}_{ \{  Y \in \mathcal{E}_{\kappa} \} }
	\right|
	y
	\right]
}{ \mathbb{P}( Y \in \mathcal{E}_{\kappa} | y ) }.
\end{align*}
Now let $E_n(x_0) =  (Y(x_0) - \widehat{Y}(x_0))/ \widehat{\sigma}(x_0) $. We have
\begin{align*}
\left| 
\frac{ \widehat{Y}(x_0) - \widehat{Y}_c(x_0) }{ \widehat{\sigma}(x_0) }
\right|  
\leqslant
\left| 
\mathbb{E}
\left[ 
\left.
E_n(x_0) ( 1 - \mathds{1}_{ \{  Y \in \mathcal{E}_{\kappa} \} }  )
\right| 
y
\right]
\right|
+
\left|
\mathbb{E}
\left[ 
\left.
E_n(x_0) \mathds{1}_{ \{  Y \in \mathcal{E}_{\kappa} \} } 
\right| 
y
\right]
\right|
\left|
1 - \frac{1}{\mathbb{P}( Y \in \mathcal{E}_{\kappa} | y )}
\right| 
\eqdef 
|A| + |B| |C|,
\end{align*}
say. By Cauchy-Schwarz's inequality, we have
\begin{align*}
|A|
& \leqslant 
\mathbb{E}[E_n(x_0)^2 | y] ^{1/2}
\mathbb{P}( Y \not \in \mathcal{E}_{\kappa} | y )^{1/2}.
\end{align*}
In the above display, the first square root is $1$ by definition and the second one is a $o_{\P | Y \in \mathcal{E}_{\kappa}}(1)$, as shown in the point \tbf{4)} in the proof of Theorem \ref{th:cond_before_var}.

By Jensen's inequality, we have $|B|  \leqslant \mathbb{E}[E_n(x_0)^2 | y] ^{1/2}= 1$. Finally, $|C| = o_{\P | Y \in \mathcal{E}_{\kappa}}(1)$, since $\mathbb{P}( Y \not \in \mathcal{E}_{\kappa} | y )   = o_{\P | Y \in \mathcal{E}_{\kappa}}(1)$, as above. This completes the proof of \eqref{eq:pred:means}.

We also have
\begin{align*}
\frac{
	\widehat\sigma(x_0)^2
	-
	\widehat\sigma_c(x_0)^2
}{
\widehat{\sigma}(x_0)^2
}
& =
\mathbb{E}
\left[ 
\left.
\frac{ (Y(x_0) - \widehat{Y}(x_0))^2 }{ \widehat{\sigma}(x_0)^2 }
-
\frac{ (Y(x_0) - \widehat{Y}_c(x_0))^2 }{ \widehat{\sigma}(x_0)^2 }
\right|
y
\right]
\\
& ~~
+
\mathbb{E}
\left[ 
\left.
\frac{ (Y(x_0) - \widehat{Y}_c(x_0))^2 }{ \widehat{\sigma}(x_0)^2 }
\right|
y
\right]
-
\mathbb{E}
\left[ 
\left.
\frac{ (Y(x_0) - \widehat{Y}_c(x_0))^2 }{ \widehat{\sigma}(x_0)^2 }
\right|
y ,  Y \in \mathcal{E}_\kappa 
\right]
\\
& \eqdef  D+E -F,
\end{align*}
say. We have
\begin{align*}
D = 
\mathbb{E}
\left[
\left.
\frac{\widehat{Y}_c(x_0) - \widehat{Y}(x_0)}{ \widehat{\sigma}(x_0) }
\left(
2 \frac{Y(x_0) - \widehat{Y}(x_0)}{ \widehat{\sigma}(x_0) }
+
\frac{\widehat{Y}(x_0) - \widehat{Y}_c(x_0)}{ \widehat{\sigma}(x_0) }
\right)
\right|
y
\right]
= -
\left(
\frac{\widehat{Y}(x_0) - \widehat{Y}_c(x_0)}{ \widehat{\sigma}(x_0) }
\right)^2
= o_{ \P | Y \in \mathcal{E}_{\kappa} } (1),
\end{align*}
from \eqref{eq:pred:means}. We also have
\begin{align*}
\left| E - F \right|
\leqslant 
\left|
\mathbb{E}
\left[ 
\left.
\frac{ (Y(x_0) - \widehat{Y}_c(x_0))^2 }{ \widehat{\sigma}(x_0)^2 }
\mathds{1}_{ \{Y \not \in \mathcal{E}_{\kappa} \}}
\right| 
y
\right]
\right|
+
\left|
\mathbb{E}
\left[ 
\left.
\frac{ (Y(x_0) - \widehat{Y}_c(x_0))^2 }{ \widehat{\sigma}(x_0)^2 }
\mathds{1}_{ \{ Y \in \mathcal{E}_{\kappa} \} }
\right| 
y
\right]
\right| 
\left|C \right|
\eqdef |G| + |H||C|,
\end{align*}
say. By Cauchy-Schwarz's inequality, we have
\begin{align*}
|G|  \leqslant &
\left| 
\mathbb{E}
\left[ 
\left.
\frac{ (Y(x_0) - \widehat{Y}(x_0))^2 }{ \widehat{\sigma}(x_0)^2 }
\mathds{1}_{ \{Y \not \in \mathcal{E}_{\kappa}\} }
\right| 
y
\right]
\right|
\\
&
+
\left|
\mathbb{E}
\left[
\left.
\mathds{1}_{\{ Y \not \in \mathcal{E}_{\kappa}\} }
\frac{\widehat{Y}(x_0) - \widehat{Y}_c(x_0)}{ \widehat{\sigma}(x_0) }
\left(
2 \frac{Y(x_0) - \widehat{Y}(x_0)}{ \widehat{\sigma}(x_0) }
+
\frac{\widehat{Y}(x_0) - \widehat{Y}_c(x_0)}{ \widehat{\sigma}(x_0) }
\right)
\right|
y
\right]
\right|
\\
\leqslant &
\mathbb{P}( Y \not \in \mathcal{E}_\kappa | y )^{1/2}
\mathbb{E}
\left[ 
\left.
\frac{ (Y(x_0) - \widehat{Y}(x_0))^4 }{ \widehat{\sigma}(x_0)^4 }
\right| 
y
\right]^{1/2}+
\left(
\frac{\widehat{Y}(x_0) - \widehat{Y}_c(x_0)}{ \widehat{\sigma}(x_0) }
\right)^2
\mathbb{P}(Y \not \in \mathcal{E}_{\kappa} | y)
\\
& + 2
\left|
\frac{\widehat{Y}(x_0) - \widehat{Y}_c(x_0)}{ \widehat{\sigma}(x_0) }
\right|
\mathbb{E}
\left[
\left|
\frac{Y(x_0) - \widehat{Y}(x_0)}{ \widehat{\sigma}(x_0) }
\right|
y
\right]
\\
= & o_{\mathbb{P} | Y \in \mathcal{E}_\kappa}(1),
\end{align*}
since  $\mathbb{P}( Y \not \in \mathcal{E}_\kappa | y )= o_{\P | Y \in \mathcal{E}_{\kappa}}(1)$, using \eqref{eq:pred:means} and the fact that  conditionally to the observation vector $y$, $(Y(x_0) - \widehat{Y}(x_0))/\widehat{\sigma}(x_0)$ is distributed as a standard Gaussian random variable. 
As $|C| = o_{\P | Y \in \mathcal{E}_{\kappa}}(1)$,  it remains to show that 
$|H|
=
O_{\mathbb{P} | Y \in \mathcal{E}_\kappa}(1)$
which is done by using that 
\[
|H|
\leqslant
\mathbb{E}
\left[
\left.
\frac{ (Y(x_0) - \widehat{Y}(x_0))^2 }{ \widehat{\sigma}(x_0)^2 }
\right| 
y
\right]
+ |D|=1+|D|,
\]
and that $D =  o_{ \P | Y \in \mathcal{E}_{\kappa} } (1)$ as established above.
\hfill $\square$ \\

\tbf{Proof of Corollary \ref{corollary:equivalence:pred:constraints}}.
It is shown in the point \tbf{3)} in the proof of Theorem \ref{th:cond_before_matern_i} that $\mathbb{P}_1( Y \not \in \mathcal{E}_{\kappa} | y )   = o_{\P | Y \in \mathcal{E}_{\kappa}}(1)$, where the conditional probability $\mathbb{P}_1( \cdot | \cdot )$ is calculated with respect to $\widetilde{k}_1$. Hence, using \eqref{eq:equivalence:prediction:pred} and \eqref{eq:equivalence:prediction:var}, one can show that Proposition \ref{proposition:prediction} remains true when $\widehat{\sigma}(x_0)$, $\widehat{Y}(x_0)$, $\widehat{Y}_c(x_0)$, $\widehat{\sigma}_c(x_0)$, and $\mathbb{E}[ \cdot | \cdot ]$ are replaced by $\widehat{\sigma}_1(x_0)$, $\widehat{Y}_1(x_0)$, $\widehat{Y}_{c,1}(x_0)$, $\widehat{\sigma}_{c,1}(x_0)$, and $\mathbb{E}_1[ \cdot | \cdot ]$. Then, the corollary is a consequence of this updated Proposition \ref{proposition:prediction} and of \eqref{eq:equivalence:prediction:pred} and \eqref{eq:equivalence:prediction:var}.
\hfill $\square$ \\

\tbf{Proof of Theorems \ref{th:cond_after_wendland}, \ref{th:cond_before_wendland_i} and \ref{th:cond_before_wendland_ii}}.
The proof is the same as in the Mat\'ern case in Theorems \ref{th:cond_after_matern}, \ref{th:cond_before_matern_i} and \ref{th:cond_before_matern_ii}. In particular, when $1+2s = \nu$, the Mat\'ern and Wendland covariance functions have the same smoothness, see  \citep[Theorem 1]{bevilacqua2019estimation}. Hence, a lemma similar as Lemma \ref{lem:lower:bounded:proba:bounded:matern} holds. We also remark that a lemma similar as Lemma \ref{lem:mater:change:rho} can be proved, by using \citep[Lemma 1]{bevilacqua2019estimation} together with the results given in the proof of  \citep[Theorem 8]{bevilacqua2019estimation} (see the online supplementary material to this paper).
\hfill $\square$ \\

\tbf{Proof of Proposition \ref{proposition:nugget}}.
Without loss of generality, we can consider that $u < + \infty$. Recall the notation $Y_c^{**} = \sup_{x \in [0,1]^d} | Y_c(x) | < + \infty$ a.s. of Appendix \ref{sec:intermediate}. Let $(x_i)_{i \in \mathbb{N}}$ be any sequence of two-by-two distinct points in $[0,1]^d$. We have
\begin{align*}
\mathbb{P}
\left(
\forall x \in [0,1]^d,
\;
\ell \leqslant Y(x) \leqslant u
\right)
&
\leqslant
\mathbb{P} 
\left(
\forall i \in \{1,\ldots,n\},
\;
Y_c(x_i) + Y_{\delta}(x_i) \leqslant u
\right) 
\\
& =
\mathbb{E} \left[
\mathbb{P}
\left(
\left.
\forall i \in \{1,\ldots,n\}, \;
Y_c(x_i) + Y_{\delta}(x_i) \leqslant u
\right|
Y_c(x_1) , \ldots , Y_c(x_n)
\right)
\right]
\\
& \leqslant
\mathbb{E} \left[
\mathbb{P}
\left(
\left.
\forall i \in \{1,\ldots,n\}, \;
Y_{\delta}(x_i) \leqslant |u| + Y_c^{**}
\right|
Y_c(x_1) , \ldots , Y_c(x_n)
\right)
\right]
\\
& =
\mathbb{E} \left[
\mathbb{P}
\left(
\left.
Y_{\delta}(x_1) \leqslant |u| + Y_c^{**}
\right|
Y_c^{**}
\right)^n
\right].
\end{align*}
The above probability goes to zero as $n \to \infty$ for any $Y_c^{**} < \infty$. Thus by dominated convergence, the above expectation goes to zero as $n \to \infty$. This concludes the proof. \hfill $\square$ \\

\tbf{Proof of Proposition \ref{prop:MLE:noisy}}.
Let $X_n \defeq n^{1/4}
\left(
\widehat \sigma_{n}^2/\widehat \rho_n
-
\sigma_0^2 / \rho_0
\right)$. For $K \geqslant 0$, we have
\begin{align*}
\limsup_{n \to \infty}
\mathbb{P} \left(
|X_n| \geqslant K
|
Y \in \mathcal{E}_0
\right)
&
\leqslant 
\limsup_{n \to \infty}
\frac{
	\mathbb{P} \left(
	|X_n| \geqslant K
	\right)
}
{
	\mathbb{P} \left(
	Y \in \mathcal{E}_0
	\right)
}
\underset{K \to \infty}{\longrightarrow} 0,
\end{align*}
from \eqref{eq:TCL:noisy} and because $\mathbb{P} \left(
Y \in \mathcal{E}_0
\right) >0$ does not depend on $K$ from Lemma \ref{lem:lower:bounded:proba:bounded:GP_var}. Hence $X_n = O_{\P | Y \in \mathcal{E}_0}(1)$. Moreover, let $\epsilon >0$. Then, 
\begin{align*}
\liminf_{n \to \infty}
\mathbb{P} \left(
|X_n| \geqslant \epsilon
|
Y \in \mathcal{E}_0
\right)
\geqslant 
\frac{  
	1 -
	\displaystyle \limsup_{n \to \infty}
	\mathbb{P} \left(
	|X_n| < \epsilon
	\right)
	- 
	\mathbb{P} \left(
	Y \not \in \mathcal{E}_0
	\right)
}
{
	\mathbb{P} \left(
	Y \in \mathcal{E}_0
	\right)
}
 \underset{\epsilon \to 0}{\longrightarrow} 
\frac{  
	1 -
	\mathbb{P} \left(
	Y \not \in \mathcal{E}_0
	\right)
}
{
	\mathbb{P} \left(
	Y \in \mathcal{E}_0
	\right)
}
 = 1
,
\end{align*}
since the asymptotic variance is non-zero in \eqref{eq:TCL:noisy}. Hence $X_n \not = o_{\P | Y \in \mathcal{E}_0}(1)$. By the same arguments, the same conclusions hold for the MLE of the noise variance. \hfill $\square$ \\

\bibliographystyle{apa}  
\bibliography{arXiv2019_cMLE}  

\end{document}